%% file: univinj.tex
\documentclass{amsart}
\newif\ifcref\creftrue
\let\sect\S
\input{decls}
\usepackage{newclude}
\usepackage[all]{xy}
\fxsetup{final}
\usetikzlibrary{decorations.pathmorphing}
\title{All $(\infty,1)$-toposes have strict univalent universes}
\author{Michael Shulman}
\thanks{This material is based
    on research sponsored by The United States Air Force Research
    Laboratory under agreement number FA9550-15-1-0053.  The
    U.S. Government is authorized to reproduce and distribute reprints
    for Governmental purposes notwithstanding any copyright notation
    thereon.  The views and conclusions contained herein are those of
    the author and should not be interpreted as necessarily
    representing the official policies or endorsements, either
    expressed or implied, of the United States Air Force Research
    Laboratory, the U.S. Government, or Carnegie Mellon University.}
\date{\today}

\def\type{\,\mathsf{type}}
\autodefs{\cGpd\cCat\nSet\ncore\cGPD\cCAT\cPsFun\cPsNat\nTop\nSET\cAdj\cRetr\dFib\dRep\dTy\dTm\nTy\nTm}

\def\iLoc{\mathsf{Loc}}
\def\pr#1#2{\pmb[{#1}\op,#2\pmb]}
\def\func#1#2{\pmb[{#1},#2\pmb]}
\let\V\sV
\let\W\sW
\let\D\sD
\let\E\sE
\let\M\sM
\let\N\sN
\let\P\sP
\let\S\sS
\let\C\sC
\let\E\sE

\def\fact#1#2{\phi_{#1,#2}}
\def\io{\ensuremath{(\infty,1)}}

\def\Ehat{\ensuremath{\mathcal{PSH}(\E)}\xspace}
\let\F\dF
\def\Fib{\dFib}
\let\cE\dE
\def\cEka{{\ensuremath{\dE^\ka}}\xspace}
\def\Fka{{\ensuremath{\dF^\ka}}\xspace}
\def\Fibka{{\ensuremath{\Fib^\ka}}\xspace}
\def\cEp{{\ensuremath{\dE_{\scriptscriptstyle\bullet}}}\xspace}

\DeclareFontFamily{U}{min}{}
\DeclareFontShape{U}{min}{m}{n}{<-> udmj30}{}

\def\prcs{\pr\C\nSet}
\def\unit{\mathbb{I}}
\def\ddup{\dDelta^{\mathrm{op}\uparrow}}

\def\phihat{\widehat{\phi}}

\def\drr{\dDelta\rtimes\cR}

\def\slice#1#2{#1/_{#2}}
\def\pscl#1{#1'}
\def\cocl#1{#1^\star}
\def\nfs{notion of fibred structure\xspace}

\def\nfss{notions of fibred structure\xspace}

\def\uly#1{|#1|}
\def\ttmt{type-theoretic model topos\xspace}
\def\ttmts{type-theoretic model toposes\xspace}
\def\Ttmts{Type-theoretic model toposes\xspace}

\def\ehom#1{\underline{\smash[b]{#1}}}
\let\cpw\odot
\let\pow\pitchfork

\def\lcirc{\mathbin{\widehat{\circ}}}
\makeatletter\def\app{\mathbin{@}}\makeatother
\def\lapp{\mathbin{\widehat{\app}}}

\def\lifts{\boxslash}
\def\bitimes{\times^{\mathrm{h}}}
\def\bilim{\operatorname{lim}^{\mathrm{h}}}
\def\bicolim{\colim^{\mathrm{h}}}

\let\cohto\rightsquigarrow

\def\zero{\underline{0}}
\def\one{\underline{1}}
\def\bar{\mathsf{B}}
\def\barp{\mathsf{B}'}
\def\sbar{\mathsf{B}_\bullet}
\def\sbarp{\mathsf{B}'_\bullet}
\def\cobar{\mathsf{C}}
\def\scobar{\mathsf{C}^\bullet}

\def\local{locally representable\xspace}

\def\locality{local representability\xspace}

\def\stratified{relatively acyclic\xspace}

\def\ce{\mathord{\centerdot}}
\def\ec{\diamond}
\def\las#1{#1_{\pmb{!}}}
\def\lasv#1{#1_{\pmb{!},\cV}}

\def\lasU#1{\iU_{\pmb{!},#1}}
\def\lasu#1{\iu_{\pmb{!},#1}}
\def\lasEl#1{\mathsf{El}_{\pmb{!},#1}}
\def\lu{\triangleleft}
\def\suc{\ensuremath{\mathsf{suc}}\xspace}
\def\Lift{\ensuremath{\mathsf{Lift}}\xspace}
\def\lift{\ensuremath{\mathsf{lift}}\xspace}
\def\low{\ensuremath{\mathsf{lower}}\xspace}
\def\U{\iU}
\def\El{\mathsf{El}}
\let\jdeq\equiv
\def\iid{\mathsf{id}}

\def\isequiv{\mathsf{isEquiv}}
\def\equiv{\mathsf{Equiv}}
\def\Eq{\mathsf{Eq}}
\def\Idtype{\mathsf{Id}}
\def\Prop{\mathsf{Prop}}
\def\Iso{\mathsf{Iso}}
\def\Retr{\mathsf{Retr}}
\def\bool{\ensuremath{\mathbf{2}}\xspace}

\def\slcc{simplicially locally cartesian closed\xspace}
\def\slcclosure{simplicial local cartesian closure\xspace}
\def\Slcclosure{Simplicial local cartesian closure\xspace}

\def\qcof{Quillen cofibration\xspace}
\def\qcofs{Quillen cofibrations\xspace}
\def\qfib{Quillen fibration\xspace}

\def\qcoft{Quillen cofibrant\xspace}
\def\qfibt{Quillen fibrant\xspace}
\def\qucoft{Quillen unit-cofibrant\xspace}
\def\qufibt{Quillen counit-fibrant\xspace}

\let\shlt\lhd

\let\shgt\rhd
\def\nshlt{\mathrel{\mathrlap{\,\not}\shlt}}

\let\diag\triangle
\def\dia{^\triangle}
\def\pb{^{\mathrm{pb}}}

\def\card#1{|#1|}

\usepackage{version}
\includeversion{concise}
\excludeversion{verbose}

\begin{document}
\maketitle

\begin{abstract}
  We prove the conjecture that any Grothendieck \io-topos can be presented by a Quillen model category that interprets homotopy type theory with strict univalent universes.
  Thus, homotopy type theory can be used as a formal language for reasoning internally to \io-toposes, just as higher-order logic is used for 1-toposes.

  As part of the proof, we give a new, more explicit, characterization of the fibrations in injective model structures on presheaf categories.
  In particular, we show that they generalize the coflexible algebras of 2-monad theory.
\end{abstract}

\DeclareRobustCommand{\gobblefive}[5]{}
\newcommand*{\skiptoc}{\addtocontents{toc}{\gobblefive}}
\setcounter{tocdepth}{1}
\tableofcontents


\include*{intro}

\include*{2cat}

\include*{nfs}
\include*{relpres}
\include*{universes}
\include*{ttmt}

\include*{bar}
\include*{injmodel}

\include*{intloc}
\include*{lexloc}

\include*{summary}

\appendix
\include*{coherence}

\bibliographystyle{alpha}
\bibliography{all}
\end{document}

%% file: decls.tex
\usepackage{amssymb,amsmath,stmaryrd,mathrsfs}

\def\definetac{\newif\iftac}    
\ifx\tactrue\undefined
  \definetac
  \ifx\state\undefined\tacfalse\else\tactrue\fi
\fi

\def\definebeamer{\newif\ifbeamer}
\ifx\beamertrue\undefined
  \definebeamer
  \ifx\uncover\undefined\beamerfalse\else\beamertrue\fi
\fi

\def\definecref{\newif\ifcref}
\ifx\creftrue\undefined
  \definecref
  \creffalse
\fi

\iftac\else\usepackage{amsthm}\fi
\usepackage{tikz,tikz-cd}
\usetikzlibrary{arrows}
\ifbeamer\else
  \usepackage{enumitem}
  \usepackage{xcolor}
  \definecolor{darkgreen}{rgb}{0,0.45,0} 
  \iftac\else\usepackage[colorlinks,citecolor=darkgreen,linkcolor=darkgreen,backref=page]{hyperref}
  \renewcommand*{\backref}[1]{}
  \renewcommand*{\backrefalt}[4]{({%
      \ifcase #1 Not cited.%
            \or cited on p.~#2%
            \else cited on pp.~#2%
      \fi%
    })}
  \fi
\fi
\usepackage{mathtools}          
\usepackage{ifmtarg}            
\usepackage{microtype}
\usepackage{braket}             
\let\setof\Set
\usepackage{url}                
\usepackage{xspace}             
\ifcref\usepackage{cleveref,aliascnt}\fi
\usepackage[status=draft,author=]{fixme}
\fxusetheme{color}
\usepackage{mathpartir}

\makeatletter
\let\ea\expandafter

\def\mdef#1#2{\ea\ea\ea\gdef\ea\ea\noexpand#1\ea{\ea\ensuremath\ea{#2}\xspace}}
\def\alwaysmath#1{\ea\ea\ea\global\ea\ea\ea\let\ea\ea\csname your@#1\endcsname\csname #1\endcsname
  \ea\def\csname #1\endcsname{\ensuremath{\csname your@#1\endcsname}\xspace}}

\newcount\foreachcount

\def\foreachletter#1#2#3{\foreachcount=#1
  \ea\loop\ea\ea\ea#3\@alph\foreachcount
  \advance\foreachcount by 1
  \ifnum\foreachcount<#2\repeat}

\def\foreachLetter#1#2#3{\foreachcount=#1
  \ea\loop\ea\ea\ea#3\@Alph\foreachcount
  \advance\foreachcount by 1
  \ifnum\foreachcount<#2\repeat}

\def\definescr#1{\ea\gdef\csname s#1\endcsname{\ensuremath{\mathscr{#1}}\xspace}}
\foreachLetter{1}{27}{\definescr}
\def\definecal#1{\ea\gdef\csname c#1\endcsname{\ensuremath{\mathcal{#1}}\xspace}}
\foreachLetter{1}{27}{\definecal}
\def\definebold#1{\ea\gdef\csname b#1\endcsname{\ensuremath{\mathbf{#1}}\xspace}}
\foreachLetter{1}{27}{\definebold}
\def\definebb#1{\ea\gdef\csname d#1\endcsname{\ensuremath{\mathbb{#1}}\xspace}}
\foreachLetter{1}{27}{\definebb}
\def\definefrak#1{\ea\gdef\csname f#1\endcsname{\ensuremath{\mathfrak{#1}}\xspace}}
\foreachletter{1}{9}{\definefrak}
\foreachletter{10}{27}{\definefrak}
\foreachLetter{1}{27}{\definefrak}
\def\definesf#1{\ea\gdef\csname i#1\endcsname{\ensuremath{\mathsf{#1}}\xspace}}
\foreachletter{1}{6}{\definesf}
\foreachletter{7}{14}{\definesf}
\foreachletter{15}{20}{\definesf}
\foreachletter{21}{27}{\definesf}
\foreachLetter{1}{27}{\definesf}
\def\definebar#1{\ea\gdef\csname #1bar\endcsname{\ensuremath{\overline{#1}}\xspace}}
\foreachLetter{1}{27}{\definebar}
\foreachletter{1}{8}{\definebar} 
\foreachletter{9}{15}{\definebar} 
\foreachletter{16}{27}{\definebar}
\def\definetil#1{\ea\gdef\csname #1til\endcsname{\ensuremath{\widetilde{#1}}\xspace}}
\foreachLetter{1}{27}{\definetil}
\foreachletter{1}{27}{\definetil}
\def\definehat#1{\ea\gdef\csname #1hat\endcsname{\ensuremath{\widehat{#1}}\xspace}}
\foreachLetter{1}{27}{\definehat}
\foreachletter{1}{27}{\definehat}
\def\definechk#1{\ea\gdef\csname #1chk\endcsname{\ensuremath{\widecheck{#1}}\xspace}}
\foreachLetter{1}{27}{\definechk}
\foreachletter{1}{27}{\definechk}
\def\defineul#1{\ea\gdef\csname u#1\endcsname{\ensuremath{\underline{#1}}\xspace}}
\foreachLetter{1}{27}{\defineul}
\foreachletter{1}{27}{\defineul}

\def\autofmt@n#1\autofmt@end{\mathrm{#1}}
\def\autofmt@b#1\autofmt@end{\mathbf{#1}}
\def\autofmt@d#1#2\autofmt@end{\mathbb{#1}\mathsf{#2}}
\def\autofmt@c#1#2\autofmt@end{\mathcal{#1}\mathit{#2}}
\def\autofmt@s#1#2\autofmt@end{\mathscr{#1}\mathit{#2}}
\def\autofmt@f#1\autofmt@end{\mathsf{#1}}
\def\autofmt@k#1\autofmt@end{\mathfrak{#1}}
\def\autofmt@u#1\autofmt@end{\underline{\smash{\mathsf{#1}}}}
\def\autofmt@U#1\autofmt@end{\underline{\underline{\smash{\mathsf{#1}}}}}
\def\autofmt@h#1\autofmt@end{\widehat{#1}}
\def\autofmt@r#1\autofmt@end{\overline{#1}}
\def\autofmt@t#1\autofmt@end{\widetilde{#1}}
\def\autofmt@k#1\autofmt@end{\check{#1}}

\def\auto@drop#1{}
\def\autodef#1{\ea\ea\ea\@autodef\ea\ea\ea#1\ea\auto@drop\string#1\autodef@end}
\def\@autodef#1#2#3\autodef@end{%
  \ea\def\ea#1\ea{\ea\ensuremath\ea{\csname autofmt@#2\endcsname#3\autofmt@end}\xspace}}
\def\autodefs@end{blarg!}
\def\autodefs#1{\@autodefs#1\autodefs@end}
\def\@autodefs#1{\ifx#1\autodefs@end%
  \def\autodefs@next{}%
  \else%
  \def\autodefs@next{\autodef#1\@autodefs}%
  \fi\autodefs@next}

\DeclareSymbolFont{bbold}{U}{bbold}{m}{n}
\DeclareSymbolFontAlphabet{\mathbbb}{bbold}
\newcommand{\dDelta}{\ensuremath{\mathbbb{\Delta}}\xspace}

\newcommand{\dtwo}{\ensuremath{\mathbbb{2}}\xspace}
\newcommand{\dthree}{\ensuremath{\mathbbb{3}}\xspace}

\mdef\delbar{\overline{\partial}}

\mdef\hf{\textstyle\frac12 }
\mdef\thrd{\textstyle\frac13 }
\mdef\qtr{\textstyle\frac14 }

\let\dn\downarrow
\newcommand{\op}{^{\mathrm{op}}}

\let\adj\dashv

\mdef\Id{\mathrm{Id}}
\mdef\id{\mathrm{id}}
\alwaysmath{ell}
\alwaysmath{infty}
\let\oo\infty
\alwaysmath{odot}
\def\frc#1/#2.{\frac{#1}{#2}}   
\mdef\ten{\mathrel{\otimes}}

\mdef\sqten{\mathrel{\boxtimes}}
\def\gt{>}

\DeclareFontFamily{U}{mathb}{\hyphenchar\font45}
\DeclareFontShape{U}{mathb}{m}{n}{
      <5> <6> <7> <8> <9> <10> gen * mathb
      <10.95> mathb10 <12> <14.4> <17.28> <20.74> <24.88> mathb12
      }{}
\DeclareSymbolFont{mathb}{U}{mathb}{m}{n}
\DeclareFontSubstitution{U}{mathb}{m}{n}
\DeclareFontFamily{U}{mathx}{\hyphenchar\font45}
\DeclareFontShape{U}{mathx}{m}{n}{
      <5> <6> <7> <8> <9> <10>
      <10.95> <12> <14.4> <17.28> <20.74> <24.88>
      mathx10
      }{}
\DeclareSymbolFont{mathx}{U}{mathx}{m}{n}
\DeclareFontSubstitution{U}{mathx}{m}{n}
\DeclareMathSymbol{\dotplus}       {2}{mathb}{"00}
\DeclareMathSymbol{\dotdiv}        {2}{mathb}{"01}
\DeclareMathSymbol{\dottimes}      {2}{mathb}{"02}
\DeclareMathSymbol{\divdot}        {2}{mathb}{"03}
\DeclareMathSymbol{\udot}          {2}{mathb}{"04}
\DeclareMathSymbol{\square}        {2}{mathb}{"05}
\DeclareMathSymbol{\Asterisk}      {2}{mathb}{"06}
\DeclareMathSymbol{\bigast}        {1}{mathb}{"06}
\DeclareMathSymbol{\coAsterisk}    {2}{mathb}{"07}
\DeclareMathSymbol{\bigcoast}      {1}{mathb}{"07}
\DeclareMathSymbol{\circplus}      {2}{mathb}{"08}
\DeclareMathSymbol{\pluscirc}      {2}{mathb}{"09}
\DeclareMathSymbol{\convolution}   {2}{mathb}{"0A}
\DeclareMathSymbol{\divideontimes} {2}{mathb}{"0B}
\DeclareMathSymbol{\blackdiamond}  {2}{mathb}{"0C}
\DeclareMathSymbol{\sqbullet}      {2}{mathb}{"0D}
\DeclareMathSymbol{\bigstar}       {2}{mathb}{"0E}
\DeclareMathSymbol{\bigvarstar}    {2}{mathb}{"0F}
\DeclareMathSymbol{\corresponds}   {3}{mathb}{"1D}
\DeclareMathSymbol{\boxleft}      {2}{mathb}{"68}
\DeclareMathSymbol{\boxright}     {2}{mathb}{"69}
\DeclareMathSymbol{\boxtop}       {2}{mathb}{"6A}
\DeclareMathSymbol{\boxbot}       {2}{mathb}{"6B}
\DeclareMathSymbol{\updownarrows}          {3}{mathb}{"D6}
\DeclareMathSymbol{\downuparrows}          {3}{mathb}{"D7}
\DeclareMathSymbol{\Lsh}                   {3}{mathb}{"E8}
\DeclareMathSymbol{\Rsh}                   {3}{mathb}{"E9}
\DeclareMathSymbol{\dlsh}                  {3}{mathb}{"EA}
\DeclareMathSymbol{\drsh}                  {3}{mathb}{"EB}
\DeclareMathSymbol{\looparrowdownleft}     {3}{mathb}{"EE}
\DeclareMathSymbol{\looparrowdownright}    {3}{mathb}{"EF}
\DeclareMathSymbol{\curvearrowleftright}   {3}{mathb}{"F2}
\DeclareMathSymbol{\curvearrowbotleft}     {3}{mathb}{"F3}
\DeclareMathSymbol{\curvearrowbotright}    {3}{mathb}{"F4}
\DeclareMathSymbol{\curvearrowbotleftright}{3}{mathb}{"F5}
\DeclareMathSymbol{\leftsquigarrow}        {3}{mathb}{"F8}
\DeclareMathSymbol{\rightsquigarrow}       {3}{mathb}{"F9}
\DeclareMathSymbol{\leftrightsquigarrow}   {3}{mathb}{"FA}
\DeclareMathSymbol{\lefttorightarrow}      {3}{mathb}{"FC}
\DeclareMathSymbol{\righttoleftarrow}      {3}{mathb}{"FD}
\DeclareMathSymbol{\uptodownarrow}         {3}{mathb}{"FE}
\DeclareMathSymbol{\downtouparrow}         {3}{mathb}{"FF}
\DeclareMathSymbol{\varhash}       {0}{mathb}{"23}
\DeclareMathSymbol{\sqSubset}       {3}{mathb}{"94}
\DeclareMathSymbol{\sqSupset}       {3}{mathb}{"95}
\DeclareMathSymbol{\nsqSubset}      {3}{mathb}{"96}
\DeclareMathSymbol{\nsqSupset}      {3}{mathb}{"97}
\DeclareMathAccent{\widecheck}    {0}{mathx}{"71}

\DeclareMathOperator\colim{colim}

\DeclareMathOperator*\hocolim{hocolim}
\DeclareMathOperator*\holim{holim}

\DeclareMathOperator\ob{ob}

\DeclareMathOperator\Ho{Ho}

\newcommand{\too}[1][]{\ensuremath{\overset{#1}{\longrightarrow}}}
\newcommand{\ot}{\ensuremath{\leftarrow}}

\let\toot\rightleftarrows

\let\toto\rightrightarrows
\let\into\hookrightarrow

\mdef\we{\overset{\sim}{\longrightarrow}}
\mdef\leftwe{\overset{\sim}{\longleftarrow}}
\let\mono\rightarrowtail

\let\cof\rightarrowtail

\let\fib\twoheadrightarrow

\let\cohto\rightsquigarrow

\def\acof{\mathrel{\mathrlap{\hspace{3pt}\raisebox{4pt}{$\scriptscriptstyle\sim$}}\mathord{\rightarrowtail}}}

\let\xto\xrightarrow
\let\xot\xleftarrow
\def\rightarrowtailfill@{\arrowfill@{\Yright\joinrel\relbar}\relbar\rightarrow}
\newcommand\xrightarrowtail[2][]{\ext@arrow 0055{\rightarrowtailfill@}{#1}{#2}}

\def\twoheadrightarrowfill@{\arrowfill@{\relbar\joinrel\relbar}\relbar\twoheadrightarrow}
\newcommand\xtwoheadrightarrow[2][]{\ext@arrow 0055{\twoheadrightarrowfill@}{#1}{#2}}

\def\slashedarrowfill@#1#2#3#4#5{%
  $\m@th\thickmuskip0mu\medmuskip\thickmuskip\thinmuskip\thickmuskip
   \relax#5#1\mkern-7mu%
   \cleaders\hbox{$#5\mkern-2mu#2\mkern-2mu$}\hfill
   \mathclap{#3}\mathclap{#2}%
   \cleaders\hbox{$#5\mkern-2mu#2\mkern-2mu$}\hfill
   \mkern-7mu#4$%
}
\def\rightslashedarrowfill@{%
  \slashedarrowfill@\relbar\relbar\mapstochar\rightarrow}
\newcommand\xslashedrightarrow[2][]{%
  \ext@arrow 0055{\rightslashedarrowfill@}{#1}{#2}}
\mdef\hto{\xslashedrightarrow{}}
\mdef\htoo{\xslashedrightarrow{\quad}}

\def\toiso{\xto{\smash{\raisebox{-.5mm}{$\scriptstyle\sim$}}}}

\def\jd#1{\@jd#1\ej}
\def\@jd#1|-#2\ej{\@@jd#1,,\;\vdash\;\left(#2\right)}
\def\@@jd#1,{\@ifmtarg{#1}{\let\next=\relax}{\left(#1\right)\let\next=\@@@jd}\next}
\def\@@@jd#1,{\@ifmtarg{#1}{\let\next=\relax}{,\,\left(#1\right)\let\next=\@@@jd}\next}
\def\jdm#1{\@jdm#1\ej}
\def\@jdm#1|-#2\ej{\@@jd#1,,\\\vdash\;\left(#2\right)}

\def\alg{\text{-}\mathcal{A}\mathit{lg}}
\def\algs{\text{-}\mathcal{A}\mathit{lg}_s}

\long\def\my@drawfill#1#2;{%
\@skipfalse
\fill[#1,draw=none] #2;
\@skiptrue
\draw[#1,fill=none] #2;
}
\newif\if@skip
\newcommand{\skipit}[1]{\if@skip\else#1\fi}
\newcommand{\drawfill}[1][]{\my@drawfill{#1}}

\newcounter{nodemaker}
\newcommand{\drpullback}[1][dr]{\ar[#1,phantom,near start,"\lrcorner"]}
\newcommand{\drpushout}[1][dr]{\ar[#1,phantom,near end,"\ulcorner"]}

\newif\ifhyperref
\@ifpackageloaded{hyperref}{\hyperreftrue}{\hyperreffalse}
\iftac
  \let\your@state\state
  \def\state#1{\my@state#1}
  \def\my@state#1.{\gdef\currthmtype{#1}\your@state{#1.}}
  \let\your@staterm\staterm
  \def\staterm#1{\my@staterm#1}
  \def\my@staterm#1.{\gdef\currthmtype{#1}\your@staterm{#1.}}
  \let\@defthm\newtheorem
  \def\switchtotheoremrm{\let\@defthm\newtheoremrm}
  \def\defthm#1#2#3{\@defthm{#1}{#2}} 
  \let\your@section\section
  \def\section{\gdef\currthmtype{section}\your@section}
  \def\currthmtype{}
  \ifhyperref
    \def\autoref#1{\ref*{label@name@#1}~\ref{#1}}
  \else
    \def\autoref#1{\ref{label@name@#1}~\ref{#1}}
  \fi
  \AtBeginDocument{%
    \let\old@label\label%
    \def\label#1{%
      {\let\your@currentlabel\@currentlabel%
        \edef\@currentlabel{\currthmtype}%
        \old@label{label@name@#1}}%
      \old@label{#1}}
  }
  \let\cref\autoref
\else\ifcref
  \def\defthm#1#2#3{%
    \newaliascnt{#1}{thm}
    \newtheorem{#1}[#1]{#2}
    \aliascntresetthe{#1}
    \crefname{#1}{#2}{#3}
  }
\else
  \ifhyperref
    \def\defthm#1#2#3{
      \newtheorem{#1}{#2}[section]%
      \expandafter\def\csname #1autorefname\endcsname{#2}%
      \expandafter\let\csname c@#1\endcsname\c@thm}
  \else
    \def\defthm#1#2#3{\newtheorem{#1}[thm]{#2}} 
    \ifx\SK@label\undefined\let\SK@label\label\fi
    \let\old@label\label
    \let\your@thm\@thm
    \def\@thm#1#2#3{\gdef\currthmtype{#3}\your@thm{#1}{#2}{#3}}
    \def\currthmtype{}
    \def\label#1{{\let\your@currentlabel\@currentlabel\def\@currentlabel%
        {\currthmtype~\your@currentlabel}%
        \SK@label{#1@}}\old@label{#1}}
    \def\autoref#1{\ref{#1@}}
  \fi
  \let\cref\autoref
\fi\fi

\newtheorem{thm}{Theorem}[section]
\ifcref
  \crefname{thm}{Theorem}{Theorems}
\else
  
\fi
\defthm{cor}{Corollary}{Corollaries}
\defthm{prop}{Proposition}{Propositions}
\defthm{lem}{Lemma}{Lemmas}
\defthm{sch}{Scholium}{Scholia}
\defthm{assume}{Assumption}{Assumptions}
\defthm{claim}{Claim}{Claims}
\defthm{conj}{Conjecture}{Conjectures}
\defthm{hyp}{Hypothesis}{Hypotheses}
\iftac\switchtotheoremrm\else\theoremstyle{definition}\fi
\defthm{defn}{Definition}{Definitions}
\defthm{notn}{Notation}{Notations}
\defthm{rmk}{Remark}{Remarks}
\defthm{eg}{Example}{Examples}
\defthm{egs}{Examples}{Examples}
\defthm{ex}{Exercise}{Exercises}
\defthm{ceg}{Counterexample}{Counterexamples}

\ifcref
  \crefformat{section}{\sect#2#1#3}
  \Crefformat{section}{Section~#2#1#3}
  \crefrangeformat{section}{\sect\sect#3#1#4--#5#2#6}
  \Crefrangeformat{section}{Sections~#3#1#4--#5#2#6}
  \crefmultiformat{section}{\sect\sect#2#1#3}{ and~#2#1#3}{, #2#1#3}{ and~#2#1#3}
  \Crefmultiformat{section}{Sections~#2#1#3}{ and~#2#1#3}{, #2#1#3}{ and~#2#1#3}
  \crefrangemultiformat{section}{\sect\sect#3#1#4--#5#2#6}{ and~#3#1#4--#5#2#6}{, #3#1#4--#5#2#6}{ and~#3#1#4--#5#2#6}
  \Crefrangemultiformat{section}{Sections~#3#1#4--#5#2#6}{ and~#3#1#4--#5#2#6}{, #3#1#4--#5#2#6}{ and~#3#1#4--#5#2#6}
  \crefformat{appendix}{Appendix~#2#1#3}
  \Crefformat{appendix}{Appendix~#2#1#3}
  \crefrangeformat{appendix}{Appendices~#3#1#4--#5#2#6}
  \Crefrangeformat{appendix}{Appendices~#3#1#4--#5#2#6}
  \crefmultiformat{appendix}{Appendices~#2#1#3}{ and~#2#1#3}{, #2#1#3}{ and~#2#1#3}
  \Crefmultiformat{appendix}{Appendices~#2#1#3}{ and~#2#1#3}{, #2#1#3}{ and~#2#1#3}
  \crefrangemultiformat{appendix}{Appendices~#3#1#4--#5#2#6}{ and~#3#1#4--#5#2#6}{, #3#1#4--#5#2#6}{ and~#3#1#4--#5#2#6}
  \Crefrangemultiformat{appendix}{Appendices~#3#1#4--#5#2#6}{ and~#3#1#4--#5#2#6}{, #3#1#4--#5#2#6}{ and~#3#1#4--#5#2#6}
  \crefformat{subappendix}{\sect#2#1#3}
  \Crefformat{subappendix}{Section~#2#1#3}
  \crefrangeformat{subappendix}{\sect\sect#3#1#4--#5#2#6}
  \Crefrangeformat{subappendix}{Sections~#3#1#4--#5#2#6}
  \crefmultiformat{subappendix}{\sect\sect#2#1#3}{ and~#2#1#3}{, #2#1#3}{ and~#2#1#3}
  \Crefmultiformat{subappendix}{Sections~#2#1#3}{ and~#2#1#3}{, #2#1#3}{ and~#2#1#3}
  \crefrangemultiformat{subappendix}{\sect\sect#3#1#4--#5#2#6}{ and~#3#1#4--#5#2#6}{, #3#1#4--#5#2#6}{ and~#3#1#4--#5#2#6}
  \Crefrangemultiformat{subappendix}{Sections~#3#1#4--#5#2#6}{ and~#3#1#4--#5#2#6}{, #3#1#4--#5#2#6}{ and~#3#1#4--#5#2#6}
  \crefname{part}{Part}{Parts}
  \crefname{figure}{Figure}{Figures}
\fi

\iftac
  \let\qed\endproof
  \let\your@endproof\endproof
  \def\my@endproof{\your@endproof}
  \def\endproof{\my@endproof\gdef\my@endproof{\your@endproof}}
  \def\qedhere{\tag*{\endproofbox}\gdef\my@endproof{\relax}}
\fi

\iftac
  \def\pr@@f[#1]{\subsubsection*{\sc #1.}}
\fi

\def\thmqedhere{\expandafter\csname\csname @currenvir\endcsname @qed\endcsname}

\ifbeamer\else

\fi

\ifbeamer\else
  \setitemize[1]{leftmargin=2em}
  \setenumerate[1]{leftmargin=*}
\fi

\iftac
  \let\c@equation\c@subsection
\else
  \let\c@equation\c@thm
\fi
\numberwithin{equation}{section}

\ifcref\else
  \@ifpackageloaded{mathtools}{\mathtoolsset{showonlyrefs,showmanualtags}}{}
\fi

\alwaysmath{alpha}
\alwaysmath{beta}
\alwaysmath{gamma}
\alwaysmath{Gamma}
\alwaysmath{delta}
\alwaysmath{Delta}
\alwaysmath{epsilon}
\mdef\ep{\varepsilon}
\alwaysmath{zeta}
\alwaysmath{eta}
\alwaysmath{theta}
\alwaysmath{Theta}
\alwaysmath{iota}
\alwaysmath{kappa}
\alwaysmath{lambda}
\alwaysmath{Lambda}
\alwaysmath{mu}
\alwaysmath{nu}
\alwaysmath{xi}
\alwaysmath{pi}
\alwaysmath{rho}
\alwaysmath{sigma}
\alwaysmath{Sigma}
\alwaysmath{tau}
\alwaysmath{upsilon}
\alwaysmath{Upsilon}
\alwaysmath{phi}
\alwaysmath{Pi}
\alwaysmath{Phi}
\mdef\ph{\varphi}
\alwaysmath{chi}
\alwaysmath{psi}
\alwaysmath{Psi}
\alwaysmath{omega}
\alwaysmath{Omega}
\let\al\alpha
\let\be\beta
\let\gm\gamma

\let\om\omega
\let\ka\kappa
\let\la\lambda

\let\ze\zeta

\makeatother


%% file: intro.tex
\section{Introduction}
\label{sec:introduction}

\subsection{Background}
\label{sec:background}

In the 1960s Grothendieck introduced \emph{toposes} (categories of sheaves on sites) as a powerful tool in the study of all kinds of geometry.
Part of their usefulness comes from the fact that they share most properties of the category of sets, so that all sorts of mathematics can be done ``internal'' to an arbitrary topos.
For instance, internal group theory yields the theory of sheaves of groups, and so on.

The direct way to express mathematics internally in a topos is to rewrite it manually in terms of objects, morphisms, and commutative diagrams; but this is tedious and verbose.
A more elegant method was found by Mitchell, Benabou, and Joyal in the early 1970s, building on the ``elementary'' toposes of Lawvere and Tierney.
Namely, the formal language of \emph{intuitionistic higher-order logic (IHOL)} can be interpreted algorithmically in any topos; thus almost any mathematics can be internalized in toposes simply by writing it in IHOL (or observing that it can be so written).
This generally requires very little modification,\footnote{Except that ``non-constructive'' principles such as the axiom of choice and the law of excluded middle must be avoided.} since the ``types'' of IHOL have ``elements'' and behave otherwise like sets.

In the 21\textsuperscript{st} century it has become clear that generalizing Grothendieck's toposes to \emph{higher categories} yields an even more powerful tool for geometry.
A \emph{Grothendieck \io-topos}~\cite{tv:hag-i,rezk:homotopy-toposes,lurie:higher-topoi} contains a subcategory that is an ordinary Grothendieck topos (a.k.a.\ 1-topos), but also many higher-dimensional objects that share the \io-categorical properties of the ``spaces'' in homotopy theory.
Thus not only set-based mathematics, but homotopy-theoretic and higher-categorical mathematics, can be done ``internal'' to any \io-topos.\footnote{One can also do homotopy-theoretic mathematics in a 1-topos, e.g.\ with internal simplicial objects or chain complexes, as indeed Grothendieck himself did.
  But compared to the \io-topos version, this approach sometimes fails to have the desired properties or to capture the desired information even about classical topological spaces; see~\cite[\sect 6.5.4]{lurie:higher-topoi} (the connection with internal simplicial objects is established by~\cite[Proposition 6.5.2.14]{lurie:higher-topoi}, \cite[Theorem 5]{jardine:fields-spre}, and~\cite{beke:sheafifiable,jardine:bool-loc}, while internal chain complexes are discussed in~\cite[Remark 6.5.4.3]{lurie:higher-topoi}).
  Moreover, there are important \io-toposes carrying homotopical information that is invisible to their underlying 1-topos, such as those whose objects are parametrized spaces~\cite{maysig:pht,abghr:thom-oocat}, parametrized spectra~\cite{joyal:logoi,hoyois:topoi-param}, generalized orbispaces~\cite{rezk:global-cohesion}, or excisive functors for Goodwillie calculus~\cite{abfj:goodwillie-htt}.}

Unsurprisingly, directly internalizing mathematics in an \io-topos is even more tedious and verbose than in a 1-topos, since there are many more coherences to keep track of.
Fortunately, an \io-categorical analogue of the Mitchell--Benabou--Joyal language emerged at about the same time as the theory of Grothendieck \io-toposes.
Building on~\cite{hs:gpd-typethy}, Awodey and Warren~\cite{aw:htpy-idtype} and Voevodsky (published a decade later as~\cite{klv:ssetmodel}) discovered that \emph{Martin-L\"{o}f intensional dependent type theory (MLTT)}~\cite{martinlof:itt,martinlof:itt-pred} has a interpretation in which its ``types'' behave like the spaces of homotopy theory, and hence like the objects of an \io-topos.
Voevodsky then formulated the \emph{univalence axiom} for MLTT, which internalizes the ``object classifiers'' of an \io-topos using a \emph{universe type}, just as the type of propositions in IHOL internalizes the subobject classifier of a 1-topos.

The study of formal systems with this sort of homotopical interpretation is now known as \emph{homotopy type theory (HoTT)}\footnote{Voevodsky's closely related ``univalent foundations'' program~\cite{vv:unimath} is, roughly, the use of such formal systems as a computer-verifiable foundation for mathematics.}.
The specific system of MLTT with univalence and also ``higher inductive types'', which internalize homotopy colimits and recursive constructions, is sometimes known as \emph{Book HoTT}\footnote{Although since~\cite{hottbook} did not formulate a general notion of ``higher inductive type'', the phrase ``Book HoTT'' does not have a fully precise meaning.} with reference to the book~\cite{hottbook} which popularized it.
In recent years it has been shown that large amounts of homotopy theory can be done inside Book HoTT, including homotopy groups of spheres~\cite{ls:pi1s1,lb:pinsn,brunerie:thesis}, ordinary and generalized cohomology~\cite{lf:emspaces,bf:cellcoh-hott,cavallo:cohom-hott,br:rp-hott}, the Freudenthal suspension and Blakers--Massey theorems~\cite{ffll:blakers-massey}, the Serre and Atiyah--Hirzebruch spectral sequences~\cite{floris:thesis}, and localization at primes~\cite{cors:loc-hott,scoccola:nilpfrac-hott}.
These proofs, which constitute the subfield of \emph{synthetic homotopy theory}, are sometimes simply rewritings of classical ones, but often they contribute new ideas even to classical homotopy theory.
More importantly, their expected validity internal to any \io-topos should lead to many new applications, some of which have already been realized (e.g.~\cite{abfj:gen-blakers-massey,abfj:goodwillie-htt}).

Unfortunately, rigorous proofs of the interpretation of type theories into \io-toposes have lagged behind these internal developments.
The problem\footnote{In addition to this coherence problem, it is necessary to relate the type-theoretic syntax (e.g.\ with variable binding) to the fully strictified ``algebraic'' version.
  Known as the ``initiality principle'', this is essentially a straightforward but tedious bookkeeping problem.
  It is well-known to be true for simpler type theories, but has not been written down completely yet for more complicated ones like MLTT and Book HoTT, though it is universally expected to be unproblematic.} is that, unlike for 1-toposes, many equalities that hold strictly in type theory hold only up to homotopy in an \io-topos, giving rise to a \emph{coherence problem}.
The straightforward solution (though one may imagine others) is to replace an \io-topos by a strict model that satisfies all the same strict equalities as the type theory.

The insight of~\cite{aw:htpy-idtype,klv:ssetmodel} was that the existing notion of \emph{Quillen model category}~\cite{quillen:htpical-alg} is already almost sufficient\footnote{In fact,~\cite{gg:idtypewfs,lumsdaine:hit-model-strux} showed that it is almost \emph{necessary} as well: any model of Book HoTT has the structure of at least the cofibrant-fibrant objects in a model category.}.
In fact, Voevodsky~\cite{klv:ssetmodel} constructed a model of all of MLTT, with univalence, in the model category of simplicial sets, which is a strictification of the \io-topos of \oo-groupoids (the \io-categorical analogue of the 1-topos of sets).
Later work such as~\cite{gb:topsimpid,ak:htmtt,shulman:invdia,gk:univlcc,ls:hits} showed that any \io-topos (and indeed any locally cartesian closed, locally presentable \io-category) can be strictified to a Quillen model category that models \emph{almost} all of Book HoTT,\footnote{To be precise, even after finding an appropriate Quillen model category there is an additional step of making all the structure strictly stable under pullback; but this was solved quite generally by~\cite{lw:localuniv} (see also~\cite{awodey:natmodels}), generalizing the technique of~\cite{klv:ssetmodel}.
In \cref{sec:coherence} we will extend this method to deal with universes.} including most higher inductive types --- but \emph{not} including the all-important univalence axiom.
In~\cite{shulman:invdia,shulman:elreedy,cisinski:elegant,shulman:eiuniv} some very special classes of \io-toposes were shown to model univalence, but a construction applicable to all \io-toposes proved elusive.

Our main goal in this paper is to solve this problem, showing that every Grothendieck \io-topos can be strictified to a Quillen model category that models MLTT with univalence, and hence nearly all of Book HoTT.
(The principal remaining open question is whether the higher inductive types can be chosen such that the univalent universes are closed under them.)

\begin{rmk}\label{rmk:cubical}
  We work only with model categories, relying on previous work such as~\cite{aw:htpy-idtype,klv:ssetmodel,lw:localuniv,shulman:invdia,shulman:elreedy} to establish the connection with type-theoretic syntax.
  (With one exception: in \cref{sec:coherence} we sketch an extension of these results to universe types, since this does not appear to be in the literature yet.)
  And we will use classical homotopy theory (notably, simplicial sets), classical logic (including the axiom of choice)\footnote{In the ``metatheory'' where we construct models, that is.  The type theory being modeled is intuitionistic, as befits the internal language of an arbitrary topos.}, and Book HoTT with univalence as an axiom (as stated by Voevodsky).

  Another popular class of formal systems used in homotopy type theory are \emph{cubical type theories}~\cite{cchm:cubicaltt,ahw:chtt-i,ah:chtt-ii,abchhl:cart-cube,cm:unif-cartcube}, which enjoy computational advantages (e.g.\ univalence is a theorem rather than an axiom), and \emph{cubical set} models, which can be studied in constructive metatheories~\cite{bch:tt-cubical} and contain univalent universes~\cite{bch:univalence-cubical,ahh:chtt-iii} closed under higher inductive types~\cite{chm:cubical-hits,ch:chtt-iv}.
  However, the class of \io-categories presented by cubical models remains unclear.
  Indeed, it seems not even to be known which kinds of cubical sets can present the standard \io-topos of \oo-groupoids.
  Thus, there are still good reasons to study Book HoTT and its models in classical homotopy theory.
\end{rmk}

\begin{rmk}\label{rmk:pfthy}
  As in~\cite{klv:ssetmodel}, we are only able to obtain universes closed under the basic type formers of MLTT by assuming that inaccessible cardinals exist in our background set theory.
  It is true that classical ZFC set theory with inaccessibles is much stronger, proof-theoretically, than \emph{predicative} MLTT 
  with universes (even with univalence); the latter is only as strong as Kripke-Platek set theory with recursive inaccessibles~\cite{rathjen:pfthy-ua}.
  However, we will show in \cref{thm:propresizing} that our models also validate the \emph{propositional resizing principle}, which is impredicative and probably increases the consistency strength of type theory to nearly that of ZFC with inaccessibles.
  Thus our metatheory is not really much stronger than necessary.

  Of course, ZFC itself does not prove that any inaccessible cardinals exist.
  Their consistency is fairly uncontroversial among modern set theorists, who routinely study much stronger large cardinal axioms.
  However, for the reader who nevertheless prefers to avoid them, we note that our construction also produces univalent universes within ZFC, albeit with weaker closure properties in the internal type theory.
\end{rmk}

\subsection{Overview}
\label{sec:overview}

The model categories we will use are familiar from \io-topos theory~\cite{lurie:higher-topoi}: they are left exact left Bousfield localizations of injective model structures on categories of simplicial presheaves.
These have always seemed a promising choice for modeling type theory, and indeed they are a subclass of those previously known to model the rest of MLTT; all they are missing is universes.

The model-categorical version of a universe is a classifier for small fibrations: a fibration $\pi:\Util\to U$ with the property that every fibration satisfying some cardinality bound on its fibers occurs as a (strict, 1-categorical) pullback of $\pi$.
The problem is that until now no explicit description of the fibrations in a general injective model structure, let alone a localization thereof, has been known.
Thus, our main task will be to establish new characterizations of these fibrations.

For a long time the dominant tool in model category theory for describing classes of fibrations was \emph{cofibrant generation}, whereby a set $\cI$ of ``generating acyclic cofibrations'' determines the class of fibrations $f:X\to Y$ having the right lifting property against all morphisms in \cI, i.e.\ for any commutative square
\begin{equation}
  \begin{tikzcd}
    A \ar[d,"i"'] \ar[r,"g"] & X \ar[d,"f"]\\
    B \ar[r,"h"'] & Y
  \end{tikzcd}\label{eq:intro-sq}
\end{equation}
with $i\in \cI$, there exists a lifting $k:B\to X$ with $k i = g$ and $f k = h$.
It was shown in~\cite[A.3.3.3]{lurie:higher-topoi} that injective model structures are cofibrantly generated, but the generating acyclic cofibrations are very inexplicit, consisting essentially of all pointwise acyclic cofibrations whose domain and codomain satisfy some cardinality bound.
This makes it difficult to say anything concrete about the injective fibrations.

However, recently a more general \emph{algebraic} approach to model category theory has arisen~\cite{gt:nwfs,garner:soa,riehl:nwfs-model,bg:awfs-i,rosicky:acc-model}, in which the fibrations are characterized by admitting the structure of an algebra for some pointed endofunctor on the category of arrows.
The algebraic approach has proven particularly useful for modeling type theory~\cite{gb:topsimpid,cchm:cubicaltt,gs:uniform,awodey:qmc-cube}.
Our characterization of the injective fibrations mixes these two approaches: they are the pointwise fibrations (these are determined by a right lifting property that is fairly explicit) that are additionally algebras for a certain pointed endofunctor.

The endofunctor in question is nothing new: in homotopy theory it is called a \emph{cobar construction}~\cite{may:goils,meyer:bar_i}, and in 2-category theory the analogous construction is called a \emph{pseudomorphism coclassifier}~\cite{bkp:2dmonads,lack:codescent-coh}.
Given a square~\eqref{eq:intro-sq} in which $i$ is a pointwise acyclic cofibration and $f$ a pointwise fibration, we can choose pointwise diagonal lifts; and while these will not generally form a natural transformation, they do always form a \emph{homotopy coherent natural transformation} $B\cohto X$ which satisfies naturality ``up to all higher homotopies''.
Thus, to obtain an actual lift we need $f$ to have the property that such homotopy coherent transformations can be ``rectified'' to strictly natural ones.
But the cobar construction of $f$, which we denote $E f$, is a representing object for such homotopy coherent transformations, so that any such $B\cohto X$ corresponds bijectively to some strictly natural $B\to E f$.
Thus, to obtain a strict map $B\to X$ we simply need a suitable strict map $E f\to X$, i.e.\ $f$ must be an algebra for the pointed endofunctor $E$.

This idea is inspired by~\cite{lack:htpy-2monads}, who characterized the cofibrant objects in the dual \emph{projective} model structure, in the simpler 2-categorical case, as those admitting a \emph{coalgebra} structure for the dual pseudomorphism classifier.
In 2-category theory the latter are known as \emph{flexible} objects; thus the injectively fibrant objects may be called \emph{coflexible}.
Informally, a coflexible presheaf $X$ is equipped with a rectification of any ``pseudo-element'' --- meaning an element $x\in X_d$ equipped with a coherent family of homotopical replacements for its restrictions along all morphisms in the domain category --- into an ordinary element of $X_d$, in a coherent and natural way.

With this characterization in hand, there are two remaining problems to be solved.
The first is to actually use this description to construct a universe of injective fibrations.
The existing techniques for constructing universes from~\cite{klv:ssetmodel,streicher:ttuniv,shulman:elreedy,cisinski:elegant} are designed for the cofibrantly generated case, but can be generalized to the algebraic case by considering the algebra operation as a structure, rather than its mere existence as a property.

More specifically, the ``unstructured'' techniques require the fibrations to be ``local on the base'', in the sense that a morphism $f:X\to Y$ into a colimit $Y = \colim_i Y_i$ is a fibration as soon as its pullback to each $Y_i$ is.
This corresponds to the generating acyclic cofibrations having representable codomains, which is what fails for general injective model structures.
But the generalized ``structured'' technique instead allows us to assume that the pullback of $f$ to each $Y_i$ is equipped with an algebra structure, and moreover that each square relating these pullbacks along a transition morphism $Y_i\to Y_j$ is a \emph{morphism} of algebras.
This makes it much easier to show that the algebra structures can be ``glued together'' into an algebra structure on $f$, yielding a universe of injective fibrations.
(Similar phenomena have been observed for the ``uniform fibrations'' used in cubical set models, which are cofibrantly generated with representable codomains only in an algebraic sense~\cite{bch:tt-cubical,gs:uniform,awodey:qmc-cube}.)

The last remaining problem is to deal with left exact localizations.
Here the basic idea was already sketched in~\cite[Remarks 3.24 and A.29]{rss:modalities}: since the local objects can be described internally in the type theory of the un-localized model category, we can construct internally a ``universe of local objects''.
In the left exact case this universe is itself local, and thus supplies a universe for the localized model structure.
At the time of~\cite{rss:modalities} we could not quite carry this out, since we were unable to show that every left exact localization of an \io-topos yields an \emph{internal} left exact localization in its internal type theory.
Fortunately,~\cite{abfj:lexloc} have recently given an explicit characterization of left exact localizations of \io-toposes, which passes to the internal type theory and thus allows us to prove this result, and hence construct universes for such localizations.
(For consistency in exposition, in this paper we will express the construction entirely at the model-categorical level without using the internal type theory, but the underlying idea is as described above.)

\subsection{\Ttmts}
\label{sec:ttmts-intro}

This much would be sufficient to prove our main result that homotopy type theory with univalence can be interpreted into Grothendieck \io-toposes.
However, with an eye towards the future, we do not restrict attention only to simplicial presheaf categories, but introduce a further abstraction.
Inspecting the proof, we see that the fact that we are working in a presheaf category is never really used (at least if we use the method of~\cite{shulman:elreedy} to construct universes by ``cofibrant replacement'', rather than the explicit methods of~\cite{klv:ssetmodel,streicher:ttuniv}).
What we do frequently need is the fact that the cofibrations satisfy certain exactness properties characteristic of the monomorphisms in a 1-topos.

Thus, we will define a \emph{\ttmt} to be a model category $\E$ whose underlying category is a \slcc Grothendieck 1-topos\footnote{Thus there is something of a pun in the name ``\ttmt'': it is a (1-)topos with a model structure that presents an (\io-)topos.}, whose model structure is proper, simplicial, and combinatorial, with its cofibrations being the monomorphisms, and that is furthermore equipped with a good notion of ``algebraic structured fibration''.
The above constructions can then be factored into the following results:
\begin{enumerate}[label=(\arabic*)]
\item Simplicial sets are a \ttmt (this is~\cite{klv:ssetmodel}).\label{item:it1}
\item \Ttmts are closed under passage to presheaves on simplicial categories, with injective model structures.\label{item:it2}
\item \Ttmts are closed under left exact localizations.\label{item:it3}
\item Every \ttmt has the structure to model homotopy type theory with univalent universes.
\end{enumerate}
Steps~\ref{item:it1}--\ref{item:it3} imply that every Grothendieck \io-topos can be strictified into some \ttmt, namely a left exact localization of an injective model structure on simplicial presheaves.
And conversely, every \ttmt presents a Grothendieck \io-topos, indeed it is a \emph{model topos}~\cite{rezk:homotopy-toposes}.

This additional step of abstraction is analogous to the step from defining toposes themselves via presentations (a category equivalent to a left exact localization of a presheaf category) to defining them intrinsically (a locally presentable category satisfing Giraud's exactness axioms).
\Ttmts also admit analogues of the basic constructions of elementary topos theory: slice categories, presheaves on internal as well as enriched categories, left exact localizations, coalgebras for suitable comonads, and Artin gluing.
In particular, they are closed under passage to ``internal sheaves on internal sites'', so that we are free to use any \ttmt as the ``ambient base topos'' instead of simplicial sets.

We have no use for this extra generality as yet, but the type theory we consider here has various extensions that we would also like to interpret in \io-toposes, and it might happen that some ``nonstandard'' \ttmts are more convenient for some such purpose.
This includes:
\begin{itemize}
\item \emph{Two-level type theories} such as~\cite{ack:2ltt}, especially with axioms such as fibrancy of the natural numbers object.
  The latter holds in enriched simplicial presheaf categories, but fails in general in their left exact localizations.
  However, it might hold in a model structure on simplicial \emph{sheaves} (though it seems doubtful that such models would suffice to present all \io-toposes).
\item \emph{Modal type theories} as in~\cite{ls:1var-adjoint-logic,lsr:multi,lsr:depdep}, which represent a collection of toposes connected by functors.
  Greater freedom to change the model categories seems likely to help in finding sufficiently strict models of such.
\item Our universes are closed under the standard type formers of MLTT, but it remains to be shown that they are closed under higher inductive types (the construction of the latter in~\cite{ls:hits} does not preserve smallness).
  It could be that some \ttmts make this easier to achieve than others.
  For instance, the method used to build higher inductive types in cubical set models~\cite{chm:cubical-hits,ch:chtt-iv} may be adaptable to simplicial sets, and perhaps to some other \ttmts with good properties.
\end{itemize}

\begin{rmk}
  The discussion above, and the rest of the paper, deals only with \emph{Grothendieck} \io-toposes, i.e.\ left exact localizations of presheaf \io-categories.
  Eventually one would also like to also model HoTT in ``elementary \io-toposes'', a hypothetical notion with proposed definitions~\cite{shulman:eleminf-talk,rasekh:eleminf} but essentially no known examples (other than Grothendieck ones).
  Also, sometimes the ``internal language correspondence'' refers to the even stronger conjecture that some ``homotopy theory of type theories''~\cite{kl:hot-tt} is \emph{equivalent} to that of elementary \io-toposes.

  The tools developed here for the Grothendieck case may be useful for proving these stronger conjectures, e.g.\ by embedding an elementary \io-topos in its presheaf category (as done in~\cite{ks:intlang-lex} to prove the strong internal language correspondence for \io-categories with finite limits).
  However, while these stronger conjectures are theoretically interesting, for most practical applications all that is needed is the interpretation of type theory into Grothendieck \io-toposes, which we establish here.
  For instance, this interpretation implies that the type-theoretic proof of the Blakers-Massey theorem~\cite{ffll:blakers-massey} \emph{is already} a proof of the corresponding \io-topos-theoretic theorem, without needing to be manually translated into \io-categorical language as in~\cite{abfj:gen-blakers-massey}.
\end{rmk}

\subsection{Outline}
\label{sec:outline}

We begin in \cref{sec:2cat} with 2-categorical preliminaries.
Then in \cref{sec:nfs} we study classes of structured morphisms, and in \cref{sec:relpres} we investigate how to restrict the cardinality of their fibers in an abstract way.
In \cref{sec:univalence} we construct universes for suitable such classes, and prove that in a suitable model category they are fibrant and univalent.
In \cref{sec:ttmt} we collect the accumulated hypotheses of these theorems into the notion of \ttmt.

We then move on to injective model structures and left exact localizations.
In \cref{sec:coh-bar} we motivate the (co)bar construction in more detail and review its basic properties, and then in \cref{sec:injmodel} we use it to give our new description of injective fibrations and prove that injective model structures are \ttmts.
Finally, in \cref{sec:int-loc} we prove some preliminary results about internal localizations (a generalization of left exact ones) and in \cref{sec:lex-loc} we show that \ttmts are preserved by left exact localization.
In \cref{sec:coherence} we deal with coherence for universes, which does not appear in any extant references. 

The main prerequisite for the reader is familiarity with model category theory; good modern references include~\cite{hovey:modelcats,hirschhorn:modelcats,mp:more-concise,riehl:cht}.
We will also use basic notions of 2-category theory, for which~\cite{steve:companion} is an excellent introduction, and some theory of locally presentable categories, whose standard reference is~\cite{ar:loc-pres}.

\subsection{Acknowledgments}
\label{sec:acknowledgments}

I am very grateful to Mathieu Anel, Steve Awodey, Dan Christensen, Eric Finster, Jonas Frey, Andr\'{e} Joyal, Peter LeFanu Lumsdaine, Emily Riehl, Bas Spitters, and Raffael Stenzel for helpful conversations, feedback, and detailed and careful reading of drafts.
And I would like to once again thank Peter May, for the two-sided simplicial bar construction.


%% file: 2cat.tex
\section{2-categorical preliminaries}
\label{sec:2cat}

We begin with some 2-categorical observations.
A morphism $f:\dX\to\dY$ in a 2-category \sK is an \textbf{internal fibration} if each induced functor $\sK(\dZ,\dX) \to \sK(\dZ,\dY)$ is a Grothendieck fibration.
Explicitly, this means for any $g:\dZ\to \dX$ and 2-cell $\alpha: h \to f\circ g$ there is a cartesian lift $\beta : k\to g$ with $f\circ k = h$ and  $f\circ \beta = \alpha$.
Similarly we have \textbf{internal discrete fibrations}, for which $\beta$ must be unique.

Let \E be a (large, locally small) category and write $\Ehat = \cPsFun(\E\op,\cGPD)$ for the (very large) 2-category of contravariant pseudofunctors from \E to large groupoids, and pseudonatural transformations between them.

\begin{defn}\label{defn:dfib}
  A \textbf{strict fibration} in \Ehat is a \emph{strictly} natural transformation $\dX\to\dY$ such that each component $\dX(A)\to \dY(A)$ is a fibration of groupoids.
  Similarly we define a \textbf{strict discrete fibration}.
\end{defn}

Of course, if \dX and \dY are discrete (i.e.\ functors $\E\op\to\nSET$), then any morphism $\dX\to\dY$ is a strict discrete fibration.
More generally, we have:

\begin{lem}\label{thm:fib-eqv}
  Any morphism $f:\dZ\to\dY$ in \Ehat is equivalent to a strict fibration over \dY, which is discrete if and only if $f$ is faithful.
\end{lem}
\begin{proof}
  Given $f:\dZ\to\dY$, each component $\dZ(A) \to \dY(A)$ is equivalent to a fibration $\dX(A) \to \dY(A)$, and we can transfer the functorial action of \dZ across these equivalences to make $\dX$ an object of $\Ehat$ and factor $f$ as an equivalence $\dZ\to\dX$ followed by a pseudonatural transformation $\dX\to\dY$ whose components are fibrations.
  Then we can use isomorphism-lifting for these fibrations to modify the pseudofunctorial action of $\dX$ to make this second transformation strictly natural.
  Finally, a fibration of groupoids is discrete if and only if it is faithful, which is invariant under equivalence.
\end{proof}

\begin{lem}
  A strict fibration $f:\dX\to\dY$ is an internal fibration in $\Ehat$, and similarly in the discrete case.
\end{lem}
\begin{proof}
  Given $g:\dZ\to \dX$ and $h:\dZ\to\dY$ with $\alpha : h\cong f \circ g$, for each $A$ there is a functor $k_A : \dZ(A) \to \dX(A)$ and isomorphism $\beta_A : k_A \cong g_A$ (unique in the discrete case) such that $f_A \circ k_A = h_A$ and $f_A \circ \beta_A = \alpha_A$.
  Transferring the pseudonaturality constraints of $g$ along the isomorphisms $\beta$ makes $k$ a pseudonatural transformation and $\beta$ a modification such that $f\circ k = h$ and $f\circ \beta = \alpha$.
\end{proof}

\begin{lem}\label{thm:dfib-pb}
  Given a strict fibration $f:\dX\to\dY$ and a morphism $g:\dZ\to\dY$, for each $A\in\E$ define $\dW(A)$ by the strict pullback of groupoids on the left:
  \[
    \begin{tikzcd}
      \dW(A) \ar[r,"h_A"] \ar[d,"k_A"'] \drpullback & \dX(A) \ar[d,"f_A"] \\
      \dZ(A) \ar[r,"g_A"'] & \dY(A)
    \end{tikzcd}
    \hspace{2cm}
    \begin{tikzcd}
      \dW \ar[r,"h"] \ar[d,"k"'] \drpullback & \dX \ar[d,"f"] \\
      \dZ \ar[r,"g"'] & \dY.
    \end{tikzcd}
  \]
  Then \dW can be made into a pseudofunctor and $h$ a pseudonatural transformation such that $k$ is a strict fibration and the square on the right commutes strictly in \Ehat and is a weak bicategorical pullback there.
\end{lem}
We will refer to this construction as a \textbf{strict pullback}.
\begin{proof}
  We define the pseudofunctorial actions of \dW and the pseudonaturality constraints of $h$ by lifting the pseudonaturality constraints of $g$ along the components of $f$.
  All the claims are straightforward to verify.
\end{proof}

Henceforth we suppose \E has pullbacks.
The following notion is fairly standard.

\begin{defn}\label{defn:rep}
  A morphism $f:\dW\to\dY$ in \Ehat is \textbf{representable} if for any $Z\in\E$ and weak bicategorical pullback
  \begin{equation*}
    \begin{tikzcd}
      \dP \ar[r] \ar[d] \ar[dr,phantom,near start,"\lrcorner"] & \dW\ar[d,"{f}"] \\
      \E(-,Z) \ar[r] & \dY
    \end{tikzcd}
  \end{equation*}
  the object \dP is equivalent to a representable $\E(-,X)$.
\end{defn}

Note that a morphism with representable codomain is a representable morphism if and only if it has representable domain.
In addition, every representable morphism is faithful, hence (by \cref{thm:fib-eqv}) equivalent to a strict discrete fibration.
Moreover, if $f$ is a representable strict discrete fibration, then in \cref{defn:rep} it suffices to consider strict pullbacks as in \cref{thm:dfib-pb}, and such a strict pullback \dP must be \emph{isomorphic} to a representable.
We therefore mainly consider representable strict discrete fibrations, which encompass all representable morphisms up to equivalence but are simpler to work with.

\begin{defn}
  Let $\cE\in\Ehat$ denote the \textbf{core of the self-indexing} of \E, where $\cE(Y)$ is the maximal subgroupoid of the slice category $\E/Y$, with pseudofunctorial action by pullback.
  Similarly, define $\cEp \in\Ehat$ such that $\cEp(Y)$ is the groupoid of morphisms $f:X\to Y$ equipped with $s:Y\to X$ such that $f s = \id_Y$.
\end{defn}

\begin{prop}\label{thm:univrep}
  The forgetful map $\varpi:\cEp\to\cE$ is the pseudo-universal representable morphism and the strictly universal representable strict discrete fibration.
  That is, the following are equivalent for any $\dY\in\Ehat$:
  \begin{enumerate}
  \item The hom-groupoid $\Ehat(\dY,\cE)$.\label{item:urep1}
  \item The groupoid of representable strict discrete fibrations with codomain \dY.\label{item:urep2}
  \item The 2-groupoid of representable morphisms with codomain \dY.\label{item:urep3}
  \end{enumerate}
  The functor~\ref{item:urep1}$\to$\ref{item:urep2} is by strict pullback and~\ref{item:urep2}$\to$\ref{item:urep3} is by inclusion, so that~\ref{item:urep1}$\to$\ref{item:urep3} is by weak bicategorical pullback.
\end{prop}
\begin{proof}
  Assuming that pullbacks in \E are globally chosen, $\varpi$ is strictly natural, and it is straightforward to check that it is a discrete fibration.
  A morphism $\E(-,Z) \to \cE$ corresponds, by the bicategorical Yoneda lemma, to a morphism $X\to Z$, and the pullback of $\cEp$ to $\E(-,Z)$ is then $\E(-,X)$.
  Thus $\varpi$ is a representable strict discrete fibration, hence so are its strict pullbacks.

  Note that~\ref{item:urep2} is a groupoid rather than a 2-groupoid since there are no nonidentity 2-cells between discrete fibrations, and it is equivalent to~\ref{item:urep3} by \cref{thm:fib-eqv}.
  Now for any representable strict discrete fibration $f:\dX\to\dY$, choose for each $y\in \dY(Z)$ a strict pullback square
  \[
    \begin{tikzcd}
      \E(-,P_y) \ar[d] \ar[r] \drpullback & \dX \ar[d,"f"] \\
      \E(-,Z) \ar[r,"y"'] & \dY.
    \end{tikzcd}
  \]
  Sending $y\mapsto (P_y\to Z)$ then defines a pseudonatural transformation $\dY\to\cE$, yielding a functor~\ref{item:urep2}$\to$\ref{item:urep1}, which can be verified to be an inverse equivalence to pullback.
\end{proof}

We conclude this section with bicategorical lifting and orthogonality.

\begin{defn}\label{defn:2liftorth}
  Let $i:\dA\to\dB$ and $p:\dX\to\dY$ be morphisms in a 2-category \sK.
  We say $i$ and $p$ have the \textbf{lifting property}, and write $i \lifts p$, if the functor
  \begin{equation}
    \sK(\dB,\dX) \to \sK(\dA,\dX) \bitimes_{\sK(\dA,\dY)} \sK(\dB,\dY)\label{eq:bicat-lift}
  \end{equation}
  is essentially surjective, where $\bitimes$ denotes a weak bicategorical pullback.
  In other words, $i \lifts p$ if given any $f:\dA\to\dX$ and $g:\dB\to\dY$ and isomorphism $\al : p f \cong g i$, there exists a morphism $h:\dB\to \dX$ and isomorphisms $\be:f \cong h i$ and $\gm : p h \cong g$ such that $(\gm i).(p\be) = \al$.

  If instead~\eqref{eq:bicat-lift} is an equivalence of categories, we say $i$ and $p$ are \textbf{orthogonal} and write $i\perp p$.
  In addition to the lifting property, this says that given $h,k:\dB\to\dX$ with 2-cells $\ph:h i \to k i$ and $\psi : p h \to p k$ such that $\psi i = p \ph$, there exists a unique $\chi : h\to k$ with $\chi i = \ph$ and $p \chi = \psi$.

  For an object \dX, we write $i\lifts \dX$ and $i\perp \dX$ to refer to lifting and orthogonality properties for the unique morphism $\dX\to 1$.
\end{defn}

\begin{lem}\label{thm:dfib-perp}
  If $p:\dX\to\dY$ is an internal discrete fibration in \sK, the following are equivalent for any $i:\dA\to\dB$ in \sK.
  \begin{enumerate}
  \item $i \perp p$.\label{item:dp1}
  \item The analogous functor
    \begin{equation}
      \sK(\dB,\dX) \to \sK(\dA,\dX) \times_{\sK(\dA,\dY)} \sK(\dB,\dY)\label{eq:bicat-strict-lift}
    \end{equation}
    to the strict pullback is an isomorphism.\label{item:dp2}
  \item The functor~\eqref{eq:bicat-strict-lift} is bijective on objects.\label{item:dp3}
  \item Given any $f:\dA\to\dX$ and $g:\dB\to\dY$ such that $p f = g i$, there exists a unique morphism $h:\dB\to \dX$ such that $p h = g$ and $h i = f$.\label{item:dp4}
  \end{enumerate}
\end{lem}
\begin{proof}
  If $p$ is a fibration, then so is the induced map $\sK(\dA,\dX) \to \sK(\dA,\dY)$, and thus the bicategorical pullback is equivalent to the strict one.
  Since isomorphisms are also equivalences,~\ref{item:dp2}$\Rightarrow$\ref{item:dp1}.

  On the other hand, assuming~\ref{item:dp1}, then if $p f = g i$ we can choose $\alpha = \id$ to obtain an $h:\dB\to \dX$ and isomorphisms $\be:f \cong h i$ and $\gm : p h \cong g$ such that $(\gm i).(p\be) = \id$.
  But now since $p$ is a fibration, there is an isomorphism $\delta:h\cong h'$ such that $p h' = g$ and $p \delta = \gm$.
  Then we have $(\delta i).\be : f \cong h' i$, and $p((\delta i).\be) = (p\delta i).(p\be) = (\gm i).(p\be) = \id$; thus since $p$ is a \emph{discrete} fibration $(\delta i).\be = \id$ and so $f = h' i$.
  Thus~\eqref{eq:bicat-strict-lift} is a \emph{surjective} equivalence.
  Furthermore, if $h,k :\dB\to \dY$ are such that $p h = p k$ and $h i = k i$, then~\ref{item:dp1} gives a unique isomorphism $\xi : h\cong k$ with $p \xi = \id$ and $\xi i = \id$.
  But since $p$ is a discrete fibration, $p \xi = \id$ implies $\xi=\id$ and $h=k$; thus~\ref{item:dp1}$\Rightarrow$\ref{item:dp2}.

  Now of course~\ref{item:dp2}$\Rightarrow$\ref{item:dp3}.
  But conversely, if~\ref{item:dp3} holds and $h,k:\dB\to\dY$ are given with $\ph:h i \to k i$ and $\psi : p h \to p k$ such that $\psi i = p \ph$, then since $p$ is a discrete fibration there is a unique $\chi : k'\to k$ with $p k' = p h$ and $p \chi = \psi$, and uniqueness furthermore implies $k' i = h i$ and $\chi i = \ph$.
  So~\ref{item:dp3} implies that in fact $k' = h$; hence the unique $\chi$ shows that~\eqref{eq:bicat-strict-lift} is fully faithful, i.e.~\ref{item:dp2} holds.

  Finally,~\ref{item:dp4} is just an explicit restatement of~\ref{item:dp3}.
\end{proof}



\begin{eg}\label{eg:orth-cancel}
  As in the 1-categorical case, morphisms having a right lifting or right orthogonality property are closed under composition and (weak bicategorical) pullback.
  Similarly, morphisms with a right orthogonality property are cancellable: if $i\perp qp$ and $i\perp q$ then $i\perp p$.
\end{eg}

\begin{eg}\label{eg:colim-orth}
  Let $Y = \colim_i Y_i$ be a diagram in a category \E, and let $\Yhat = \bicolim_i \E(-,Y_i)$ be the corresponding weak bicategorical colimit of representables in \Ehat.
  Its universal property is that
  \[\Ehat(\Yhat,\dX) \simeq \bilim_i \Ehat(\E(-,Y_i),\dX) \simeq \bilim_i \dX(Y_i)\]
  where $\bilim_i$ denotes the weak bicategorical limit: an object of $\bilim_i \dX(Y_i)$ consists of objects $x_i \in \dX(Y_i)$ and coherent isomorphisms $\dX_\iota(x_j) \cong x_i$ for all $\iota:i\to j$.

  In particular, the given colimit cocone in \E induces a map $q:\Yhat \to \E(-,Y)$.
  Then for $\dX\in\Ehat$ we have $q \perp \dX$ if and only if \dX preserves the colimit $Y = \colim_i Y_i$, in that the induced map $\dX(Y) \to \bilim_i \dX(Y_i)$ is an equivalence.
  Note that all representables $\E(-,Z) \in \Ehat$ preserve all colimits.
\end{eg}


%% file: nfs.tex
\section{Notions of fibred structure}
\label{sec:nfs}

We now introduce the notions of ``structured morphism'' that our universes will classify.
Let \E, \Ehat, and $\cE$ be as in \cref{sec:2cat}.

\begin{defn}\label{defn:fcos}
  A \textbf{\nfs} on \E is a strict discrete fibration $\phi:\F\to\cE$ with codomain \cE in $\Ehat$ that has small fibers.
\end{defn}

Explicitly, this means that for any morphism $f:X\to Y$ we have a set (possibly empty, but not a proper class) of \textbf{\F-structures} on $f$, and for any pullback square
\[
  \begin{tikzcd}
    X' \ar[d,"f'"'] \ar[r] \ar[dr,phantom, near start,"\lrcorner"] & X \ar[d,"f"]\\
    Y' \ar[r] & Y
  \end{tikzcd}
\]
we have a function from \F-structures on $f$ to \F-structures on $f'$, which is pseudofunctorial in pullback squares.
When $f$ is equipped with a chosen \F-structure we call it an \textbf{\F-algebra}, and if in the above pullback square $f'$ has the \F-structure induced from $f$ we call the square an \textbf{\F-morphism}.

\begin{eg}\label{eg:full-fcos}
  If $\F\into\cE$ is the inclusion of a subfunctor, then \F is just a pullback-stable class of morphisms in \E.
  We call this a \textbf{full} \nfs.
  This includes most previous work on universes, e.g.~\cite{klv:ssetmodel,shulman:elreedy,cisinski:elegant,stenzel:thesis}.

  In particular, when $\F=\cE$ we have the \textbf{trivial} \nfs.
\end{eg}

\begin{eg}\label{eg:pb-fcos}
  If $\F_1$ and $\F_2$ are \nfss, then so is the pullback $\F_1\times_\cE \F_2$.
  An $(\F_1\times_\cE \F_2)$-structure on a morphism is just a pair consisting of an $\F_1$-structure and an $\F_2$-structure.
\end{eg}

\begin{eg}\label{eg:retr-fcos}
  The forgetful morphism $\varpi:\cEp \to \cE$ from \cref{thm:univrep} makes $\cEp$ into a \nfs: a $\cEp$-structure on a morphism is a section of it.
\end{eg}


\begin{eg}\label{eg:preimage-fcos}
  Let $G:\E_1\to\E_2$ be a pullback-preserving functor and \F a \nfs on $\E_2$.
  Then there is a \textbf{preimage} \nfs $G^{-1}(\F)$ on $\E_1$, where a $G^{-1}(\F)$-structure on $f$ is by definition an \F-structure on $G(f)$.
\end{eg}

\begin{eg}\label{eg:fend-fcos}
  Let $H$ be a \emph{fibred core-endofunctor} of \E, i.e.\ a family of endofunctors $H_Y$ of the core of $\E/Y$ that commute with pullback up to coherent isomorphism.
  Then $H$ induces a map $\cE\to\cE$ in $\Ehat$, and any \nfs $\F\to\cE$ yields a \textbf{pullback} \nfs $H^*(\F)$, where an $H^*(\F)$-structure on a morphism $X\to Y$ is an \F-structure on $H_Y(X) \to Y$.
\end{eg}

\begin{eg}\label{eg:ff-fcos}
  Suppose given a \textbf{functorial factorization} on \E, i.e.\ a functor $E:\E^\dtwo \to \E^\dthree$ sending each morphism $f:X\to Y$ (regarded as an object of $\E^\dtwo$) to a composable pair $X \xto{\lambda_f} E f \xto{\rho_f} Y $ such that $f = \rho_f \lambda_f$.
  The functorial aspect factors every commutative square on the left as a pair of such on the right:
  \begin{equation}\label{eq:commsqfact}
    \begin{tikzcd}
      X' \ar[d,"f'"'] \ar[r,"g"] & X \ar[d,"f"]\\
      Y' \ar[r,"h"] & Y
    \end{tikzcd}
    \qquad = \qquad
    \begin{tikzcd}
      X' \ar[d,"{\lambda_{f'}}"'] \ar[r,"g"] & X \ar[d,"{\lambda_f}"]\\
      E f' \ar[d,"{\rho_{f'}}"'] \ar[r,"{\fact g h}"] & E f \ar[d,"{\rho_f}"]\\
      Y' \ar[r,"h"] & Y.
    \end{tikzcd}
  \end{equation}
  Define an \textbf{$\dR_E$-structure} on $f:X\to Y$ to be a retraction $r_f:E f \to X$ such that $r_f \lambda_f = \id_X$ and $f r_f = \rho_f$, exhibiting $f$ as a retract of $\rho_f$ in $\slice{\E}{Y}$.
  For example, if $E$ underlies a weak factorization system, then the morphisms that admit $\dR_E$-structures are those in the right class (but a given morphism in the right class will generally admit more than one $\dR_E$-structure).

  If the left square in~\eqref{eq:commsqfact} is a pullback and $f$ is an $\dR_E$-algebra, then the pullback universal property gives a unique $r_{f'} : E f' \to X'$ such that $f' r_{f'} = \rho_{f'}$ and $g r_{f'} = r_f \fact g h$; then $r_{f'} \lambda_{f'} = \id_{X'}$ follows by uniqueness.
  This defines the functorial action making $\dR_E$ a \nfs.
\end{eg}

\begin{eg}\label{eg:rep-fcos}
  For any morphism $\pi:\Util\to U$ in \E, its classifying map $\E(-,U) \to \cE$ 
  is faithful, hence equivalent to a strict discrete fibration $\dRep_\pi \to U$, which
  we call its \textbf{represented} \nfs.
  A $\dRep_\pi$-structure on $f:X\to Y$ is a morphism $Y\to U$ exhibiting $f$ as a pullback of \pi.
\end{eg}

It turns out that the property of being ``local on the base'', which we need to construct universes, is captured exactly by representability.
(Representable morphisms also have other uses in modeling type theory; see~\cite{awodey:natmodels} and \cref{sec:coherence}.)

\begin{defn}\label{defn:local}
  A \nfs \F is \textbf{\local}
  if the strict discrete fibration $\phi:\F\to\cE$ is representable (\cref{defn:rep}).
\end{defn}


Explicitly, this means that given any map $X\to Z$, there is a map $\phi^\F_X : \F_X \to Z$ in \E such that for any $g:Y\to Z$, there is a natural bijection between \F-structures on $g^*X$ and lifts of $g$ to $\F_X$.

\begin{eg}
  Representable morphisms are closed under pullback and composition in \Ehat.
  Thus, if $\F_1$ and $\F_2$ are \local so is $\F_1\times_\cE \F_2$; and if \F is \local so is $H^*(\F)$ for any fibred core-endofunctor $H$.
\end{eg}

\begin{eg}
  If $G:\E_1\to\E_2$ preserves pullbacks and has a right adjoint $H$, and $\F$ is a \local \nfs on $\E_2$, then $G^{-1}(\F)$ is a \local \nfs on $\E_1$.
  For given $X\to Z$ in $\E_1$ and $g:Y\to Z$, to give $g^*X$ a $G^{-1}(\F)$-structure is equivalently (since $G$ preserves pullbacks) to give $(G g)^*(G X)$ an \F-structure, i.e.\ to lift $G g : G Y \to G Z$ to the representing object $\F_{G X}$; but this is the same as to lift $g$ to the pullback
  \[
    \begin{tikzcd}
      (G^{-1}(\F))_X \ar[r] \ar[d]\drpullback & H (\F_{G X}) \ar[d]\\
      Z \ar[r] & H G Z.
    \end{tikzcd}
  \]
\end{eg}

\begin{eg}\label{thm:ep-local}
  By \cref{thm:univrep}, the \nfs $\cEp$ is \local.
  Moreover, every \local \nfs arises as $H^*(\cEp)$ for some fibred core-endofunctor $H$, namely $H(X) = \F_X$ in the above notation.
\end{eg}

\begin{eg}\label{eg:cart-ff-local}
  Let \E be locally cartesian closed, and call a functorial factorization \textbf{cartesian} if it preserves pullback squares, in that if the left-hand square in~\eqref{eq:commsqfact} is a pullback, so is the lower right-hand square (hence so also is the upper one).
  Under these assumptions, the \nfs $\dR_E$ from \cref{eg:ff-fcos} is \local.
  To see this, given $f:X\to Z$ with factorization $X\xto{\lambda_f} E f\xto{\rho_f} Z$, using the locally cartesian closed structure we can build an object $\Retr_Z(\lambda_f)$ of $\E/Z$ such that for any $g:Y\to Z$, lifts of $g$ to $\Retr_Z(\lambda_f)$ are naturally bijective with retractions of $g^*(\lambda_f)$ in $\E/Y$.
  Since $E$ is cartesian, $g^*(\lambda_f) = \lambda_{g^*(f)}$, so these are $\dR_E$-structures on $g^*X \to Y$, i.e.\ $\Retr_Z(\lambda_f)$ is the desired representing object.
\end{eg}

\begin{eg}\label{eg:rep-local}
  If \E is locally cartesian closed, then the represented \nfs (\cref{eg:rep-fcos}) determined by any map $\pi:\Util\to U$ is \local.
  To see this, let $\E(-,Z) \to \cE$ classify a map $f:X\to Z$.
  Using local cartesian closure, there is an object $\Iso(X,\Util)$ of $\E/(Z\times U)$ such that for any $g:Y\to Z$ and $h:Y\to U$, lifts of $(g,h)$ to $\Iso(X,\Util)$ are naturally bijective with isomorphisms $g^*X \cong h^*\Util$.
  Put differently, for any $g:Y\to Z$, lifts of $g$ along the composite $\Iso(X,\Util) \to Z\times U \to Z$ are naturally bijective with pullback squares
  \[
    \begin{tikzcd}
      g^* X \ar[d]\ar[r] \drpullback & \Util \ar[d] \\
      Y \ar[r] & U.
    \end{tikzcd}
  \]
  But this says exactly that $\Iso(X,\Util)$ is the desired representing object $(\dRep_\pi)_X$.
\end{eg}

\begin{eg}\label{eg:pshf-can}
  Let \F be the full \nfs on a pullback-stable class (\cref{eg:full-fcos}) in a presheaf category $\E=\prcs$.
  Then \F is \local if and only if a morphism $f:X\to Y$ belongs to \F as soon as all its pullbacks to representables $\C(-,c)$ do.
  (Thus our \locality includes the \emph{strongly proper} classes of fibrations from~\cite[Definition 3.7]{cisinski:elegant}.)

  For ``only if'', fullness of \F means that $\F_X\to Z$ is a monomorphism for any $f:X\to Z$.
  But if all pullbacks of $f$ to representables are in \F, all maps from representables into $Z$ factor through $\F_X$; thus $\F_X\cong Z$, hence $f\in \F$.

  For ``if'', given $f:X\to Z$, let $\F_X$ be the sub-presheaf of $Z$ containing all $z\in Z(c)$ such that the pullback of $f$ along $z:\C(-,c)\to Z$ is in \F.
  Then $g:Y\to Z$ factors through $\F_X$ if and only if the pullback of $g^*X$ along all $y:\C(-,c) \to Y$ lies in \F, which by assumption is the same as $g^*X \in \F$.
\end{eg}

\begin{eg}\label{eg:rep-cod}
  Let $\cI$ be a set of morphisms in $\prcs$ with representable codomains.
  Then by \cref{eg:pshf-can}, the class of morphisms with the right lifting property with respect to \cI is \local, since any lifting problem against $i\in \cI$ can be solved by first pulling back to the codomain of $i$, which is representable.
  For instance, 
  the class of \emph{Kan fibrations} in simplicial sets is \local.

  Another way to prove \locality in presheaf categories can be found in~\cite[Theorem 3.14]{cisinski:elegant}.
\end{eg}


Since \locality is a representability property, we expect that in good cases it can be ensured by an adjoint functor theorem.

\begin{prop}\label{thm:local}
  Let $\phi:\F\to\cE$ be a \nfs on a locally presentable category \E, and let $\sQ$ denote the class of morphisms $q:\Yhat \to \E(-,Y)$ from \cref{eg:colim-orth} for all small colimits $Y = \colim_i Y_i$ in \E.
  Then the following are equivalent.
  \begin{enumerate}
  \item \F is \local.\label{item:local1}
  \item $\sQ \perp \phi$.\label{item:local2}
  \item For any morphism $f:X\to Y$ and small colimit $Y \cong \colim_i Y_i$ with coprojections $q_i : Y_i \to Y$ and pullbacks\label{item:local3}
  \begin{equation}
    \begin{tikzcd}
      X_{j} \ar[r,"{p_{j,i}}"] \ar[d,"{f_{j}}"'] \drpullback &
      X_i \ar[d,"{f_i}"'] \ar[r,"{p_i}"] \drpullback & X \ar[d,"f"]\\
      Y_{j} \ar[r,"{q_{j,i}}"'] &
      Y_i \ar[r,"{q_i}"'] & Y,
    \end{tikzcd}\label{eq:colim-to-pb}
  \end{equation}
  if each morphism $f_i: X_i \to Y_i$ is given an \F-structure such that the squares on the left in~\eqref{eq:colim-to-pb} are \F-morphisms, then $f$ has a unique \F-structure such that the squares on the right in~\eqref{eq:colim-to-pb} are \F-morphisms.
  \end{enumerate}
\end{prop}
\begin{proof}
  Recall that \locality of \F means that in any pullback
  \begin{equation}
    \begin{tikzcd}
      \dP \ar[r] \ar[d,"{\phi^\F_X}"'] \ar[dr,phantom,near start,"\lrcorner"] & \F\ar[d,"{\phi}"] \\
      \E(-,Z) \ar[r,"X"'] & \cE.
    \end{tikzcd}\label{eq:local-pb-2}
  \end{equation}
  the object \dP is representable.
  Since $\phi$ is a small discrete fibration, so is $\phi^\F_X$; and since $\E(-,Z)$ is a small discrete object, \dP is also small and discrete, i.e.\ a presheaf $\E\op\to\nSet$.
  Thus, by the adjoint functor theorem, \dP is representable if and only if preserves small colimits, i.e.\  if and only if $\sQ \perp \dP$.
  And as noted in \cref{eg:colim-orth} representables preserve colimits, so $\sQ\perp \E(-,Z)$.
  Thus by \cref{eg:orth-cancel}, $\sQ \perp \dP$ if and only if $\sQ \perp \phi^\F_X$.

  Since right orthogonality is preserved by pullback, this is implied by~\ref{item:local2}.
  But the converse also holds, since each morphism in \sQ has a representable codomain, so that any lifting problem relating it to $\phi$ factors through some $\phi^\F_X$.
  Finally,~\ref{item:local3} is just an unraveling of~\ref{item:local2} according to \cref{thm:dfib-perp}\ref{item:dp4}.
\end{proof}

\begin{verbose}
It may be helpful to unravel the ``only if'' direction of \cref{thm:local} more explicitly.
If \F is \local, then given $f:X\to Y$ and $Y \cong \colim_i Y_i$ as in \cref{thm:local}, let $\F_X$ be the classifying object for \F-structures on $f$.
Then each given \F-structure on $f_i$ induces a map $Y_i \to \F_X$, and the left-hand squares in~\eqref{eq:colim-to-pb} being \F-morphisms mean that these maps assemble into a cone under the diagram $\{Y_i\}$.
Thus they induce a unique morphism $Y\to \F_X$, corresponding to a unique \F-structure on $f$ making the right-hand squares in~\eqref{eq:colim-to-pb} \F-morphisms.
\end{verbose}

\begin{rmk}
  For full \nfss on presheaf categories, the equivalence between the conditions of \cref{thm:local,eg:pshf-can} was observed in~\cite[Remark 4.4]{sattler:eqvext}.
  More generally, by \cref{thm:local} our \locality includes the \emph{locality} condition of~\cite[(A.2)]{sattler:eqvext}.
\end{rmk}


%% file: relpres.tex
\section{Relatively \ka-presentable morphisms}
\label{sec:relpres}

For size reasons, we cannot expect a universe to classify \emph{all} fibrations, only those with ``bounded cardinality''.
In~\cite{klv:ssetmodel,shulman:elreedy,cisinski:elegant} such a bound was imposed by explicit reference to the underlying sets of presheaves.
We will work more abstractly, and thus more generally, with ``internal'' categorical notions of size.

In this section \E will be a locally presentable category (often locally cartesian closed), and $\ka,\la,\mu,\nu$ will be regular cardinals.
Recall that $X\in\E$ is \textbf{\ka-presentable} (also called \textbf{\ka-compact}) if $\E(X,-):\E\to\nSet$ preserves \ka-filtered colimits.
A category \sC is \textbf{\ka-small} if $\ka>\card\C$ (the cardinality of the set of arrows of \C).

\begin{eg}\label{eg:pshf-pres}
  By~\cite[Example 1.31]{ar:loc-pres}, if $\E=\prcs$ is a presheaf category where \C is \ka-small, then $X\in\E$ is \ka-presentable if and only if it is a \ka-small colimit of representables, if and only if $\sum_c \card{X_c} < \ka$, and if and only if each $\card{X_c}<\ka$.
  (The hypothesis that \C is \ka-small is essential, however.)
  In~\cite{klv:ssetmodel,shulman:elreedy} these are called \textbf{\ka-small} objects.
\end{eg}

More generally, we have the following.

\begin{lem}[{\cite[Proposition 2.23]{low:univ-ct}}]\label{thm:pw-pres}
  If \E is locally \ka-presentable and \C is a \ka-small category, then the functor category $\func\C\E$ is locally \ka-presentable, and an object of $\func\C\E$ is \ka-presentable if and only if it is pointwise \ka-presentable in \E.
\end{lem}
\begin{proof}
  The first statement is~\cite[Corollary 1.54]{ar:loc-pres}.

  For the ``only if'' direction of the second, note that for any $c_0\in \C$, the ``evaluation at $c_0$'' functor $\mathrm{ev}_{c_0} : \func\C\E \to \E$ has a right adjoint given by $X \mapsto \{ X^{\C(c,c_0)} \}_{c\in \C}$.
  Since $\C(c,c_0)$-fold powers are \ka-small limits (as \C is \ka-small), they commute with \ka-filtered colimits, so this right adjoint is \ka-accessible.
  Hence its left adjoint $\mathrm{ev}_{c_0}$ preserves \ka-presentable objects.

  For the ``if'' direction, note that each $\mathrm{ev}_{c_0}$ also has a left adjoint given by $X \mapsto \{ \C(c_0,c) \cdot X \}_{c\in \C}$.
  Since $\mathrm{ev}_{c_0}$ preserves all colimits, this left adjoint preserves \ka-presentable objects.
  Now we note that any $A\in \func\C\E$ can be written as $A = \int^{c_0} \C(c_0,-) \cdot A_{c_0}$.
  Hence if $A$ is pointwise \ka-presentable, it is a \ka-small colimit of \ka-presentable objects in $\func\C\E$, hence \ka-presentable.
\end{proof}

We might hope to construct a universe of \ka-presentable objects for all sufficiently large regular cardinals \ka.
However, like other facts about locally presentable categories, this seems to only be possible in general when \ka has a certain ``large cofinality'' property.
For the reader's convenience we recall the basic characterizations of that property, and also add a new one that appears not to be in the literature.

\begin{prop}\label{thm:shrp}
  For regular cardinals $\la<\mu$, the following are equivalent.
  \begin{enumerate}
  \item Every \la-accessible category is \mu-accessible.\label{item:shrp1}
  \item For any set $X$ with $\card X <\mu$, the poset $P_\la(X)$ of subsets of $X$ of cardinality $<\la$ has a \mu-small cofinal subset.\label{item:shrp2}
  \item Every \mu-presentable object of a locally \la-presentable category is a \mu-small \la-filtered colimit of \la-presentable objects.\label{item:shrp3a}
  \end{enumerate}
\end{prop}
\begin{proof}
  \ref{item:shrp1}$\Leftrightarrow$\ref{item:shrp2} is~\cite[Theorem 2.11]{ar:loc-pres} and~\cite[\sect 2.3]{mp:accessible}, and~\ref{item:shrp2}$\Rightarrow$\ref{item:shrp3a} is~\cite[Proposition 2.3.11]{mp:accessible}. 
  To show~\ref{item:shrp3a}$\Rightarrow$\ref{item:shrp2}, note that $\nSet$ is locally \la-presentable, and its \ka-presentable objects are those of cardinality $<\ka$.
  Thus if $\card X<\mu$, then $X$ is \mu-presentable in \nSet, so by assumption we have $X \cong \colim_{i\in I} X_i$ with $I$ \mu-small and \la-filtered and each $\card {X_i}<\la$.
  Let $\cA = \{ q_i(X_i) \mid i\in I \}$ be the set of images of the coprojections $q_i:X_i\to X$; then $\cA \subseteq P_\la(X)$ and $\card\cA<\mu$.
  And for any $Y\subseteq X$ with $\card Y<\la$, the set $Y$ is \la-presentable, so the inclusion $Y\into X$ factors through some object $X_i$ in the \la-filtered colimit.
  But then $Y \subseteq q_i(X_i) \in \cA$; so \cA is cofinal.
\end{proof}

\noindent
When these conditions hold, one writes 
$\la\shlt\mu$. 
Then:
\begin{itemize}
\item The relation $\shlt$ is transitive (by \cref{thm:shrp}\ref{item:shrp1}).
\item 
  For any set of regular cardinals $\{\la_i\}$ there is a regular cardinal $\mu$ such that $\la_i\shlt\mu$ for all $i$ (and hence there are arbitrarily large such \mu).
  If \ka is inaccessible and each $\la_i<\ka$, the class $\{\mu\mid \forall i .(\la_i\shlt\mu)\}$ is unbounded below \ka.
\item If $\la<\ka$ and \ka is inaccessible, then $\la\shlt\ka$.
\item By~\cite[Theorem 2.4.9]{mp:accessible} or~\cite[Theorem 2.19]{ar:loc-pres},
for any accessible functor $F$ there is a \la such that for any $\mu\shgt\la$, the functor $F$ is \mu-accessible and preserves \mu-presentable objects.
\end{itemize}

\begin{rmk}
  Since $\aleph_1\nshlt\aleph_{\om+1}$ by~\cite[Example 2.13(8)]{ar:loc-pres}, \cref{thm:shrp}\ref{item:shrp3a} also fails in this case.
  In particular, a set of cardinality $\aleph_\om$ is $\aleph_{\om+1}$-presentable in \nSet, but is not a $\aleph_{\om+1}$-small $\aleph_1$-filtered colimit of $\aleph_1$-presentable objects.
  In addition, there exist accessible functors $F$ (e.g.\ the endofunctor of \nSet defined by $F(X) = X^I$ for some infinite set $I$, cf.~\cite[Remark 3.2(4)]{br:aec-acc}) for which there are arbitrarily large regular cardinals \mu such that $F$ does not preserve \mu-presentable objects.
  Thus, it seems we cannot avoid the relation $\shlt$ or something like it.
  The relation $\ll$ used in~\cite{lurie:higher-topoi} is \textit{a priori} stronger than $\shlt$, but coincides with it if the Generalized Continuum Hypothesis holds~\cite[Fact 2.5]{lrv:intsize}.
\end{rmk}

\begin{prop}\label{thm:pres-pb}
  For any locally presentable category \E, there is a \la such that for any $\ka\shgt \la$, $\E$ is locally \ka-presentable and the \ka-presentable objects in \E are closed under finite limits.
\end{prop}
\begin{proof}
  By \cref{thm:pw-pres}, if \E is locally \ka-presentable, then the \ka-presentable objects of the category $\E^{(\to\ot)}$ of cospans are the pointwise \ka-presentable ones.
  Thus, as soon as the pullback functor $\E^{(\to\ot)} \to \E$ preserves \ka-presentable objects, the \ka-presentable objects of \E are closed under pullbacks.
  But since this functor is a right adjoint, it is accessible, so there is a \la such that this occurs for all $\ka\shgt\la$; and we can choose \la large enough that the terminal object is also \la-presentable.
\end{proof}

\begin{prop}\label{thm:pres-sub}
  For any locally presentable category \E, there is a \la such that for any $\ka\shgt \la$, $\E$ is locally \ka-presentable and the \ka-presentable objects in \E are closed under finite limits and subobjects.
\end{prop}
\begin{proof}
  We first prove the result when \E is a Grothendieck 1-topos.
  In this case it suffices to take \la satisfying \cref{thm:pres-pb} and such that the subobject classifier $\Omega$ is \la-presentable, since any subobject of $X$ occurs as a pullback to $X$ of the universal subobject $1\to \Omega$.

  Now an arbitrary locally presentable category \E is a reflective subcategory of some Grothendieck (indeed presheaf) 1-topos \sT, say with inclusion functor $U$ and reflector $L\adj U$.
  Both $L$ and $U$ are accessible, so there is a \la such that for any $\ka\shgt \la$, the statement holds for \sT and $L$ and $U$ preserve \ka-presentable objects.
  Now if $A$ is a subobject of a \ka-presentable $X\in \E$, then $UA$ is a subobject of the \ka-presentable $UX\in \sT$.
  Hence $U A$ is \ka-presentable in \sT, so $A \cong L U A$ is \ka-presentable in \E.
\end{proof}

For constructing universes, we also need a ``fiberwise'' notion of size for morphisms.

\begin{defn}[{\cite[Definition 6.1.6.4]{lurie:higher-topoi}}]
  A morphism $X\to Y$ in \E is \textbf{relatively \ka-presentable} if $Z\times_Y X$ is \ka-presentable for any morphism $Z\to Y$ where $Z$ is \ka-presentable.
\end{defn}

Of course, a relatively \ka-presentable morphism with \ka-presentable codomain also has \ka-presentable domain.
Conversely, if \ka-presentable objects are closed under finite limits in \E (cf. \cref{thm:pres-pb}), then every morphism between \ka-presentable objects is relatively \ka-presentable, and an object $X$ is \ka-presentable just when the map $X\to 1$ is relatively \ka-presentable.

\begin{prop}\label{thm:relpres-detect}
  Let \E be locally \la-presentable and locally cartesian closed, with \cG a strong generating set of \la-presentable objects, and let $\ka\ge\la$.
  Then $f:X\to Y$ in \E is relatively \ka-presentable if and only if $Z\times_Y X$ is \ka-presentable for any morphism $g:Z\to Y$ where $Z\in\cG$.
\end{prop}
\begin{proof}
  Every \la-presentable object is \ka-presentable, so ``only if'' is trivial.
  Conversely, the \ka-presentable objects are the closure of \cG under \ka-small colimits; for as in~\cite[Theorem 1.11]{ar:loc-pres}, the \la-small colimits of \cG form a dense generator whose canonical diagrams are \la-filtered, and then as in~\cite[Remark 1.30]{ar:loc-pres} every \ka-presentable object is (a retract of) the colimit of a \ka-small subdiagram of its canonical diagram with respect to these.
  Thus, it suffices to show that if $Y = \colim_i Y_i$ is a \ka-small colimit, with each $Y_i$ (hence also $Y$) being \ka-presentable, and $f:X\to Y$ is such that each $Y_i \times_Y X$ is \ka-presentable, then $X$ is \ka-presentable.
  But since \E is locally cartesian closed, colimits are stable under pullback; so $X = \colim_i (Y_i \times_Y X)$ is a \ka-small colimit of \ka-presentable objects, hence \ka-presentable.
\end{proof}

\begin{cor}
  For any morphism $f:X\to Y$ in a locally presentable and locally cartesian closed category, there exists a regular cardinal \ka such that $f$ is relatively \ka-presentable.
\end{cor}
\begin{proof}
  Let \E be locally \la-presentable and \cG a set of representatives for isomorphism classes of \la-presentable objects.
  Then there are only a small set of morphisms $Z\to Y$ for objects $Z\in \cG$, hence there is a $\ka\ge \la$ such that all objects $Z\times_Y X$ are \ka-presentable.
  Hence by \cref{thm:relpres-detect}, $f$ is relatively \ka-presentable.
\end{proof}

\begin{eg}\label{eg:pshf-relpres}
  As noted in \cref{eg:pshf-pres}, in a presheaf category $\E=\prcs$ where \C is \ka-small, the \ka-presentable objects are the \ka-small colimits of representables.
  Since the representables are a strong generating set of \om-presentable objects, in this case $f:X\to Y$ is relatively \ka-presentable if and only if its pullback to any representable is \ka-presentable, hence if and only if all its fibers are \ka-small sets.
  Thus, in a presheaf category, for sufficiently large \ka the relatively \ka-presentable morphisms coincide with the \textbf{\ka-small morphisms} of~\cite{klv:ssetmodel,shulman:elreedy}.
\end{eg}

We now study the preservation properties of relatively \ka-presentable morphisms, which will yield the closure of universes under type forming operations.

\begin{lem}\label{thm:relpres-comp}
  For any regular cardinal \ka, the composite of relatively \ka-presentable morphisms is relatively \ka-presentable.
\end{lem}
\begin{proof}
  If $X\to Y\to Z$ are relatively \ka-presentable and we have $W\to Z$ with $W$ \ka-presentable, then $W\times_Z X \cong (W\times_Z Y) \times_Y X$ is also \ka-presentable.
\end{proof}

\begin{lem}\label{thm:pres-slice}
  If \E is a locally \la-presentable category, then for any $A\in \E$, a morphism $X\to A$ is a \la-presentable object of $\E/A$ if and only if $X$ is a \la-presentable object of \E.
\end{lem}
Note that $A$ is not required to be \la-presentable.
If $A$ \emph{is} \la-presentable, and \la-presentable objects are closed under pullbacks, then the \la-presentable objects of $\E/A$ will also coincide with the relatively \la-presentable morphisms of \E with codomain $A$; but in the general case this is not true.
\begin{proof}[Proof of \cref{thm:pres-slice}]
  First suppose $X$ is \la-presentable in \E, and $Y = \colim_i Y_i$ is a \la-filtered colimit in $\E/A$.
  Then any map $X\to Y$ in $\E/A$ is in particular a map $X\to Y$ in \E, hence factors through some $Y_i$, and the factorization lies in $\E/A$.
  Similarly, any two maps $X\to Y_i$ and $X\to Y_j$ in $\E/A$ are in particular maps in \E, hence coincide in some $Y_k$, and so the same is true in $\E/A$.
  Thus, $X\to A$ is \la-presentable in $\E/A$.

  Conversely, suppose $X\to A$ is \la-presentable in $\E/A$, and write $X = \colim_i X_i$ as a \la-filtered colimit of \la-presentable objects in \E.
  Then we can make each $X_i$ into an object of $\E/A$ by the composite $X_i \to X\to A$, and the colimit $X = \colim_i X_i$ then lies in $\E/A$ too.
  Thus, since $X\to A$ is \la-presentable, the identity $\id_X : X\to X$ factors through some $X_i$ in $\E/A$ and hence also in \E.
  So $X$ is a retract of some $X_i$, and is therefore \la-presentable like it.
\end{proof}

\begin{lem}\label{thm:relpres-pow}
  If \E is locally presentable and enriched with powers and copowers over some category \sV, then for any $K\in \sV$ there is a \la such that for any $\ka\shgt\la$ and $Y\in \E$ the power functor $(K \pow_Y -):\E/Y \to \E/Y$ preserves relatively \ka-presentable morphisms.
\end{lem}
\begin{proof}
  Let \E be locally \mu-presentable; by \cref{thm:relpres-detect}, it suffices to ensure that $Z\times_Y (K\pow_Y X)$ is \ka-presentable for any \mu-presentable $Z$ and relatively \ka-presentable $X\to Y$.
  But $Z\times_Y (K\pow_Y X) \cong K\pow_Z (Z\times_Y X)$, so for this it suffices to ensure that $(K\pow_Z -):\E/Z\to \E/Z$ preserves objects with \ka-presentable domain for all \mu-presentable $Z$.
  But by \cref{thm:pres-slice}, this is the same as preserving \ka-presentable objects of $\E/Z$.
  And each functor $(K\pow_Z -)$ is accessible (being a right adjoint), so there is a $\la_Z$ such that it preserves \ka-presentable objects for any $\ka\shgt\la_Z$.
  Finally, there is only a set of isomorphism classes of \mu-presentable objects $Z$, so there is a $\la$ with $\la\shgt\la_Z$ for all such $Z$.
\end{proof}

\begin{lem}\label{thm:relpres-mono-dp}
  Let \E be locally presentable and locally cartesian closed.
  Then there is a \la such that for any monomorphism $i:A\mono B$ and any $\ka\shgt\la$, the functor $i_* : \E/A \to \E/B$ preserves relatively \ka-presentable morphisms.
\end{lem}
\begin{proof}
  By \cref{thm:pres-sub}, there is a \mu such that \E is locally \mu-presentable and \mu-presentable objects are closed under subobjects.
  For each morphism $j$ between \mu-presentable objects, the functor $j_*$ is accessible, and there is only a set of such functors; thus there is a \la such that $\la\ge\mu$ and for any $\ka\shgt\la$ all these functors preserve \ka-presentable objects.

  Now let $i:A\mono B$ be any monomorphism and $f:X\to A$ be relatively \ka-presentable, where $\ka\shgt\la$.
  By \cref{thm:relpres-detect}, to show that $i_*(X)$ is relatively \ka-presentable it suffices to show that $Z\times_B i_*(X)$ is \ka-presentable for any morphism $Z\to B$ where $Z$ is \mu-presentable.
  Let $Y$ be the pullback $Z\times_B A$; then $j:Y\mono Z$ is a monomorphism, so $Y$ is also \mu-presentable by our choice of \mu.
  And by the Beck-Chevalley condition, we have $Z\times_B i_*(X) \cong j_*(Y\times_A X)$.
  But $Y \times_A X$ is \ka-presentable since $X\to A$ is relatively \ka-presentable, while $j_*$ preserves \ka-presentable objects since $\ka\shgt\la$; thus $Z\times_B i_*(X)$ is \ka-presentable.
\end{proof}

\begin{lem}\label{thm:relpres-dp}
  Let \E be locally presentable and locally cartesian closed.
  Then there is a regular cardinal $\la$ such that for any inaccessible cardinal $\ka>\la$ and any relatively \ka-presentable morphisms $X \xto{g} Y \xto{f} Z$, the dependent product $f_*(X) \to Z$ is relatively \ka-presentable.
\end{lem}
\begin{proof}
  Let \la satisfy \cref{thm:pres-pb}; then the regular cardinals \mu such that the \mu-presentable objects of \E are closed under finite limits are then unbounded below any inaccessible $\ka>\la$.

  Let \ka be such an inaccessible and $X \xto{g} Y \xto{f} Z$ be relatively \ka-presentable; we must show that for any morphism $W\to Z$ with $W$ \ka-presentable, the pullback $W\times_Z f_*(X)$ is \ka-presentable.
  As in \cref{thm:relpres-mono-dp}, let $V = W\times_Z Y$ with projection $h:V\to W$; then $V$ and $V\times_Y X$ are \ka-presentable, while $W\times_Z f_*(X) \cong g_*(V\times_Y X)$.
  Thus, it suffices to prove that $g_*$ preserves \ka-presentable objects for any morphism $g:V\to W$ between \ka-presentable objects.

  Since \ka is a limit cardinal, any \ka-presentable object is \mu-presentable for some $\mu<\ka$.
  Indeed, every \ka-presentable object is a \ka-small colimit of \la-presentable objects; but this diagram has some cardinality $<\ka$, hence is \mu-small for some $\mu<\ka$, and so its colimit is \mu-presentable.
  In particular, for any \ka-presentable object $Q$ and morphism $Q\to V$, there is a $\mu<\ka$ such that $Q$, $V$, and $W$ are all \mu-presentable and \mu-presentable objects are closed under pullback.
  Thus, $g^*:\E/W \to \E/V$ preserves \mu-presentable objects, hence by the proof of~\cite[Proposition 2.23]{ar:loc-pres} its right adjoint $g_*$ is \mu-accessible.
  Therefore, there is a $\nu\ge \mu$ with $\nu<\ka$ such that $g_*$ preserves \nu-presentable objects; hence $g_*(Q)$ is \nu-presentable and thus \ka-presentable.
\end{proof}

\begin{rmk}
  The asymmetry in hypotheses between \cref{thm:relpres-comp,thm:relpres-dp} (corresponding to $\Sigma$- and $\Pi$-types respectively) is due to our adherence to to the traditional use of only \emph{regular} cardinals to bound the size of objects in a locally presentable category.
  However, it should also be possible to study objects and morphisms bounded in size by singular cardinals (cf.~\cite{lrv:intsize}), enabling \cref{thm:relpres-dp} to apply to any strong limit \ka satisfying a $\shlt$-like property.
\end{rmk}

Finally, the following two facts show that relatively \ka-presentable morphisms serve our desired purpose.

\begin{prop}\label{thm:relpres-small}
  Let \E be locally \ka-presentable and locally cartesian closed.
  Then for any $Y\in\E$ the subcategory of $\E/Y$ determined by the relatively \ka-presentable morphisms is essentially small.
\end{prop}
\begin{proof}
  Write $Y = \colim_i Y_i$ as a colimit of \ka-presentable objects.
  Since colimits are stable under pullback, any $X\in \E/Y$ is the colimit in \E of the diagram of pullbacks $X \cong \colim_i (Y_i\times_Y X)$.
  Thus, if for $X$ and $X'$ the corresponding diagrams of pullbacks are isomorphic over $\{Y_i\}$, then $X\cong X'$.
  And if $X\to Y$ is relatively \ka-presentable, then each object $Y_i \times_Y X$ must be \ka-presentable, so the result follows since the full subcategory of \ka-presentable objects is essentially small.
\end{proof}

\begin{prop}[{cf.~\cite[6.1.6.5--6.1.6.7]{lurie:higher-topoi}}]\label{thm:relpres-local}
  Let \E be locally presentable and locally cartesian closed.
  Then there is a \la such that for any $\ka\shgt\la$, the relatively \ka-presentable morphisms are a \local full \nfs.
\end{prop}
\begin{proof}
  Let \la satisfy \cref{thm:pres-pb}, and $\ka\shgt\la$.
  By \cref{thm:local}, we must show that if $Y = \colim_{i\in I} Y_i$, then $f:X\to Y$ is relatively \ka-presentable as soon as its pullback to each $Y_i$ is.
  Since every colimit is a \ka-filtered colimit of \ka-small colimits, it will suffice to consider those two cases separately.

  When $I$ is \ka-filtered, any morphism $Z\to Y$ where $Z$ is \ka-presentable must factor through some $Y_i$.
  Then $Z\times_Y X$ is isomorphic to $Z\times_{Y_i} (Y_i \times_Y X)$, which is \ka-presentable since $Y_i \times_Y X \to Y_i$ is relatively \ka-presentable.

  When $I$ is \ka-small, note that by~\cite[Corollary 1.54]{ar:loc-pres}, $\E^I$ is also locally \ka-presentable.
  Thus we can write $\{Y_i\}\in\E^I$ as a \ka-filtered colimit of \ka-presentable diagrams, $Y_i = \colim_{j\in J} W_{j i}$ where $J$ is \ka-filtered and each $W_j\in \E^I$ is \ka-presentable.
  By \cref{thm:pw-pres}, each $W_{j i}\in \E$ is also \ka-presentable.
  But by assumption, the pullback of $f$ to each $Y_i$ is relatively \ka-presentable, and hence each $W_{j i} \times_Y X$ is \ka-presentable.

  Since $Y \cong \colim_i \colim_j W_{j i} \cong \colim_j \colim_i W_{j i}$, and $J$ is \ka-filtered, it will suffice to show that the pullback of $X$ to each $\colim_i W_{j i}$ is relatively \ka-presentable.
  But $\colim_i W_{j i}$ is a \ka-small colimit of \ka-presentable objects, hence \ka-presentable, so since \ka-presentable objects are closed under pullbacks it will suffice to show that the object $(\colim_i W_{j i}) \times_Y X$ is \ka-presentable.
  Finally, since \E is locally cartesian closed, colimits are stable under pullback, so this object is isomorphic to $\colim_i (W_{j i}\times_Y X)$, which is a \ka-small colimit of \ka-presentable objects, hence also \ka-presentable.
\end{proof}

We will denote this \nfs by $\cEka$.
More generally, for any \nfs \F we write $\Fka = \F \times_\cE \cEka$.


%% file: universes.tex
\section{Universes in model categories}
\label{sec:univalence}

Now let \E be a model category and \Fib the full \nfs determined by its fibrations, so that $\Fibka = \Fib\times_\cE \cEka$ denotes the relatively \ka-presentable fibrations.
In an ideal world, the object $\Fibka\in\Ehat$ would be (pseudonaturally equivalent to) a representable presheaf $\E(-,U)$.
But since \E is itself a 1-category rather than an \io-category, this is unreasonable to expect. 

Instead, we will replace \Fibka by a representable presheaf that is ``weakly equivalent'' in some sense.
We do not have a model structure on $\Ehat$ with which to make sense of this, but we can at least use the Yoneda embedding $\E\to\Ehat$ to lift the weak factorization systems of \E.

\begin{defn}\label{defn:afib}
  Let \E be a model category.
  A morphism $\dX\to\dY$ in $\Ehat$ is an \textbf{acyclic fibration} if it has the right lifting property (\cref{defn:2liftorth}) for all morphisms $\E(-,j) : \E(-,A) \to \E(-,B)$, where $j:A\to B$ is a cofibration in \E.
\end{defn}

\begin{rmk}
  If $f:\dX\to\dY$ is a representable morphism, then it is an acyclic fibration in the sense of \cref{defn:afib} if and only if in any pullback
  \begin{equation*}
    \begin{tikzcd}
      \E(-,W) \ar[r] \ar[d] \ar[dr,phantom,near start,"\lrcorner"] & \dX\ar[d,"{f}"] \\
      \E(-,Z) \ar[r] & \dY
    \end{tikzcd}
  \end{equation*}
  the induced map $W\to Z$ is an acyclic fibration in \E.
  This fits a standard pattern for extending pullback-stable properties of morphisms in \E to properties of representable morphisms in \Ehat.
  
  Note that this notion of (acyclic) fibration based on the model structure of \E is unrelated to the 2-categorical notions of (strict, discrete) fibration defined in \cref{sec:2cat}.
\end{rmk}

\begin{defn}
  If \F is a \nfs on \E, a \textbf{universe} for \F is a cofibrant object $U\in\E$ equipped with an acyclic fibration $\E(-,U) \to \F$ in \Ehat.
\end{defn}

That is, a universe is a sort of ``cofibrant replacement'' of $\Fibka$.
(This perspective was introduced informally in~\cite[\sect 3]{shulman:elreedy}.)

\begin{rmk}\label{rmk:universe}
  When $U$ is a universe, the morphism $\E(-,U) \to \F$ corresponds by the Yoneda lemma to an \F-algebra $\pi:\Util\to U$.
  And the fact that the morphism $\E(-,U) \to \F$ is an acyclic fibration means that given the solid arrows below, where $i:A\mono B$ is a cofibration, $f:X\to B$ is an \F-algebra, and both squares of solid arrows are pullbacks and \F-morphisms:
  \begin{equation*}
    \begin{tikzcd}[row sep=small, column sep=small]
      i^*(X) \arrow[rr] \arrow[dd, "g"'] \arrow[rd] &  & \Util \arrow[dd, "\pi"] \\
      & X \arrow[ru, dashed] &  \\
      A \arrow[rd, "i"', tail] \arrow[rr, near start, "h" description] &  & U \\
      & B \arrow[ru, dashed] \arrow[from=uu, near start, "f" description, crossing over] & 
    \end{tikzcd}
  \end{equation*}
  there exist the dashed arrows rendering the diagram commutative and the third square also a pullback and an \F-morphism.
  This property of a universe was first noted in the proof of~\cite[Theorem 2.2.1]{klv:ssetmodel} and isolated more abstractly in~\cite[(2$'$)]{shulman:elreedy}, \cite[Corollary 3.11]{cisinski:elegant}, and~\cite{stenzel:thesis} under varying names.

  If the initial object $\emptyset$ is strict (i.e.\ every morphism with codomain $\emptyset$ is an isomorphism) and $\id_\emptyset$ has a unique \F-structure, this implies that every \F-algebra with cofibrant codomain is a pullback of $\pi$ (though not in a unique way).
\end{rmk}

As usual, when \E is cofibrantly generated we can hope to produce such a cofibrant replacement by the small object argument.
However, the colimits in \E used to build cell complexes are no longer colimits in \Ehat; thus we have to restrict to the objects of \Ehat that preserve these particular colimits.

\begin{defn}
  Let \E be a model category.
  We say $\dX\in\Ehat$ is a \textbf{stack for cell complexes} if as a pseudofunctor $\dX:\E\op\to\cGPD$ it preserves (in the weak bicategorical sense) coproducts, pushouts of cofibrations, and transfinite composites of cofibrations.
\end{defn}

\begin{eg}\label{eg:topos-stack}
  If \E is a Grothendieck 1-topos and all cofibrations are monomorphisms, then the trivial \nfs \cE is a stack for cell complexes.
  This because any topos is infinitary extensive~\cite{clw:ext-dist}, adhesive~\cite{ls:adhesive,ls:topadh}, and exhaustive~\cite{nlab:exhaustive,shulman:elreedy}; see~\cite[\sect3]{shulman:elreedy} and~\cite[Lemma 7.5]{sattler:eqvext}.
  More generally,
  \cE being a stack for cell complexes is one of the conditions for the (cofibration, acyclic fibration) weak factorization system of \E to be \emph{suitable} as in~\cite[Definition 3.2]{sattler:eqvext}.
\end{eg}

\begin{lem}\label{thm:nfs-stack}
  If \cE is a stack for cell complexes and $\phi:\F\to\cE$ is a \local \nfs, then \F is also a stack for cell complexes.
\end{lem}
\begin{proof}
  Let \sQ be the class of morphisms $q:\Yhat \to \E(-,Y)$ from \cref{eg:colim-orth} where $Y= \colim_i Y_i$ ranges over coproducts, pushouts of cofibrations, and transfinite composites of cofibrations.
  Since \cE is a stack for cell complexes, $\sQ \perp \cE$; and since \F is \local, by \cref{thm:local} we have $\sQ\perp\phi$.
  Hence $\sQ\perp \F$.
\end{proof}

\begin{lem}\label{thm:icell-acyc}
  If \E is cofibrantly generated with \cI a set of generating cofibrations, \dX and \dY are stacks for cell complexes, and $f:\dX\to\dY$ has the right lifting property in \Ehat against all morphisms $\E(-,j):\E(-,A)\to\E(-,B)$ where $j:A\to B$ is in \cI, then $f$ is an acyclic fibration.
\end{lem}
\begin{proof}
  As usual, any cofibration is a retract of a transfinite composite of pushouts of coproducts of elements of \cI.
  Thus, given that \dX and \dY are stacks for cell complexes, the lifting property carries through all these operations in the usual way.
\end{proof}

We call a pseudofunctor $\dZ\in\Ehat$ \textbf{small-groupoid-valued} if each groupoid $\dZ(A)$ is essentially small.
Note that by definition, for any \nfs \F the map $\phi:\F\to\cE$ has small \emph{fibers}, i.e.\ any given morphism $f:X\to Y$ has a small set of \F-structures; but \F is only small-groupoid-valued if any given object $Y\in\E$ there is a small set of isomorphism classes of \F-algebras with codomain $Y$.
In general we will achieve this by considering $\Fka = \dF\times_\cE \cEka$ as in \cref{sec:relpres}.

\begin{thm}\label{thm:2cat-soa}
  Let \E be a combinatorial model category, and $\dZ\in\Ehat$ a small-groupoid-valued stack for cell complexes.
  Then any morphism $f:\E(-,X) \to\dZ$ in \Ehat factors, up to isomorphism, as $\E(-,X) \xto{\E(-,j)} \E(-,Y) \xto{p} \dZ$, where $j$ is a cofibration in \E and $p$ is an acyclic fibration in \Ehat.
\end{thm}
\begin{proof}
  This is just a bicategorical adaptation of the small object argument.
  Let \cI be a set of generating cofibrations for \E; we will define an \cI-cell complex sequence $X_0 \to X_1 \to \cdots$ in \E, along with maps $f_n : \E(-,X_n) \to \dZ$ and coherent isomorphisms $f_n \circ j_{m,n} \cong f_m$.
  We start with $X_0 = X$ and $f_0 = f$.
  For limit $n$ we let $X_n = \colim_{m<n} X_m$, with $f_n : \E(-,X_n) \to \dZ$ and attendant isomorphisms induced by the fact that \dZ preserves this colimit.

  At a successor stage $n+1$, we let $S_n$ be a set of representatives for isomorphism classes of pseudo-commutative squares
  \[
    \begin{tikzcd}
      \E(-,A) \ar[r] \ar[d,"i"'] \ar[dr,phantom,"\scriptstyle\Downarrow\cong"] & \E(-,X_n) \ar[d] \\
      \E(-,B) \ar[r] & \dZ
    \end{tikzcd}
  \]
  where $i\in \cI$.
  This is a small set, since \dZ is small-groupoid-valued and \E is locally small.
  Now let $X_{n+1}$ be the pushout
  \[
    \begin{tikzcd}
      \coprod_{s\in S_n} A_s \ar[r] \ar[d] \drpushout & X_n \ar[d]\\
      \coprod_{s\in S_n} B_s \ar[r] & X_{n+1}
    \end{tikzcd}
  \]
  Since \dX preserves these coproducts and pushouts, there is an essentially unique induced map $f_{n+1} : X_{n+1}\to Z$ with attendant isomorphisms.

  Finally, since \E is locally presentable, there is a regular cardinal \la such that all domains of morphisms in \cI are \la-presentable.
  Thus, in any square
  \[
    \begin{tikzcd}
      \E(-,A) \ar[r] \ar[d,"i"'] \ar[dr,phantom,"\scriptstyle\Downarrow\cong"] & \E(-,X_\la) \ar[d] \\
      \E(-,B) \ar[r] & \dZ
    \end{tikzcd}
  \]
  the top morphism $A\to X_\la$ factors through $X_n$ for some $n<\la$, and hence there is a lift $B\to X_{n+1} \to X_\la$.
  Therefore, the map $f_{\la} : X_\la \to \dZ$ has right lifting for \cI, and is thus an acyclic fibration by \cref{thm:icell-acyc}.
\end{proof}

\begin{cor}\label{thm:nfs-universe}
  If \E is a Grothendieck 1-topos with a combinatorial model structure in which all cofibrations are monomorphisms, then any small-groupoid-valued \local \nfs \F on \E has a universe.
\end{cor}
\begin{proof}
  By \cref{eg:topos-stack,thm:nfs-stack}, \F is a stack for cell complexes; thus we can apply \cref{thm:2cat-soa} to factor the map $\E(-,\emptyset) \to \F$.
\end{proof}

In some cases such as \cref{eg:pshf-can,eg:rep-cod}, \Fib is \local and hence so is \Fibka.
This includes the universes constructed in~\cite{klv:ssetmodel,shulman:elreedy,cisinski:elegant}.
However, in the general case we need a different approach: we will suppose given a non-full \nfs \F that \emph{is} \local, and an acyclic fibration $\F\to\Fib$. Thus we will be able to apply \cref{thm:nfs-universe} to \F instead.

More generally, for a \nfs \F, let $\uly\F$ denote the image of the map $\phi:\F\to \cE$.
Thus $\uly\F$ is a full \nfs (though not generally \local, even if \F is), and the $\uly\F$-algebras are the morphisms that admit some \F-structure.

\begin{defn}\label{defn:stratified}
  A \nfs \F on a model category \E is \textbf{\stratified} if the map $\F \to \uly\F$ is an acyclic fibration.
  That is, for any pullback
  \[
    \begin{tikzcd}
      X' \ar[d,"f'"'] \ar[r,"g"] \ar[dr,phantom, near start,"\lrcorner"] & X \ar[d,"f"]\\
      Y' \ar[r,"i",tail] & Y
    \end{tikzcd}
  \]
  with $f$ and $f'$ \F-algebras and $i$ a cofibration, there exists a new \F-structure on $f$ making the square an \F-morphism.
\end{defn}

\begin{prop}\label{thm:pre-u2p}
  Let \F be a \local, \stratified, small-groupoid-valued \nfs on a Grothendieck 1-topos that is a combinatorial model category whose cofibrations are monomorphisms.
  Then $\uly\F$ has a universe.
\end{prop}
\begin{proof}
  Apply \cref{thm:nfs-universe} to \F, and observe that the composite $\E(-,U) \to \F \to \uly\F$ of acyclic fibrations is again an acyclic fibration.
\end{proof}

\begin{rmk}
  By \cref{rmk:universe}, if all objects are cofibrant, $\emptyset$ is strict, and $U$ is a universe for a full \nfs \F such that $\id_\emptyset\in\F$, then in fact $\F = \uly{\dRep_\pi}$ (with $\dRep_\pi$ as in \cref{eg:rep-fcos}).
  Conversely, if \E is locally cartesian closed and $\pi:\Util\to U$ is a universe for $\uly{\dRep_\pi}$, then $\dRep_\pi$ is \local (by \cref{eg:rep-local}), \stratified, and small-groupoid-valued.
  Thus the hypotheses of \cref{thm:pre-u2p} are basically optimal.
\end{rmk}

\begin{eg}
  Full \nfss are always \stratified, as is $\F_1 \times_\cE \F_2$ if $\F_1$ and $\F_2$ are.
\end{eg}

\begin{eg}
  If \F is a \stratified \nfs on $\E_2$ and $G:\E_1\to\E_2$ preserves pullbacks (hence also monomorphisms), then $G^{-1}(\F)$ is also \stratified.
\end{eg}

\begin{eg}\label{thm:sec-afib-strat}
  Let \E be a model category and $H$ be a fibred core-endofunctor of \E such that whenever $H_Y(X)\to Y$ has a section, it is an acyclic fibration.
  Then the \nfs $H^*(\cEp)$ is \stratified.
  For given a pullback square
  \begin{equation}
    \begin{tikzcd}
      X' \ar[d,"f'"'] \ar[r,"g"] \ar[dr,phantom, near start,"\lrcorner"] & X \ar[d,"f"]\\
      Y' \ar[r,"i",tail] & Y
    \end{tikzcd}\label{eq:sec-strat-sq}
  \end{equation}
  with $i$ a cofibration, along with sections $s'$ and $s$ of $H_{Y'}(X')\to Y'$ and $H_Y(X)\to Y$, the assumption implies $H_Y(X)\to Y$ is an acyclic fibration.
  Thus we can find a lift in the square
  \[
    \begin{tikzcd}
      Y' \ar[d,"i"',tail] \ar[r,"s"] & H_{Y'}(X') \ar[r,"{H(g,i)}"] & H_Y(X) \ar[d,"\sim",two heads]\\
      Y \ar[rr,equals] & & Y
    \end{tikzcd}
  \]
  giving an $H^*(\cEp)$-structure on $f$ making~\eqref{eq:sec-strat-sq} an $H^*(\cEp)$-morphism.
\end{eg}

Recall that a \textbf{Cisinski model category}~\cite{cisinski:topos,cisinski:local-acyc} is a Grothendieck 1-topos with a combinatorial model structure whose cofibrations are \emph{precisely} the monomorphisms.

\begin{prop}\label{eg:cof-ff-fcos}
  Let \E be a Cisinski model category, \F a \stratified \nfs on \E, and $E$ a cartesian functorial factorization on \E that factors every \F-algebra as an acyclic cofibration followed by a fibration.
  Then the \nfs $\F\times_\cE \dR_E$ (where $\dR_E$ is as in \cref{eg:ff-fcos}) is also \stratified.
\end{prop}
\begin{proof}
  Suppose given the pullback square on the left:
  \[
    \begin{tikzcd}
      X' \ar[d,"f'"'] \ar[r,"g"] \ar[dr,phantom, near start,"\lrcorner"] & X \ar[d,"f"]\\
      Y' \ar[r,"i",tail] & Y\\
      \phantom{E f}
    \end{tikzcd}
    \hspace{2cm}
    \begin{tikzcd}
      X' \ar[r,"g"] \ar[d,tail,"{\lambda_{f'}}"'] \ar[dr,phantom,near end,"\ulcorner"] & X \ar[ddr,tail,"{\lambda_f}"] \ar[d] \\
      E f' \ar[r] \ar[drr,tail,"{\fact g i}"'] & P \ar[dr,"j" description]\\
      && E f.
    \end{tikzcd}
  \]
  where $f$ and $f'$ are \F-algebras with $\dR_E$-structures $r_{f} : E f \to X$ and $r_{f'} : E f' \to X'$.
  Since \F is \stratified, $f$ has a new \F-structure making $(g,i)$ an \F-morphism; so it remains to find a new $\dR_E$-structure $\rtil_{f} : E f \to X$ (so that $f \circ \rtil_f = \rho_f$ and $\rtil_f \circ \lambda_f = \id_X$) such that $(g,i)$ is also an $\dR_E$-morphism, i.e.\ $\rtil_f \circ \fact g i = g \circ r_{f'}$.

  Define $P$ and $j$ by the pushout as on the right above.
  Since $f$ and $f'$ are \F-algebras, $\lambda_{f}$ and $\lambda_{f'}$ are acyclic cofibrations, and in particular monomorphisms.
  By cartesianness, $\fact g i$ is also a monomorphism (being a pullback of $i$) and $X'\cong E f' \times_{E f}X$.
  Thus, $P$ is a union of subobjects of $E f$ in the 1-topos \E, hence $j:P\to E f$ is also a monomorphism.
  Moreover, since $X\to P$ is a pushout of $\lambda_{f'}$, it is also an acyclic cofibration; hence by 2-out-of-3 $j$ is also acyclic.

  Now since $f$ is an \F-algebra, $\rho_f$ is a fibration.
  But $f$ is a retract of $\rho_f$ (by its $\dR_E$-structure $r_f$), hence also a fibration.
  Thus we can find a lift in the square:
  \[
    \begin{tikzcd}
      P \ar[d,"j"'] \ar[r] & X \ar[d,"f"] \\
      E f \ar[r,"{\rho_f}"'] \ar[ur,dotted,"{\rtil_f}" description] & Y
    \end{tikzcd}
  \]
  where the top arrow is induced by $g \circ r_{f'}$ and $\id_X$.
  Such a lift is then an $\dR_E$-structure on $f$ such that $(g,i)$ is an $\dR_E$-morphism, as desired.
\end{proof}

It remains to ensure that our universes are fibrant and univalent.

\begin{defn}\label{defn:nfs-hoinvar}
  Let \E be a locally presentable category with a model structure.
  A \nfs \F is \textbf{homotopy invariant} if every \F-algebra is a fibration, and given any commutative square
  \[
    \begin{tikzcd}
      X' \ar[d,"f'"',two heads] \ar[r,"\sim"] & X \ar[d,"f",two heads]\\
      Y' \ar[r,"\sim"'] & Y
    \end{tikzcd}
  \]
  where $f$ and $f'$ are fibrations and the horizontal maps are weak equivalences, $f$ admits an \F-structure if and only if $f'$ does.
\end{defn}

Of course, if $\uly\F=\Fib$ is the class of all fibrations, then \F is homotopy invariant.
More generally, homotopy invariance is a condition only on $\uly\F$.

We now recall the fundamental ``equivalence extension'' property.
To my knowledge, a version of this property first appeared in~\cite[Theorem 3.4.1]{klv:ssetmodel} in the case of simplicial sets.
It was observed in~\cite[Theorem 3.1]{shulman:elreedy} and~\cite[Remark 2.19]{cisinski:elegant} that the proof generalizes to any simplicial Cisinski model category.
A similar construction for cubical sets appeared under the name ``gluing'' in~\cite{cchm:cubicaltt}, which was then placed in a more abstract setting by~\cite{sattler:eqvext}.

\begin{thm}\label{thm:u4p}
  Let $\E$ be a simplicial Cisinski model category and \F a homotopy invariant \nfs on \E.
  Then there is a \la such that for any $\ka\shgt\la$ and any cofibration $i: A\cof B$, relatively \ka-presentable $\uly\F$-algebras $D_2 \fib B$ and $E_1 \fib A$, and weak equivalence $w: E_1 \toiso E_2$ over $A$, where $E_2 \coloneqq i^* D_2$:
  \begin{equation}
  \begin{tikzcd}[column sep=small]
      E_1 \arrow[rdd, two heads] \arrow[rrd, "w"] \arrow[rrr, dashed] &  &  & D_1 \arrow[rdd, two heads, dashed] \arrow[rrd, "v", dashed] &  &  \\
      &  & \mathllap{E_2=\;}i^* D_2 \arrow[ld, two heads] \arrow[rrr,crossing over]
      &  &  & D_2 \arrow[ld, two heads] \\
      & A \arrow[rrr, tail,"i"'] &  &  & B & 
    \end{tikzcd}\label{eq:u4p}
  \end{equation}
  there exists a relatively \ka-presentable $\uly\F$-algebra $D_1\fib B$ and an equivalence $v: D_1 \toiso D_2$ over $B$ such that $i^*(v)=w$.
\end{thm}
\begin{proof}
  Largely identical to that of~\cite[Theorem 3.1]{shulman:elreedy}.
  The latter statement assumes that \E is a presheaf category, but this is only used to obtain a notion of ``\ka-small morphism'' that is preserved by $i_*$ and by pullback; using \cref{thm:pres-pb,thm:relpres-mono-dp} instead allows \E to be any Grothendieck 1-topos.\footnote{The author claimed in~\cite{shulman:elreedy} that when $\E=\prcs$ it suffices to take $\ka > \card\C$, but Raffael Stenzel has pointed out that this is not enough to ensure that $i_*$ preserves \ka-small morphisms; even in the presheaf case we need some analogue of the relation $\shlt$.}
  (The proof uses that \E is a simplicial model category and has effective unions, to extend deformation retractions along $i$.)
  We conclude $D_1$ is an $\uly\F$-algebra by homotopy invariance, since it is equivalent to $D_2$ over $B$.
\end{proof}

Univalence of our universes will follow from \cref{thm:u4p} as in~\cite{klv:ssetmodel,shulman:elreedy,cisinski:elegant}.
To show that $U$ is a fibrant object,~\cite[Theorem 2.2.1]{klv:ssetmodel} and~\cite[Proposition 2.21]{cisinski:elegant} use minimal fibrations, while~\cite[Lemma 6.3]{shulman:elreedy} uses a Reedy induction; but in fact fibrancy of $U$ is almost immediate from \cref{thm:u4p}.
A similar fact in the restricted situation of cubical-type model structures (where an explicit description of the fibrations is available) appears in~\cite{cchm:cubicaltt,sattler:eqvext}, while the general case was observed in~\cite{stenzel:thesis}.

\begin{thm}\label{thm:uf-fibrant}
  Let $\E$ be a right proper simplicial Cisinski model category, and \F a \local, \stratified, and homotopy invariant \nfs on \E.
  Then there is a regular cardinal \la such that for any regular cardinal $\ka\shgt\la$, there exists a morphism $\pi:\Util \to U$ such that:
  \begin{enumerate}
  \item The \ka-presentable objects in \E are closed under finite limits.\label{item:uf0}
  \item $\pi:\Util \to U$ is a relatively \ka-presentable $\uly\F$-algebra (in particular, a fibration).\label{item:uf1}
  \item Every relatively \ka-presentable $\uly\F$-algebra is a pullback of $\pi$.\label{item:uf2}
  \item The object $U$ is fibrant.\label{item:uf4}
  \item $\pi$ satisfies the univalence axiom.\label{item:uf3}
  \end{enumerate}
\end{thm}
\begin{proof}
  Let $\la_0$ satisfy \cref{thm:u4p}, let $\la_1$ be such that \E has a generating set of acyclic cofibrations with $\la_1$-presentable domains and codomains, let $\la_2$ be such that \E has functorial factorizations that preserve \ka-presentable objects for any $\ka\shgt\la_2$ (such exists since these factorizations are accessible functors), and let $\la_3$ be such that for any $\ka\shgt\la_3$ the \ka-presentable objects are closed under finite limits (which exists by \cref{thm:pres-pb}).
  Let $\la$ be such that $\la\shgt \la_j$ for $j=0,1,2,3$, and assume $\ka\shgt \la$; then $\ka\shgt\la_j$ for all $j$ as well, and in particular~\ref{item:uf0} holds.

  Since \F and $\cEka$ are \local and \stratified, so is $\Fka = {\F\times_\cE \cEka}$.
  Let $\pi:\Util\to U$ be the universe for $\uly\Fka$ obtained from \cref{thm:pre-u2p}; then~\ref{item:uf1} holds trivially.
  And since \cE and \F are stacks for cell complexes, in particular they preserve the initial object; so by \cref{rmk:universe} we have~\ref{item:uf2}.

  Since $\ka\gt\la_1$, to show that $U$ is fibrant~\ref{item:uf4} it suffices to show that it has right lifting for all acyclic cofibrations between \ka-presentable objects.
  Let $i:A\acof B$ be an acyclic cofibration with $A$ and $B$ \ka-presentable, let $h:A\to U$ be a map, and let $E_1 \fib A$ be the pullback of $\pi$ along $h$.
  Since $\pi$ is relatively \ka-presentable, $E_1$ is \ka-presentable.
  Thus since $\ka\shgt\la_2$, we can factor the composite $E_1\fib A \xto{i} B$ as an acyclic cofibration $E_1 \acof D_2$ followed by a fibration $D_2 \fib B$, where $D_2$ is \ka-presentable.
  So since $\ka\shgt\la_3$, $D_2\fib B$ is a relatively \ka-presentable fibration, and by homotopy invariance it is an $\uly\F$-algebra, hence an $\uly\Fka$-algebra.

  Let $E_2\coloneqq i^*(D_2)$; then by right properness the map $E_2 \to D_2$ is a weak equivalence, hence by 2-out-of-3 so is the induced map $E_1 \to E_2$.
  Thus since $\ka\shgt\la_0$, by \cref{thm:u4p} there is an $\uly\Fka$-algebra $D_1 \fib B$ with $i^*(D_1)\cong E_1$.
  Finally, since $U$ is a universe for $\uly\Fka$, by \cref{rmk:universe} there is a map $k:B\to U$ pulling $\pi$ back to $D_1$ such that $k i = h$; so $U$ has right lifting for $i$.

  For univalence~\ref{item:uf3}, we follow~\cite[Theorem 3.4.1]{klv:ssetmodel},~\cite[\sect 2]{shulman:elreedy}, and~\cite[Theorem 3.12]{cisinski:elegant}.
  Let $\Eq(\Util)$ be the universal space of auto-equivalences of $\pi$, as in~\cite[\sect 4]{shulman:elreedy}; it suffices to show that the composite projection
  \(\Eq(\Util) \to U\times U \to U\)
  is an acyclic fibration.
  Now a square
  \[
    \begin{tikzcd}
      A \ar[d,tail,"i"'] \ar[r] & \Eq(\Util) \ar[d]\\
      B \ar[r] & U
    \end{tikzcd}
  \]
  with $i$ a monomorphism yields a diagram of solid arrows~\eqref{eq:u4p} where all fibrations are $\uly\Fka$-algebras.
  Thus, since $\ka\shgt\la_0$ we can fill out the dashed arrows in~\eqref{eq:u4p} with $D_1 \fib B$ also a $\uly\Fka$-algebra; so by \cref{rmk:universe} we can classify it by a map to $U$ extending the given classifying map of $E_1$.
  But this is precisely what we need to specify a lift $B\to \Eq(\Util)$.
\end{proof}

Thus, to build fibrant and univalent universes for relatively \ka-presentable fibrations in a right proper simplicial Cisinski model category, it suffices to find a \local and \stratified \nfs \F such that $\uly\F=\Fib$ is the class of all fibrations.
We have essentially already seen one way to do this: if \E has a set of generating acyclic cofibrations with representable codomains, then $\Fib$ itself has these properties.
This was the approach of~\cite{klv:ssetmodel,shulman:elreedy}; but to deal with the general case we will have to use non-full \nfss.

\begin{rmk}
  Although our primary interest is in constructing univalent universes for \emph{all} (relatively \ka-presentable) fibrations, it is potentially useful that \cref{thm:uf-fibrant} also yields univalent universes for subclasses of fibrations.
  For instance, the \emph{left fibrations} in bisimplicial sets~\cite{vk:yoneda-css,pbb:groth-segal,rasekh:yoneda-ss} are a subclass of the Reedy fibrations, which by~\cite[Remark 2.1.4(a)]{vk:yoneda-css} are characterized by right lifting against a generating set with representable codomains; thus they admit fibrant and univalent universes.
  Such a universe of left fibrations is essentially the ``\oo-category of spaces'' constructed in~\cite{vk:yoneda-css}, although they do not explain how to make it a strict presheaf.
  In fact it is a complete Segal space (this is shown in~\cite[Theorem 2.2.11]{vk:yoneda-css}, and can also be deduced from~\cite[Theorem 4.8]{rasekh:yoneda-ss}), and could be useful for the programme of~\cite{rs:stt} to use that model for ``synthetic \io-category theory''.
  
  We will see another class of examples in \cref{thm:flf-nfs,rmk:modal-univ}.
\end{rmk}


%% file: ttmt.tex
\section{Type-theoretic model toposes}
\label{sec:ttmt}

Numerous authors have defined classes of model categories that are well-adapted to type theory; \cref{thm:uf-fibrant} suggests a particularly strong one.

\begin{defn}\label{defn:ttmt}
  A \textbf{\ttmt} is a model category $\E$ such that:
  \begin{enumerate}
  \item \E is a Grothendieck 1-topos.
  \item The model structure is right proper, 
    simplicial, and combinatorial, and its cofibrations are the monomorphisms (hence it is also left proper).
    That is, \E is a right proper simplicial Cisinsiki model category.
  \item \E is \slcc.\footnote{That is, each pullback functor $f^* :\E/Y \to \E/X$ has a \emph{simplicially enriched} right adjoint.
      Since it is always a simplicial functor and has an ordinary right adjoint (since toposes are locally cartesian closed), this is equivalent to its preserving simplicial copowers.}
  \item There is a \local and \stratified \nfs \F on \E such that $\uly\F$ is the class of all fibrations.
  \end{enumerate}
\end{defn}

By \cref{thm:uf-fibrant}, for any \ttmt \E there is a \la such that \E has a fibrant univalent universe of relatively \ka-presentable fibrations for any $\ka\shgt\la$.
It also has all the other requisite structure to model type theory:
\begin{itemize}
\item 
  \E is a \emph{logical model category} in the sense of~\cite{ak:htmtt}, hence models $\Sigma$- and $\Pi$-types (and also a unit type and identity types, although these were not discussed in~\cite{ak:htmtt}).
  Categorically, $\Sigma$-types correspond to composition of fibrations, the unit type corresponds to the identity map as a fibration, $\Pi$-types correspond to dependent product of one fibration along another, and identity types correspond to path objects as in~\cite{aw:htpy-idtype}.

  One also needs a coherence theorem to strictify \E into an actual model of type theory.
  With inaccessible universes, we can use the method of~\cite{klv:ssetmodel}; otherwise we can use the ``local universes'' technique of~\cite{lw:localuniv} (see also~\cite{awodey:natmodels}).
\item 
  \E is a \emph{type-theoretic model category} in the sense of~\cite{shulman:invdia}, hence in particular its $\Pi$-types satisfy function extensionality.
  Categorically, function extensionality means that for any fibration $f$, the adjunction $f^*\adj f_*$ is a Quillen adjunction, with $f_*$ preserving both fibrations and acyclic fibrations.
\item \E is a \emph{good model category} as in~\cite{ls:hits}, hence models higher inductive pushouts (obtained as fibrant replacements of explicit simplicial homotopy pushouts) and other non-recursive (higher) inductive types such as the empty type, the boolean type, coproduct types, circles, spheres, tori, and so on.
\item \E is also an \emph{excellent model category}\footnote{No relation to the ``excellent model categories'' of~\cite[A.3.2.16]{lurie:higher-topoi}.} in the sense of~\cite{ls:hits}, hence models many other higher inductive types, including the natural numbers, $W$-types, truncations, and localizations.
  These are obtained by mixing a fibrant replacement monad with a ``cell monad'' built from polynomial endofunctors.
  Strictification for these types is discussed in~\cite{ls:hits} using the method of~\cite{lw:localuniv}, but adapts easily to the universe method of~\cite{klv:ssetmodel}.
\end{itemize}

We also expect universes in type theory to be \emph{closed} under all the relevant type constructors.
Since our universes classify all relatively \ka-presentable fibrations, as in~\cite{klv:ssetmodel} this will be true if the corresponding categorical operations preserve relatively \ka-presentable fibrations.
\begin{itemize}
\item By \cref{thm:relpres-comp}, $\Sigma$-types preserve relatively \ka-presentable morphisms for any regular cardinal \ka.
\item Likewise, identity maps (i.e.\ the unit type) are always relatively \ka-presentable.
\item If we define identity types as powers by $\Delta[1]$, then by \cref{thm:relpres-pow} there is a \la such that they preserve relatively \ka-presentable morphisms for any $\ka\shgt\la$.
\item By \cref{thm:relpres-dp}, $\Pi$-types preserve relatively \ka-presentable morphisms if $\ka$ is sufficiently large and inaccessible.
\item We can also choose \ka large enough that the universe will contain any fixed collection of (higher) inductive types, such as the empty type, the boolean type, the natural numbers, circles, spheres, tori, and so on.
\end{itemize}
It is not yet known how to obtain universes closed under parametrized (higher) inductive types such as $W$-types and pushouts, since fibrant replacement need not preserve relatively \ka-presentable morphisms.
However, in the special case of binary coproducts, once we have universes we can use the trick of defining $A+B = \sum_{x:\bool} \mathsf{rec}_\bool(U,A,B)$, where $\bool$ is the boolean type:\footnote{I am indebted to Bas Spitters for pointing this out.}

\begin{prop}\label{thm:relpres-coprod}
  In a \ttmt, there exists a regular cardinal \la such that for any $\ka\shgt\la$, if $X\fib Z$ and $Y\fib Z$ are relatively \ka-presentable fibrations with fibrant codomain, then their copairing $X\amalg Y \to Z$ factors as an acyclic cofibration $X\amalg Y \acof P$ followed by a fibration $P\fib Z$ that is again relatively \ka-presentable.
\end{prop}
\begin{proof}
  Let $\bool$ denote a fibrant replacement of $1\amalg 1$, so we have an acyclic cofibration $1\amalg 1\acof \bool$.
  We let \la be such that $\bool$ is \la-presentable and \E has a fibrant univalent universe of relatively \ka-presentable fibrations for any $\ka\shgt\la$.
  Let $\pi:\Util\fib U$ be such a universe, and let $x,y:Z \toto U$ be classifying maps for $X$ and $Y$ respectively.
  Then since $Z$ is fibrant, the map $Z\amalg Z \cong Z\times (1\amalg 1) \acof Z\times \bool$ is again an acyclic cofibration, so the map $[x,y] : Z\amalg Z \to U$ extends to a map $w:Z\times \bool \to U$.

  Let $P = w^* \Util$; then the composite $P \fib Z\times \bool \fib Z$ is a fibration (since $\bool$ is fibrant).
  And it is relatively \ka-presentable, since $\pi$ is relatively \ka-presentable, \bool is \ka-presentable, and relatively \ka-presentable morphisms are closed under pullback and composition.
  Moreover, in the following diagram:
  \[
    \begin{tikzcd}
      X\amalg Y \ar[d,two heads] \ar[r] \drpullback & P \ar[d,two heads] \ar[r]\drpullback & \Util \ar[d,two heads,"\pi"]\\
      Z\amalg Z \ar[r,tail,"\sim"'] &  Z\times \bool \ar[r,two heads,"w"'] & U
    \end{tikzcd}
  \]
  the outer rectangle and right-hand square are pullbacks (the former since \E is extensive), hence so is the left-hand square.
  Since $P\to Z\times \bool$ is a fibration, this implies that $X\amalg Y \to P$ is, like $Z\amalg Z\acof Z\times \bool$, an acyclic cofibration.
\end{proof}

In summary, we have the following.

\begin{thm}\label{thm:ttmt-models}
  For any \ttmt \E, there is a regular cardinal \la such that \E interprets Martin-L\"{o}f type theory with the following structure:
  \begin{enumerate}
  \item $\Sigma$-types, a unit type, $\Pi$-types with function extensionality, identity types, and binary sum types.\label{item:sigpiid}
  \item The empty type, the natural numbers type, the circle type $S^1$, the sphere types $S^n$, and other specific ``cell complex'' types such as the torus $T^2$.\label{item:unparam-hits}
  \item As many universe types as there are inaccessible cardinals larger than $\la$, all closed under the type formers~\ref{item:sigpiid} and containing the types~\ref{item:unparam-hits}, and satisfying the univalence axiom.\label{item:univ}
  \item $\mathsf{W}$-types, pushout types, truncations, localizations, James constructions, and many other recursive higher inductive types.\label{item:hits}
  \end{enumerate}
\end{thm}
\begin{proof}
  We have already noted that \E has all the structure in~\ref{item:sigpiid}, \ref{item:unparam-hits}, and~\ref{item:hits}.
  Let \la satisfy \cref{thm:uf-fibrant,thm:relpres-pow,thm:relpres-coprod} for \E, and also be such that the unparametrized higher inductive types in~\ref{item:unparam-hits} are \la-presentable.
  Then if $\ka>\la$ is inaccessible (hence in particular $\ka\shgt\la$), the universe for relatively \ka-presentable fibrations is univalent and closed under~\ref{item:sigpiid} by the above remarks.
  Moreover, since each type $X$ in~\ref{item:unparam-hits} is \la-presentable, hence \ka-presentable, and the \ka-presentable objects are closed under finite limits, the morphism $X\to 1$ is relatively \ka-presentable and hence can be classified by all these universes.

  To complete the proof it is necessary to apply a coherence theorem to replace \E by a strict model of type theory (such as a category with families or contextual category).
  Unfortunately, the coherence theorem of~\cite{klv:ssetmodel} deals with only one internal universe and requires an additional inaccessible outside the model, while that of~\cite{lw:localuniv} does not mention universes at all.
  Thus, in \cref{sec:coherence} we extend the latter to handle an arbitrary family of universes.
\end{proof}

A \ttmt is also a \emph{combinatorial type-theoretic model category} as in~\cite{gk:univlcc}, hence the \io-category it presents is locally presentable and locally cartesian closed.
In fact, it is a Grothendieck \io-topos, as we now show.

\begin{thm}\label{thm:descent}
  A \ttmt has descent in the sense of~\cite{rezk:homotopy-toposes}.
\end{thm}
\begin{proof}
  As defined in~\cite{rezk:homotopy-toposes}, descent consists of two conditions.
  Condition (P1) says that homotopy colimits are stable under homotopy pullback.
  Since \E is right proper, it suffices to show that for any fibration $f:X\fib Y$, the pullback functor $f^* : \E/Y \to \E/X$ preserves homotopy colimits.
  But it preserves cofibrations (as these are the monomorphisms) and weak equivalences (by right properness), and has a right adjoint (since toposes are locally cartesian closed).
  Thus it is a left Quillen functor, hence preserves homotopy colimits.

  Condition (P2) says that if $f:X\to Y$ is a map between homotopy colimits $X = \hocolim_i X_i$ and $Y = \hocolim_i Y_i$ induced by a natural transformation $f_i : X_i\to Y_i$ such that each square on the left below is a homotopy pullback (i.e.\ the transformation is ``cartesian'' or ``equifibered''):
  \begin{equation}
    \begin{tikzcd}
      X_i \ar[r] \ar[d,"f_i"'] & X_j \ar[d,"f_j"]\\
      Y_i \ar[r] & Y_j
    \end{tikzcd}
    \hspace{2cm}
    \begin{tikzcd}
      X_i \ar[r] \ar[d,"f_i"'] & X \ar[d,"f"]\\
      Y_i \ar[r] & Y,
    \end{tikzcd}\label{eq:descent}
  \end{equation}
  then the squares on the right above are also homotopy pullbacks.
  To prove this, it suffices to consider coproducts and pushouts. 

  For coproducts, we may assume given fibrations $f_i : X_i \fib Y_i$; equifiberedness is vacuous, and all coproducts are homotopy colimits since all objects are cofibrant.
  Choosing an \F-structure on each $f_i$, \locality of \F induces an \F-structure on $f$; hence it is a fibration.
  And the squares on the right in~\eqref{eq:descent} are strict pullbacks in \E since it is a 1-topos; thus they are also homotopy pullbacks.

  For pushouts, suppose given the solid arrows below, where the vertical maps $f_i:X_i\to Y_i$ are fibrations, and the maps $X_0\cof X_2$ and $Y_0\cof Y_2$ are cofibrations.
  \[ \begin{tikzcd}[row sep=small,column sep=small]
      X_0 \arrow[rr, tail] \arrow[dd, two heads,"f_0"'] \arrow[rd] &  & X_2 \arrow[dd, two heads,near end,"f_2"] \arrow[rd, dashed] &  \\
      & X_1 \arrow[rr, dashed, crossing over]  &  & X \arrow[dd, two heads, dashed,"f"] \\
      Y_0 \arrow[rr, tail] \arrow[rd] &  & Y_2 \arrow[rd, dashed] &  \\
      & Y_1 \arrow[rr, dashed] \arrow[from=uu,crossing over, two heads,near start,"f_1"] &  & Y
    \end{tikzcd}
  \]
  Equifiberedness means the two squares of solid arrows are homotopy pullbacks, i.e.\ the maps $X_0 \to X_1\times_{Y_1} Y_0$ and $X_0 \to X_2\times_{Y_2} Y_0$ are equivalences.
  Up to homotopy, we can therefore replace $X_0$ with $X_1\times_{Y_1} Y_0$ to make the left-hand face of the cube a strict pullback, and then by \cref{thm:u4p} we can replace $f_2$ by an equivalent fibration making the back face of the cube also a strict pullback.

  Now we take the strict pushouts, which are also homotopy pushouts since $X_0\cof X_2$ and $Y_0\cof Y_2$ are cofibrations.
  Since toposes are adhesive~\cite{ls:topadh}, the front and right-hand faces of the resulting cube are also strict pullbacks.

  Choose an \F-structure on $f_1$, inducing one on $f_0$ making the left-hand face an \F-morphism.
  Since \F is \stratified, we can give $f_2$ an \F-structure making the back face also an \F-morphism.
  Thus, since \F is \local, by \cref{thm:local} there is an induced \F-structure on $f$, so that in particular $f$ is a fibration.
  So the front and right-hand faces of the cube, being strict pullbacks of $f$, are homotopy pullbacks.
\end{proof}

Recall from~\cite{rezk:homotopy-toposes} that a \textbf{model topos} is a model category that is Quillen equivalent to a left exact left Bousfield localization of the projective model structure on a category of simplicial presheaves.
This implies that the \io-category it presents is a Grothendieck \io-topos.

\begin{cor}\label{thm:ttmt-mt}
  Every \ttmt is a model topos.
\end{cor}
\begin{proof}
  Since it is combinatorial, by~\cite{dug:pres} it has a small presentation, and by \cref{thm:descent} it has descent.
  Thus we can apply~\cite[Theorem 6.9]{rezk:homotopy-toposes} or~\cite[Theorem 6.1.0.6]{lurie:higher-topoi}.
\end{proof}

\begin{rmk}
  On one hand, we have seen that a \ttmt \E has univalent universes classifying the relatively \ka-presentable fibrations for arbitrarily large regular cardinals \ka.
  On the other hand, since the \io-category $\Ho_\oo(\E)$ presented by \E is a Grothendieck \io-topos, by~\cite[Theorem 6.1.6.8]{lurie:higher-topoi} it contains object classifiers for the relatively \ka-presentable morphisms for arbitrarily large regular cardinals \ka.
  One expects the univalent universes in \E to present the object classifiers in $\Ho(\E)$, but this is a subtle question because the relationship between \ka-presentability in \E and in $\Ho_\oo(\E)$ is nontrivial; see~\cite{MO:cptobjio,stenzel:thesis}.
\end{rmk}

\cref{thm:ttmt-mt} shows that \ttmts are, as the term suggests, a subclass of model toposes that are particularly well-adapted to model type theory.
Our main goal is to show that up to homotopy, this subclass involves no loss of generality: every model topos is Quillen equivalent to a type-theoretic one.
We begin with some easy cases.

\begin{prop}\label{thm:ss-ttmt}
  The category \S of simplicial sets, with its Kan-Quillen model structure, is a \ttmt.
\end{prop}
\begin{proof}
  It is a presheaf topos and its tmodel structure is right proper, simplicial, and combinatorial, with its cofibrations the monomorphisms.
  It is \slcc because its copowers are just cartesian products, which are preserved by pullback.
  Finally, $\F=\dFib$ is \local by \cref{eg:rep-cod}, and trivially \stratified.
\end{proof}

Combining \cref{thm:ss-ttmt,thm:uf-fibrant} reproduces Voevodsky's construction~\cite{klv:ssetmodel} of a univalent universe in simplicial sets.

\begin{rmk}
  Specializing the argument for pushouts in \cref{thm:descent} to $\E=\S$, we obtain a proof of the ``cube theorem''~\cite{puppe:cube} for \S.
  A similar proof appears in~\cite[Lemma 6.1.3.12]{lurie:higher-topoi} using minimal fibrations; \cref{thm:u4p} avoids these 
  by providing a different way to turn homotopy pullbacks into strict ones.
\end{rmk}

\begin{prop}\label{thm:slice-ttmt}
  If \E is a \ttmt and $X\in \E$, then $\E/X$ is also a \ttmt.
\end{prop}
\begin{proof}
  A slice of a topos is again a topos, and the slice model structure inherits properness, simplicial-ness, combinatoriality, and cofibrations being the monomorphisms.
  A slice of a slice is a slice, so it inherits \slcclosure.
  And if \F is the \local and \stratified \nfs for the fibrations of \E, then since the forgetful functor $U: \E/X\to \E$ creates pullbacks, colimits, and fibrations, the preimage \nfs $U^{-1}(\F)$ is a \local and \stratified \nfs for the fibrations of $\E/X$.
\end{proof}

\begin{prop}\label{thm:prod-ttmt}
  If $\{\E_i\}_{i\in I}$ is a small family of \ttmts, then the product category $\prod_i \E_i$ is also a \ttmt.
\end{prop}
\begin{proof}
  A product of toposes is a topos; all the structure is inherited pointwise.
\end{proof}

In the rest of the paper we study two more basic constructions of \ttmts that together suffice to obtain our desired generality.

\begin{rmk}\label{rmk:design-space}
  I would not claim that \ttmts are the last word in ``model categories that interpret type theory'', but they do occupy a fairly stable point in the design space: they have a nearly maximal set of good properties one can assume of a model category, are sufficient to model all \io-toposes, interpret (most of) type theory, and are closed under many constructions.

  Note that not every \ttmt is a cartesian \emph{monoidal} model category (e.g.\ this fails already for slices of simplicial sets), but every \io-topos is presented by some \ttmt that does have this property (e.g.\ a left exact localization of an injective model structure on enriched simplicial presheaves).
  Allowing \E to be an arbitrary 1-topos is also somewhat unnecessary generality; all the examples we will construct in this paper are in fact \emph{presheaf} 1-toposes.

  Finally, while I have chosen to stick with simplicial enrichments for simplicity and to facilitate the connection with classical homotopy theory, it should be possible to formulate a more general notion of \ttmt that is enriched over some other monoidal model category, such as some variety of cubical sets.
\end{rmk}

\begin{verbose}
\begin{rmk}
  I would not claim that \ttmts are the last word in definitions of ``model categories that interpret type theory'', but they do seem to occupy a fairly stable point in the design space.
  
  On one hand, they already include a nearly maximal set of the good properties that can be assumed of a model category.
  The one axiom I could imagine adding is that \E is a cartesian \emph{monoidal} model category.
  This is satisfied by left exact localizations of injective model structures on enriched simplicial presheaf categories (which suffice to model all \io-toposes).
  But it fails for slice categories and presheaves on internal categories, and I have not found any particular use for it yet.

  On the other hand, most properties of the definition are necessary for one or another of our arguments.
  \begin{itemize}
  \item Of course, the \local and \stratified \nfs is the whole point, which ensures the existence of strict univalent universes.
  \item The requirement that \E is a 1-topos whose cofibrations are the monomorphisms (a ``Cisinski model category''~\cite{cisinski:presheaves,cisinski:local-acyc}) is extremely useful to relate homotopy theory to 1-categorical descent properties; it appears in \cref{thm:u4p,thm:descent,eg:wo-strat,eg:cof-ff-fcos,eg:span}.
  \item To model $\Pi$-types by 1-categorical dependent products (which appears necessary if we want them to satisfy $\eta$-conversion as well as $\beta$-reduction) we seem to need $f_*$ to preserve fibrations whenever $f$ is a fibration, which by adjointness means that pullback along fibrations preserves acyclic cofibrations.
    Since acyclic fibrations are always stable under pullback, by factorization this implies that pullback along any fibration preserves weak equivalences, i.e.\ that \E is right proper; and the converse also holds if the cofibrations are stable under pullback (as if they are the monomorphisms).
    Note that the \io-categories presentable by right proper Cisinski model categories are precisely the locally cartesian closed locally presentable ones~\cite{cisinski:lccc-rpcmc,gk:univlcc}.
  \item Simplicial enrichment of the model structure is of course very convenient (it is used in \cref{thm:u4p}, for instance), and in \cref{sec:injmodel} it appears necessary to obtain good behavior of cobar constructions.
  \item \Slcclosure appears necessary to control cobar constructions for presheaves on internal categories, and is also used in~\cite{ls:hits} to construct higher inductive types.
  \item  Combinatoriality is used in~\cite{ls:hits} for algebraic fibrant replacements, and also implies the existence of the Bousfield localizations used in \cref{sec:lex-loc}.
  \end{itemize}
  One might imagine weakening combinatoriality to accessibility~\cite{rosicky:acc-model}, but in fact any accessible model structure on a 1-topos whose cofibrations are the monomorphisms is automatically combinatorial.
    For the monomorphisms in a 1-topos are always cofibrantly generated by a set (cf.~\cite[Proposition 1.2.27]{cisinski:presheaves},~\cite[Proposition 1.2.2]{cisinski:local-acyc}, and~\cite[Proposition 1.12]{beke:sheafifiable}); so the weak equivalences are accessible and accessibly embedded in $\E^\dtwo$, being the preimage of the small-injectivity class of acyclic fibrations~\cite[Proposition 3.3]{rosicky:comb-model} under an accessible fibration-replacement functor.
    The other hypotheses of Smith's theorem (\cite[Theorem 1.7]{beke:sheafifiable} or~\cite[Proposition 1.7]{barwick:enr-localization}) or~\cite[Theorem 1.4.3]{cisinski:presheaves} are automatically satisfied if we already have a model structure.
\end{rmk}
\end{verbose}


%% file: bar.tex
\section{Coherent transformations and bar constructions}
\label{sec:coh-bar}

In preparation for our treatment of injective model structures using cobar constructions in \cref{sec:injmodel}, in this section we give some intuition for how (co)bar constructions arise and review some of their formal properties.
Let \E be a model category, and suppose for simplicity that \D is an ordinary small category.
The most obvious notions of weak equivalence, fibration, and cofibration in $\pr\D\E$ are induced pointwise from \E, and functorial factorizations in \E can also be applied pointwise.
However, there is a problem with the lifting properties: suppose we have a square in $\pr\D\E$
\[
  \begin{tikzcd}
    A \ar[r,"f"] \ar[d,"i"'] & X \ar[d,"v"] \\
    B \ar[r,"g"'] & Y
  \end{tikzcd}
\]
in which $i$ is a pointwise cofibration and $v$ a pointwise acyclic fibration.
For each $d\in\D$ we have a lift $h_d : B_d \to X_d$ with $h_d \circ i_d = f_d$ and $v_d \circ h_d = g_d$, but it may not be possible to choose these \emph{naturally} in $d$: for $\delta : d_1\to d_2$ we may not have $X_{\delta} \circ h_{d_2} = h_{d_1}\circ B_\delta$.
However, both $X_{\delta} \circ h_{d_2}$ and $h_{d_1}\circ B_\delta$ are lifts in the square
\[
  \begin{tikzcd}[column sep=large]
    A_{d_2} \ar[rr,"{X_\delta \circ f_{d_2} = f_{d_1}\circ A_\delta}"] \ar[d,"{i_{d_2}}"'] & & X_{d_1} \ar[d,"{v_{d_1}}"] \\
    B_{d_2} \ar[rr,"{Y_\delta \circ g_{d_2} = g_{d_1}\circ B_\delta}"'] & & Y_{d_1},
  \end{tikzcd}
\]
and the space of such lifts is contractible, being a fiber of the acyclic fibration
\[\ehom\E(B_{d_2},X_{d_1}) \to \ehom\E(A_{d_2},X_{d_1}) \times_{\ehom\E(A_{d_2},Y_{d_1})} \ehom\E(B_{d_2},Y_{d_1}).\]
Thus we have a \emph{homotopy} $X_{\delta} \circ h_{d_2} \sim h_{d_1}\circ B_\delta$ over $v_{d_1}$ and under $i_{d_3}$.
Similarly, given $\delta' : d_2\to d_3$ there is a 2-simplex in $\ehom\E(B_{d_3},X_{d_1})$ relating these homotopies for $\delta$, $\delta'$, and $\delta'\circ\delta$, and so on, yielding a \emph{homotopy coherent natural transformation} $h:B\cohto X$ such that $v \circ h = g$ and $h\circ i = f$.

There is an entire theory of homotopy coherent transformations (see e.g.~\cite{cp:hcct}), but as our purpose at the moment is motivational we only sketch it.
For $X,Y\in\pr\D\E$, a homotopy coherent transformation $h:X\cohto Y$ consists of:
\begin{itemize}
\item For every $d\in \D$, a morphism $h_d:X_d \to Y_d$.
\item For every $d_1\xto{\delta} d_2$ in \D, a homotopy $h_\delta: \Delta[1]\cpw X_{d_2} \to Y_{d_1}$ between $Y_{\delta} \circ h_{d_2}$ and $h_{d_1}\circ X_\delta$, such that $h_{\id_d}$ is constant.
\item For every $d_1\xto{\delta} d_2 \xto{\delta'} d_3$ in \D, a 2-simplex $h_{\delta,\delta'}:\Delta[2]\cpw X_{d_3} \to Y_{d_1}$ whose boundaries involve $h_\delta$, $h_{\delta'}$, and $h_{\delta'\circ\delta}$, satisfying similar constancy conditions.
\item And so on.
\end{itemize}
Let $U:\pr\D\E \to \E^{\ob\D}$ denote the forgetful functor, with $F$ its left adjoint defined by $(F W)_d = \coprod_{d'} \D(d,d') \cpw W_{d'}$; note both preserves simplicial copowers.
Thus:
\begin{itemize}
\item A collection of morphisms $h_d:X_d \to Y_d$ constitutes a morphism $U X \to U Y$, or equivalently $F U X \to Y$.
\item A collection of morphisms $h_\delta: \Delta[1]\cpw X_{d_2} \to Y_{d_1}$ constitutes a morphism $\Delta[1] \cpw U F U X \to U Y$, or equivalently $\Delta[1] \cpw F U F U X \to Y$.
\item A collection of morphisms $h_{\delta,\delta'}:\Delta[2]\cpw X_{d_3} \to Y_{d_1}$ constitutes a morphism $\Delta[2]\cpw U F U F U X \to U Y$, or equivalently $\Delta[2]\cpw F U F U F U X \to Y$.
\item And so on.
\end{itemize}
The fact that all these transformations have the right boundaries and constancy conditions means precisely that they assemble into a single \emph{strict} natural tranformation $\bar(F,UF,UX) \to Y$, where $\bar(F,UF,UX)$ is the geometric realization of the ``two-sided simplicial monadic bar construction''~\cite{may:goils,meyer:bar_i} $\sbar(F,UF,UX)$:
\[\small
  \begin{tikzcd}
    \cdots \ar[r,-] \ar[r,-,shift left=1] \ar[r,-,shift left=2] \ar[r,-,shift left=3]
     \ar[r,-,shift right=1] \ar[r,-,shift right=2] \ar[r,-,shift right=3]
    &
    FUFUFUX \ar[r] \ar[r,shift left=4] \ar[r,shift right=4] &
    FUFUX \ar[r,shift left=2] \ar[r, shift right=2] \ar[l,shift left=2] \ar[l,shift right=2] &
    FUX \ar[l]
  \end{tikzcd}
\]
whose face and degeneracy maps are defined by the unit and counit of the adjunction $F\adj U$.
In other words, $\bar(F,UF,UX)$ is a ``classifier'' for homotopy coherent transformations with domain $X$: we have a natural bijection between homotopy coherent transformations $X\cohto Y$ and strict transformations $\bar(F,UF,UX)\to Y$.

In particular, the identity map of $\bar(F,UF,UX)$ corresponds to a universal homotopy coherent transformation $p_X:X\cohto \bar(F,UF,UX)$, which is universal in that any homotopy coherent transformation $h:X\cohto Y$ can be written as\footnote{To compose two homotopy coherent transformations we need fibrancy and cofibrancy conditions to obtain horn-fillers in hom-spaces, but there is no trouble composing a homotopy coherent transformation with a strict one on either side.} $h = \overline{h}\circ p_X$ for a unique strict transformation $\overline{h}:\bar(F,UF,UX)\to Y$.
The identity map of $X$ corresponds to a canonical strict transformation $q_X = \overline{\id_X}:\bar(F,UF,UX)\to X$ such that $q_X \circ p_X = \id_X$, and it is a fact (see \cref{thm:bar-repl}) that we also have a homotopy $p_X \circ q_X \sim \id_{\bar(F,UF,UX)}$, so that $p_X$ and $q_X$ are inverse simplicial homotopy equivalences.
Moreover, this correspondence is natural with respect to strict transformations $v:Y\to Z$, i.e.\ we have $\overline{v\circ h} = v\circ \overline{h}$, and in particular for a strict $w:X\to Z$ we have $\overline{w} = w\circ q_X$.

\begin{rmk}\label{rmk:2mnd}
  This description of the bar construction as a classifier for homotopy coherent transformations is a generalization to homotopy theory of the \emph{pseudomorphism classifiers} of 2-monad theory.
  As shown in~\cite{bkp:2dmonads}, for any suitable 2-monad $T$, the inclusion $T\algs \into T\alg$ of the 2-category of $T$-algebras and strict morphisms into the category of $T$-algebras and pseudomorphisms has a left adjoint, traditionally denoted $A \mapsto \pscl A$.
  Thus pseudo $T$-morphisms $A\cohto B$ are in bijection with strict $T$-morphisms $\pscl A\to B$.
  In~\cite{lack:codescent-coh} the pseudomorphism classifier $\pscl A$ is constructed as a codescent object, which is really just a 2-truncated bar construction.
  (At present we are thinking only about the monad $U F$ on $\E^{\ob\D}$ whose category of algebras is $\pr\D\E$, but like the pseudomorphism classifier, the bar construction makes sense for any monad.)
\end{rmk}

The above discussion suggests that to obtain lifting properties in $\pr\D\E$, we need to be able to ``rectify'' homotopy coherent natural transformations to strict ones.
But by universality, if $p_X$ is homotopic to a strict transformation $s:X\to \bar(F,UF,UX)$ then the same is true of every homotopy coherent transformation with domain $X$.
And since $q_X$ is strict and a simplicial homotopy inverse of $p_X$, this is equivalent to saying that $q_X$ has a simplicial homotopy inverse in $\pr\D\E$.

If $q_X$ has a homotopy inverse in $\pr\D\E$ and $X$ is also pointwise cofibrant, we may call it \emph{projectively semi-cofibrant}.
In this case, if $v:Y\to Z$ is a pointwise acyclic fibration and we have $g:X\to Z$, then by choosing lifts and homotopies as above we can produce a homotopy coherent lift $h:X\cohto Y$ such that $v \circ h = g$, and then rectify it to an equivalent strict transformation $k:X\to Y$ with $k\sim h$.
But then we have only a homotopy $v\circ k \sim v\circ h = g$, i.e.\ although $k$ is a strict transformation, it is only a lift up to homotopy.

Thus we actually need a stronger property: that $s:X\to \bar(F,UF,UX)$ is a \emph{strict} section of $q_X$, i.e.\ $q_X \circ s = \id_X$.
(It is then automatically a simplicial homotopy inverse.)
In this case, the strict transformation $k$ is defined by $\overline{h} \circ s$, and we have
\[v \circ k = v\circ \overline{h}\circ s = \overline{v\circ h} \circ s = \overline{g} \circ s = g\circ q_X \circ s = g \]
so that $k$ is also a strict lift.
This discussion has been somewhat informal, but we will show more carefully in \cref{thm:injfib} below that indeed, under suitable hypotheses, $\emptyset\to X$ has left lifting for pointwise fibrations if and only if it is pointwise cofibrant and $q_X$ has a strict section.
Thus these are the \emph{projectively cofibrant} objects, which indeed are the cofibrant objects in a model structure whose weak equivalences and fibrations are pointwise.

\begin{rmk}\label{rmk:flexible}
  Continuing \cref{rmk:2mnd}, in~\cite{bkp:2dmonads} an algebra $A$ for a 2-monad $T$ was defined to be \emph{semi-flexible} if the map $q : \pscl A\to A$ is an equivalence in the 2-category $T\algs$ of $T$-algebras and strict morphisms, and \emph{flexible} if this $q$ has a section (making it automatically also an equivalence).
  In~\cite{lack:htpy-2monads} it was shown that for a suitable 2-monad $T$, the category $T\algs$ admits a model structure whose ``homotopy 2-category'' is $T\alg$; whose weak equivalences and fibrations are the strict $T$-algebra morphisms that are equivalences and fibrations, respectively, in the underlying 2-category; and in which the cofibrant objects are precisely the flexible ones.
  Thus, when $T$ is the monad whose algebras are presheaves, this model structure is the projective one.
  Our characterization of projectively cofibrant objects is directly inspired by, and generalizes, this result of~\cite{lack:htpy-2monads}.
\end{rmk}

In fact this entire discussion works for any suitable adjunction $F\adj U$, with the case of most interest being when the adjunction is monadic (as in the case of presheaves).
However, in the case of presheaves, the forgetful functor $U$ is also \emph{comonadic}, which means that the entire discussion can be dualized, producing a characterization of the \emph{injectively fibrant} objects.
And with a little extra work we can generalize this to a characterization of all the injective fibrations using a \local and \stratified \nfs.

In our formal treatment, however, it is easier not to explicitly discuss homotopy coherent natural transformations at all, but rather focus on the bar construction.
In the rest of this section we give an abstract formulation of the bar construction, following~\cite{rv:hcadj-ftm,rv:elements}, including proofs of its (well-known) basic properties.

By~\cite{ss:free-adj,rv:hcadj-ftm}, the \textbf{free adjunction} $\cAdj$ is a 2-category freely generated by two objects $\zero$ and $\one$ (called $+$ and $-$ in~\cite{rv:hcadj-ftm}) and a morphism $f:\zero\to \one$ with right adjoint $u:\one\to \zero$.
Its hom-categories are:
\begin{itemize}
\item $\cAdj(\zero,\zero) = \dDelta_+$, the augmented simplex category, which can be identified with the category of (possibly-empty) finite ordinals and monotone maps.
\item $\cAdj(\one,\one) = \dDelta_+\op$, which can be identified with the category of non-empty finite ordinals and monotone maps that preserve the top and bottom elements.
  This description of its full subcategory $\dDelta\op$, corresponding to the ordinals with \emph{distinct} top and bottom elements, is well-known.
\item $\cAdj(\zero,\one) = \dDelta_\top$, the category of non-empty finite ordinals and monotone maps that preserve the top element.\footnote{This is the convention of~\cite{ss:free-adj}, which is the opposite of that of~\cite{rv:hcadj-ftm}.
  The choice is essentially arbitrary, since $\dDelta_\top\cong\dDelta_\bot$ by reversing the order of finite ordinals.}
\item Similarly, $\cAdj(\one,\zero) = \dDelta_\bot$, the category of non-empty finite ordinals and monotone maps that preserve the bottom element.
\end{itemize}
Similarly to the identification of $\dDelta_+\op$ as the subcategory of $\dDelta_+$ consisting of maps preserving the top and bottom elements, we have an isomorphism $\dDelta_\top \cong \dDelta_\bot\op$; see~\cite[Observation 3.3.6]{rv:hcadj-ftm} for a detailed discussion.

\begin{defn}
  Given an adjunction $F:\N \toot \M:U$, determining a 2-functor $V : \cAdj\to\cCat$ with $V(\zero)=\N$ and $V(\one)=\M$, the \textbf{two-sided simplicial bar construction} is the composite
  \[ \dDelta\op \into \dDelta_+\op = \cAdj(\one,\one) \xto{V} \cCat(\M,\M) \]
  rearranged into a functor
  \[ \sbar(F,UF,U-) : \M \to \pr\dDelta\M \]
  sending each object $X\in \M$ to a simplicial object $\sbar(F,UF,UX)\in \pr\dDelta\M$.
\end{defn}

By inspection,
\( \bar_n(F,UF,UX) = \overbrace{(FU)\dotsm(FU)}^{n+1} X \)
with faces and degeneracies obtained respectively from the counit $FU\to\Id_\M$ and unit $\Id_\N\to UF$.
Thus our bar construction specializes to the classical ones for monads~\cite{may:goils,meyer:bar_i}, enriched categories~\cite{cp:hcct,shulman:htpylim}, and internal categories~\cite{may:csf,horel:model-intsscat}.
The extension to $\dDelta_+\op$ gives an augmentation that is just the counit $FUX \to X$.

Now recall that for any category \N, the category $\pr\dDelta\N$ of simplicial objects in \N has a \textbf{canonical enrichment} over the category $\S=\pr\dDelta\nSet$ of simplicial sets, 
defined most easily in terms of its powers and copowers, which exist whenever \N is complete and cocomplete:
\[
  (K\cpw X_\bullet)_n = K_n \cdot X_n \qquad
  (K\pow Y_\bullet)_n = Y_n^{K_n}.
\]

\begin{lem}\label{thm:bar-simpcontr}
  For any adjunction $F\adj U$, the composite
  \[ \M \xto{\sbar(F,UF,U-)} \pr\dDelta\M \xto{U} \pr\dDelta\N \]
  is naturally simplicially contractible to $U$.
  In other words, the augmentation
  \[ \ep : U\sbar(F,UF,U-) = \sbar(UF,UF,U-) \too U \]
  has a natural section $\sigma:U\to \sbar(UF,UF,U-)$ (so that $\ep\sigma=\id$), with a homotopy $H:\sigma\ep \sim \id$ in the canonical enrichment of $\pr\dDelta\N$ such that $H\sigma$ is constant.
\end{lem}
\begin{proof}
  This is a classical fact (e.g.~\cite[Proposition 9.8]{may:goils} or~\cite[\sect6--7]{meyer:bar_i}), traditionally proven by exhibiting an ``extra degeneracy'', defining generators for a homotopy, and checking some simplicial identities (or leaving them to the reader).
  We repackage this a bit more abstractly, following~\cite{rv:hcadj-ftm} and~\cite[Chapters 9 and 10]{rv:elements} but working with simplicially enriched categories rather than quasicategories.

  We start by applying the 2-functor $V : \cAdj\to\cCat$ to a different hom-category:
  \[ \cAdj(\one,\zero) \xto{V} \cCat(\M,\N). \]
  This extends the augmented simplicial object $\sbar(UF,UF,U-) \to U$ to a diagram on the category $\cAdj(\one,\zero) = \dDelta_\bot$ of non-empty finite ordinals and bottom-preserving monotone maps.
  (The maps that also preserve the top element are those of the original augmented simplicial object; the rest are the ``extra degeneracy''.
  The whole $\dDelta_\bot$-diagram is the image of the comparison functor of~\cite[\sect 7]{rv:hcadj-ftm} sending each $X\in\M$ to the ``homotopy coherent $U F$-algebra''~\cite[eq.~(6.1.12)]{rv:hcadj-ftm} associated to the strict $U F$-algebra $U X$.)
  Thus it suffices to prove the following lemma.
\end{proof}

\begin{lem}
  The underlying simplicial object of \emph{any} $X\in \func{\dDelta_\bot}{\N}$ is naturally simplicially contractible to its augmentation $X_{-1}$.
\end{lem}
\begin{proof}
  Again, this is a classical fact in homotopy theory (more recently stated \io-categorically as~\cite[Lemma 6.1.3.16]{lurie:higher-topoi} and~\cite[Theorem 5.3.1]{rv:2cat-quasi}); we simply repackage it abstractly.
  Let $\cR$ be the simplicially enriched category generated by a strong deformation retraction.
  It has two objects $\iI$ and $\iR$ with hom-objects:
  \begin{itemize}
  \item $\cR(\iI,\iR) = \Delta[0]$, the discrete simplicial set with one vertex $s$.
  \item $\cR(\iR,\iI) = \Delta[0]$, discrete on one vertex $e$.
  \item $\cR(\iI,\iI) = \Delta[0]$, discrete on one vertex $\id_\iI$ (so in particular $e s = \id_{\iI}$).
  \item $\cR(\iR,\iR) = \Delta[1]$, a 1-simplex from $\id_\iR$ to the composite $s e$.
  \end{itemize}
  Thus the objects of the simplicially enriched presheaf category $\pr\cR\S$ are pairs of simplicial sets related by a strong deformation retraction.
  In particular, the representable $\cR(-,\iI)$ is the identity retraction from $\Delta[0]$ to itself, while the representable $\cR(-,\iR)$ is the retraction of $\Delta[1]$ onto its right-hand vertex.

  Now, as for any simplicially enriched category \cR (see e.g.~\cite[C2.5.3]{ptj:elephant}), $\pr\cR\S$ can be identified with an \emph{unenriched} presheaf category $\pr{(\drr)}\nSet$.
  In our case, $\drr$ is the full subcategory of $\pr\cR\S$ on the objects $\cR(-,\iI)\times \Delta[n]$ (the identity retraction from the $n$-simplex $\Delta[n]$ to itself) and $\cR(-,\iR)\times \Delta[n]$ (the retraction of the prism $\Delta[n]\times\Delta[1]$ onto its right-hand $n$-simplex face).

  Here we use the combinatorial input: $\Delta[n]\times\Delta[1]$ contains $\Delta[n+1]$ as a retract, where the left-hand face of the prism is also a face of $\Delta[n+1]$, and the final vertex of $\Delta[n+1]$ lies in the right-hand face.
  Thus, we can define a functor from $\drr$ to the category of based simplices $(\Delta[n],n)$ whose basepoint is their last vertex, by
  \begin{align*}
    \cR(-,\iR)\times \Delta[n] &\;\mapsto\; (\Delta[n+1],n+1)\\
    \cR(-,\iI)\times \Delta[n] &\;\mapsto\; (\Delta[0],0).
  \end{align*}
  However, this category of based simplices is isomorphic to $\dDelta_\top$, which as we noted above is isomorphic to $\dDelta_\bot\op$.
  Thus by composition with it we obtain a functor
  \[ \func{\dDelta_\bot}{\nSet} \simeq \func{\dDelta_\top\op}{\nSet} \to \pr{(\drr)}{\nSet} \simeq \pr\cR\S \]
  extracting from any $\dDelta_\bot$-diagram of sets a strong deformation retraction of simplicial sets.
  Finally, by composing with the Yoneda embedding of \N we obtain
  \begin{equation*}\small
  \func{\dDelta_\bot}{\N} \to \func{\dDelta_\bot}{\pr\N\nSet}
    \simeq \pr\N{\func{\dDelta_\bot}{\nSet}}
    \to \pr\N{\pr\cR\S} \simeq \pr\cR{\pr\N\S}
  \end{equation*}
  sending any $X\in\func{\dDelta_\bot}{\N}$ to a strong deformation retraction of \S-valued presheaves on \N.
  The two \S-valued presheaves involved in this retraction are respectively representable by the underlying simplicial object of $X$ and the constant simplicial object at $X_{-1}$; thus the strong deformation retraction lies entirely in $\pr\dDelta\N$.
\end{proof}

If \M is simplicially enriched and cocomplete with copowers, the \textbf{geometric realization} of $X_\bullet \in \pr\dDelta\M$ is the coend
\( |X_\bullet| = \int^{n\in\dDelta} \Delta[n] \cpw X_n \),
and the \textbf{(realized) two-sided bar construction} is the geometric realization of the simplicial one:
\[ \bar(F,UF,U-) = |\sbar(F,UF,U-)|. \]
By~\cite[Proposition 5.4]{rss:simp}, geometric realization is a simplicially enriched functor $\pr\dDelta\M \to \M$ when \M has its given enrichment while $\pr\dDelta\M$ has its above {canonical} enrichment. 
Thus we obtain:

\begin{cor}\label{thm:bar-repl}
  For any simplicially enriched adjunction $F\adj U$ such that $U$ preserves geometric realizations, the two-sided bar construction defines a functor
  \[ \bar(F,UF,U-) : \M \to \M \]
  with a natural augmentation $\bar(F,UF,UX)\to X$ whose image under $U$ is a simplicial strong deformation retraction in \N, and hence a weak equivalence if \N is a simplicial model category.
\end{cor}
\begin{proof}
  Since $U$ preserves geometric realization, $U\bar(F,UF,U-)$ is the geometric realization of $U\sbar(F,UF,U-)$.
  By \cref{thm:bar-simpcontr} the latter is simplicially contractible to $U$ in $\pr\dDelta\N$.
  Thus, since geometric realization is a simplicial functor, $U\bar(F,UF,U-)$ is simplicially contractible to $U$ in $\N$.
  The last statement follows since simplicial homotopy equivalences in a simplicial model category are weak equivalences.
\end{proof}

\begin{rmk}
  If the adjunction is monadic, then $U$ preserves geometric realizations just when the corresponding monad does.
  For us this will usually be because it is a simplicial left adjoint, but there are non-left-adjoint monads that preserve geometric realizations; notably, those associated to operads, as in~\cite[Theorem 12.2]{may:goils} where the monadic two-sided bar construction was first introduced.
\end{rmk}


%% file: injmodel.tex
\section{Injective model structures}
\label{sec:injmodel}

In this section we will show that \ttmts are closed under passage to presheaf categories with injective model structures.
Although for our main theorem it would suffice to use \emph{enriched} presheaf categories $\pr\D\E$, 
where \D is a small simplicially enriched category, we will also include \emph{internal} presheaf categories $\pr\cD\E$, where \cD is an internal category in \E.
Combined with \cref{sec:lex-loc}, this will imply that \ttmts are closed under passage to ``internal sheaves on internal sites''.

\begin{defn}
  For any functor $U:\M\to\N$ where \N is a model category, we define a morphism in \M to be:
  \begin{itemize}
  \item A \textbf{$U$-weak equivalence}, \textbf{$U$-cofibration}, or \textbf{$U$-fibration} if its image under $U$ belongs to the respective class in \N.
  \item A \textbf{projective cofibration} if it has left lifting for all $U$-acyclic $U$-fibrations.
  \item A \textbf{projective acyclic cofibration} if it has left lifting for all $U$-fibrations.
  \item An \textbf{injective fibration} if it has right lifting for all $U$-acyclic $U$-cofibrations.
  \item An \textbf{injective acyclic fibration} if it has right lifting for all $U$-cofibrations.
  \end{itemize}
  The \textbf{projective} or \textbf{right-lifted model structure}, if it exists, consists of the $U$-weak equivalences, $U$-fibrations, and projective cofibrations.
  Similarly, the \textbf{injective} or \textbf{left-lifted model structure}, if it exists, consists of the $U$-weak equivalences, $U$-cofibrations, and injective fibrations.
\end{defn}

The following is a compilation of some well-known and more recent results.

\begin{prop}\label{thm:projmodel-gen}
  Let \N be a model category and $U:\M\to\N$ a functor with both adjoints $F\adj U \adj G$, such that the adjunction $U F \adj U G$ is Quillen.
  Then:
  \begin{enumerate}
  \item If \N is cofibrantly generated, the projective model structure exists and is cofibrantly generated.
  \item If \N is an accessible model category and \M is locally presentable, then both projective and injective model structures exist and are accessible.
  \item If \N is a combinatorial model category and \M is locally presentable, then both projective and injective model structures exist and are combinatorial.
  \end{enumerate}
  Moreover:
  \begin{itemize}
  \item Every projective cofibration is a $U$-cofibration, and every injective fibration is a $U$-fibration.
  \item The projective and injective model structures are right or left proper if \N is.
  \item If \N is a \V-model category for some monoidal model category \V, and \M is a \V-category and the adjunctions are \V-adjunctions, then the projective and injective model structures are also \V-model categories.
  \end{itemize}
\end{prop}
\begin{proof}
  We begin by constructing factorizations.
  The cofibrantly generated case is standard: we take the generating cofibrations and acyclic cofibrations of \M to be the $F$-images of those of \N, whose domains are small in \M because $U$ preserves all colimits.
  The accessible case follows from~\cite{gkr:liftacc-modelstr}, and the combinatorial case from~\cite[Remark 3.8]{mr:cellular}.

  We sketch how to complete the proof in the projective case (the injective is dual).
  By retract arguments, it suffices to show every projective acyclic cofibration $f$ is a $U$-weak equivalence.
  We show $U f$ is an acyclic cofibration, i.e.\ that it has left lifting for any fibration $h$.
  This is equivalent to $f$ having left lifting for $G h$, so it will suffice to show $G h$ is a $U$-fibration, i.e.\ that $U G h$ is a fibration; but $U G$ preserves fibrations since it is right Quillen.

  A similar argument shows that projective cofibrations are $U$-cofibrations.
  Both properness claims now follow since $U$ preserves pullbacks, pushouts, cofibrations, and fibrations, and creates weak equivalences.
  Finally, enrichment of the model structure is easy using the characterization in terms of powers, since $U$ preserves powers and pullbacks and creates fibrations and acyclic fibrations.
\end{proof}

\begin{rmk}
  \cite[Remark 3.8]{mr:cellular} uses intricate calculations as in~\cite[A.3.3.3]{lurie:higher-topoi}, but in our case there is a simpler argument.
  If the cofibrations of \N are the monomorphisms, $U$ is faithful, and \M is a topos, then the $U$-cofibrations are also the monomorphisms, hence cofibrantly generated by~\cite[Proposition 1.12]{beke:sheafifiable}.
  And the $U$-weak equivalences are accessible and accessibly embedded, being the $U$-preimage of the weak equivalences of \N, so Smith's theorem applies.
\end{rmk}

Often $U$ will be both monadic and comonadic, so that \M is the category of $UF$-algebras and also of $UG$-coalgebras.

\begin{eg}\label{eg:enrcat}
  If \V is a monoidal category, \E a complete and cocomplete \V-category, and \D a small \V-category, then the forgetful functor $U : \pr\D\E \to \E^{\ob\D}$ is both monadic and comonadic. 
  In the accessible case this specializes to~\cite{moser:injproj}, and we have an analogous result in the combinatorial case: the projective and injective model structures exist if the copowers $(\D(x,y)\cpw -)$ preserves cofibrations and acyclic cofibrations.
  For instance, this occurs if \V is a monoidal model category, \E is a \V-model category, and each $\D(x,y)$ is cofibrant.
  The earliest reference I know for both model structures (when $\E=\S$) is~\cite{heller:htpythys}; later references for projective model structures include~\cite{piacenza:hodia,hirschhorn:modelcats,dro:enrfunc} and for injective ones~\cite{lurie:higher-topoi,moser:injproj}.
\end{eg}

\begin{eg}\label{eg:intcat}
  If \E is a locally cartesian closed category and $\cD = (D_1 \toto D_0)$ an internal category in \E, then the forgetful functor $\pr \cD \E \to \E/D_0$ is monadic and comonadic, where $\pr \cD \E$ is the ``internal presheaf category''.
  Thus, if \E is also an accessible model category, then $\pr \cD \E$ has projective and injective model structures if $(D_1 \times_{D_0} -)$ preserves cofibrations and acyclic cofibrations.
  In particular, this occurs if the cofibrations in \E are the monomorphisms and the target map $D_1 \to D_0$ is a \textbf{sharp morphism}~\cite{rezk:sharp-maps}, i.e.\ the pullback functor $\E/D_0 \to \E/D_1$ preserves weak equivalences (e.g.\ \E is right proper and it is a fibration).
  For $\E=\S$, this projective model structure was constructed in~\cite[Proposition 6.6]{horel:model-intsscat}.
\end{eg}

Now a bar construction has always been thought of as a ``cofibrant resolution'' of a sort; we will show that it is in fact a \emph{projective} cofibrant replacement.
Related results include~\cite[Proposition 14.8.8]{hirschhorn:modelcats} and~\cite{gambino:wgtlim} as well as~\cite[Propositions 6.9 and 6.10]{horel:model-intsscat}.
However, we need to strengthen the hypotheses of \cref{thm:projmodel-gen} a little.

Following~\cite[\sect 4]{rv:reedy}, for a two-variable functor $\Asterisk:\M\times\N\to\P$ we write its \textbf{Leibniz} (or ``pushout-product'') two-variable functor as $\widehat{\Asterisk} :\M^\dtwo\times\N^\dtwo\to\P^\dtwo$.
In order to use this notation for the ``application'' functor $\func\M\N \times \M\to \N$, we denote the latter by $\app$; thus $F\app X = F(X)$.

\begin{defn}
  Let $S,T:\M\to \N$ be functors between a pair of model categories, and $\alpha:S\to T$ a natural transformation.
  We say $\alpha$ is a \textbf{\qcof} if for any cofibration $i:A\to B$ in \M, the ``Leibniz application'' $\alpha \lapp i$:
  \begin{equation}
    \begin{tikzcd}
      S A \ar[r,"\alpha_A"] \ar[d,"S i"'] \ar[dr,phantom,near end,"\ulcorner"] & T A \ar[d] \ar[ddr,"T i"]\\
      S B \ar[r] \ar[drr,"\alpha_B"'] & \bullet 
      \ar[dr,dashed,"{\alpha\lapp i}" description] \\
      & & T B
    \end{tikzcd}\label{eq:lapp}
  \end{equation}
  is a cofibration that is acyclic if $i$ is.
\end{defn}

\begin{lem}\label{thm:leibcomp}
  Given \qcofs $\al:S\to T$ and $\be:P\to Q$ between pushout-preserving functors $S,T:\M\to\N$ and $P,Q:\N\to\P$, their ``Leibniz composite''
  \[ \be\lcirc\al :  (Q\circ S) \amalg_{P\circ S} (P\circ T) \too (Q\circ T) \]
  is again a \qcof.
\end{lem}
\begin{proof}
  By the usual arguments (cf.~\cite[Observation 4.7]{rv:reedy}), we have
  $(\be\lcirc\al)\lapp i \cong \be \lapp (\al \lapp i)$,
  so we can apply the assumptions on $\al$ and $\be$ in succession.
\end{proof}

\begin{rmk}
  There is a dual notion of a \textbf{\qfib}.
  If $S$ and $T$ have right adjoints $S^*$ and $T^*$, then $\alpha:S\to T$ is a \qcof if and only if its mate $\alpha^* : T^* \to S^*$ is a \qfib.
  \qcofs can also be defined by a left lifting property in $\func\M\N$ against ``Leibniz right Kan extensions'' of a cofibration and a fibration one of which is acyclic, using the fact that the two-variable application functor
  \( \func\M\N \times \M \to \N \)
  has an adjoint on one side.
\end{rmk}

\begin{eg}
  A functor is \textbf{\qcoft} (i.e.\ $\emptyset \to T$ is a \qcof, where $\emptyset$ is the functor constant at the initial object) if and only if it preserves cofibrations and acyclic cofibrations.
  Hence \qcoft left adjoints are precisely left Quillen functors (hence the name ``\qcof'').
\end{eg}

\begin{defn}
  If $T$ is a monad on a model category whose unit $\Id\to T$ is a \qcof, we say that $T$ is \textbf{\qucoft}.
  Dually, we have the notion of a \textbf{\qufibt} comonad.
\end{defn}

Since the identity functor is \qcoft, any \qucoft monad is also a \qcoft functor, but not conversely.

\begin{eg}\label{thm:monmodel-lcof}
  If \V is a monoidal model category, \N is a \V-model category, and $j:K\to L$ is a cofibration in \V, then $(j\cpw -) : (K\cpw -) \to (L\cpw -)$ is a \qcof.
  In particular, the monad on \N induced by a monoid $D$ in \V, whose algebras are $D$-modules, is \qucoft if the unit map $\unit\to D$ is a cofibration in \V.
\end{eg}

\begin{eg}\label{eg:mat}
  Let \V be a monoidal model category and \E a \V-model category, and let $\sO$ be a small set.
  Let $\W=\V^{\sO\times\sO}$ with the ``matrix multiplication'' monoidal structure, $(A\otimes B)(x,y) = \coprod_z A(x,z) \otimes B(z,y)$ with unit $\unit_\W(x,y) = \unit_\V$ if $x=y$ and $\emptyset$ otherwise.
  Let $\N = \E^\sO$ with the usual product model structure, with $\W$-enriched hom-objects $\ehom{\N}(X,Y)(x,y) = \ehom{\E}(X(x),Y(y))$, so that the copower is ``matrix multiplication'' $(A\cpw X)(x) = \coprod_y A(x,y) \otimes X(y)$.
  Then \W is a 
  monoidal model category and \N is a \W-model category, so \cref{thm:monmodel-lcof} applies.

  In particular, note that a small \V-category \D can be regarded as a monoid in \W, and the algebras for the resulting monad on \N are enriched presheaves on \D.
  This monad is \qucoft if each hom-object $\D(x,y)$ is cofibrant and each unit map $\unit \to\D(x,x)$ is a cofibration in \V.
\end{eg}

\begin{eg}\label{eg:span}
  Let \E be a simplicial model category that is a Grothendieck topos and whose cofibrations are the monomorphisms, and $O\in \E$ an object.
  Let $\W= \E/(O\times O)$, regarded as the category of spans $O \xot{s} A \xto{t} O$ with the span-composition monoidal structure (pullback over $O$).
  Let $\N = \E/O$ with the usual slice model structure, with $\W$-enrichment $\ehom{\N}(X,Y)$ the local exponential in $\W$ of $\pi_1^* X$ and $\pi_2^* Y$, so that the copower is again pullback over $O$:
  \[
    \begin{tikzcd}[row sep=small,column sep=small]
      && A\cpw X\ar[dl] \ar[dr] \ar[dd,near start,phantom,"{\rotatebox{-45}{$\lrcorner$}}"] \\
      & A \ar[dl,"s"] \ar[dr,"t"'] && X \ar[dl]\\
      O & {} & O
    \end{tikzcd}
  \]
  A monoid in \W is an internal category in \E, and the algebras for the induced monad on \N are the internal presheaves.

  Now \W is \emph{not} in general a monoidal model category.
  But if $i:A_1\to A_2$ is a monomorphism for which the target morphisms $A_1 \to O$ and $A_2\to O$ are sharp morphisms, then the induced morphism of endofunctors of \N is a \qcof.
  To see this, let $j:B_1\to B_2$ be a cofibration in \N and construct the map:
  \[
    \begin{tikzcd}
      A_1\times_O B_1 \ar[d,"{i\times_O B_1}"'] \ar[r,"{A_1 \times_O j}"] \ar[dr,phantom, near end,"\ulcorner"] &
      A_1\times_O B_2 \ar[d] \ar[ddr,"{i\times_O B_2}"] \\
      A_2 \times_O B_1 \ar[r] \ar[drr,"{A_2 \times_O j}"'] & \bullet \ar[dr,"{i \mathbin{\widehat{\times_O}} j}" description]\\
      && A_2\times_O B_2
    \end{tikzcd}
  \]
  But the outer square is a pullback, hence the pushout is a union of subobjects in the topos $\E/O$; so $i \mathbin{\widehat{\times_O}} j$ is a monomorphism.
  If $j$ is acyclic, then so are $A_1 \times_O j$ and $A_2 \times_O j$ since $A_1$ and $A_2$ are sharp over $O$; hence the pushout of the former is again an acyclic cofibration, and so by 2-out-of-3 $i \mathbin{\widehat{\times_O}} j$ is also acyclic.

  Since identity maps are sharp, and the inclusion of the identities $O\to A$ of an internal category is monic, any internal category with sharp target-morphism $A\to O$ induces a \qucoft monad on $\E/O$.
\end{eg}

\begin{thm}\label{thm:cobar-injfib}
  Let \N be a simplicial model category, \M a simplicial category, and $F : \N \toot \M : U$ a simplicial adjunction such that $U$ preserves pushouts
  and $U F$ is \qucoft.
  If $g:A\to B$ is a $U$-cofibration in \M, then
  \[\bar(F,UF,Ug):\bar(F,UF,UA)\to \bar(F,UF,UB)\]
  is a projective cofibration, and a projective acyclic cofibration if $g$ is $U$-acyclic.
\end{thm}
\begin{proof}
  We will consider the $U$-cofibration case; the acyclic one is entirely analogous.
  Thus, let $g$ be a $U$-cofibration; we must show that $\bar(F,UF,Ug)$ has the left lifting property against any $U$-acyclic $U$-fibration $p:X\to Y$.
  Since geometric realization is left adjoint to the ``total singular object'' functor $Y \mapsto \Delta[\bullet] \pow Y$, this means we must lift in any square
  \[
    \begin{tikzcd}
      \sbar(F,UF,UA) \ar[r]\ar[d] & \Delta[\bullet] \pow X \ar[d] \\
      \sbar(F,UF,UB) \ar[r] & \Delta[\bullet] \pow Y
    \end{tikzcd}
  \]
  in $\pr\dDelta\M$.
  Using the Reedy structure on $\dDelta$, this means we must inductively find lifts in squares of the form
  \[\small
    \begin{tikzcd}
      \bar_n(F,UF,UA) \amalg_{L_n\bar(F,UF,UA)} {L_n\bar(F,UF,UB)} \ar[r]\ar[d] & \Delta[n] \pow X \ar[d] \\
      \bar_{n}(F,UF,UB) \ar[r] & (\Delta[n] \pow Y) \times_{(\partial\Delta[n] \pow Y)} (\partial\Delta[n] \pow X).
    \end{tikzcd}
  \]
  in \M.
  Now, none of the degeneracy maps involved in building the colimits on the left-hand side of this square involve the outermost copy of $F$ in
  \[ \bar_n(F,UF,UZ) = \overbrace{(FU)\dotsm(FU)}^{n+1} Z = F \overbrace{(UF)\dotsm(UF)}^{n} U Z. \]
  Since $F$ preserves colimits, this left-hand morphism is therefore $F$ applied to an analogous colimit construction.
  Specifically, let $\ddup$ be the category of degeneracy maps in $\dDelta\op$; then the restriction of $\sbar(F,UF,UZ)$ to $\ddup$ can be written as $F$ applied pointwise to a diagram $\sbarp(UZ) \in \func{\ddup}{\N}$, where
  \[\barp_n(W) = \overbrace{(UF)\dotsm(UF)}^{n} W.\]
  Thus, since $F$ preserves colimits, by the adjunction $F\adj U$ it suffices to find a lift in the adjunct square
  \[\small
    \begin{tikzcd}
      \barp_n(UA) \amalg_{L_n \barp(UA)} {L_n \barp(UB)} \ar[r]\ar[d] & \Delta[n] \pow U X \ar[d] \\
      \barp_n(UB) \ar[r] & (\Delta[n] \pow U Y) \times_{(\partial\Delta[n] \pow U Y)} (\partial\Delta[n] \pow U X).
    \end{tikzcd}
  \]
  where on the right we have also used the fact that $U$ preserves pullbacks and simplicial powers, since it is a simplicial right adjoint.
  And this right-hand morphism is an acyclic fibration, since $U p: U X \to U Y$ is an acyclic fibration by assumption and \N is a simplicial model category; thus it will suffice to show the left-hand map is a cofibration in \N.

  Since $U g$ is a cofibration in \N, it will suffice to show that the functor $\sbarp : \N \to \func{\ddup}{\N}$ takes cofibrations to Reedy cofibrations.
  Inspecting the definition, this is equivalent to saying that each latching map transformation $L_n\barp \to \barp_n$ is a \qcof between endofunctors of $\N$, i.e.\ that $\barp$ is ``Reedy \qcoft''.
  At this point the argument essentially reduces to a classical one showing that bar constructions are Reedy cofibrant (a.k.a.\ ``good'' or ``proper'') when the ``unit maps'' are cofibrations (e.g.~\cite[Proposition A.10]{may:goils} and~\cite[Proposition 23.6]{shulman:htpylim}); again we simply package it abstractly.

  We can express $\ddup$ as the category of finite ordinals with distinct endpoints and injective endpoint-preserving monotone maps.
  But injectivity means we can also discard the endpoints and identify $\ddup$ with the category of finite (possibly empty) ordinals and injective monotone maps.
  In particular, the reduced slice category $\ddup\sslash [n]$, over which the colimit defining $L_n \barp$ is taken, is isomorphic to the poset of proper subsets of an $n$-element set.

  Let $\cA$ be the category $(a \ot b \to c)$, and regard $[n+1]\in\ddup$ as the $(n+1)$-element ordinal $\{\underline{0},\dotsc,\underline{n}\}$.
  There is a functor $q:(\ddup\sslash [n+1]) \to \cA$ which sends the subset $\{\underline{1},\dotsc,\underline{n}\}$ to $a$, all of its proper subsets to $b$, and all other proper subsets of $[n+1]$ (those containing $\underline{0}$) to $c$.
  This functor is an opfibration, so that colimits over $\ddup\sslash [n+1]$ (such as that defining $L_{n+1} \barp$) can be computed by first taking colimits over the fibers of $q$ and then a colimit over \cA (i.e.\ a pushout).

  The induced functor $q^{-1}(b) \to q^{-1}(c)$ is an isomorphism (given by adding $\underline{0}$ to a subset), and these fibers $q^{-1}(b)$ and $q^{-1}(c)$ are both isomorphic to $\ddup\sslash[n]$.
  Moreover, when the diagram whose colimit is $L_{n+1} \barp$ is restricted to these fibers, it becomes the diagrams whose colimits are $L_n \barp$ and $UF \circ L_n\barp$ respectively.
  Thus we obtain an expression of $L_{n+1} \barp$ as a pushout:
  \[
    \begin{tikzcd}
      L_n \barp \ar[r] \ar[d] \ar[dr,phantom,near end,"\ulcorner"] & UF \circ L_n \barp  \ar[d] \ar[ddr] \\
      \barp_n \ar[r] \ar[drr] & L_{n+1} \barp \ar[dr,dashed] \\
      && (UF\circ \barp_{n}) \mathrlap{\;\cong \barp_{n+1}}
    \end{tikzcd}
  \]
  which identifies the latching map $L_{n+1}\barp \to \barp_{n+1}$ as a Leibniz composite of the latching map $L_n \barp \to \barp_n$ with the unit $\Id \to UF$.
  Since by assumption all the functors preserve pushouts and $UF$ is \qucoft, it follows by \cref{thm:leibcomp} and induction on $n$ that each map $L_n\barp \to \barp_n$ is a \qcof.
  (The base case $L_0\barp \to \barp_0$ is $\emptyset \to \Id$, a \qcof as remarked above.)
\end{proof}

\begin{cor}
  Let \N be a simplicial model category, \M a simplicial category, and $F : \N \toot \M : U$ a simplicial adjunction such that $U$ preserves pushouts and geometric realizations and $U F$ is \qucoft.
  If $X\in\M$ is $U$-cofibrant, then the augmentation $\bar(F,U F,U X) \to X$ is a projective cofibrant replacement (i.e.\ a weak equivalence from a projective cofibrant object).\qed
\end{cor}

Since our main interest is in injective model structures, we henceforth dualize everything.
Thus, a forgetful simplicial functor $U:\M\to\N$ preserving totalizations (the dual of geometric realization) and having a simplicial right adjoint $G$ yields a \textbf{cobar construction} $\cobar(G,UG,U-)$ with a natural coaugmentation $\nu_X:X \to \cobar(G,UG,UX)$ such that $U\nu_X$ is a simplicial strong deformation coretraction (hence $\nu_X$ is a $U$-weak equivalence).
If $U$ also preserves pullbacks and the counit $U G \to \Id$ is a \qfib, then $\cobar(G,UG,U-)$ takes $U$-fibrations to injective fibrations and $U$-acyclic $U$-fibrations to injective acyclic fibrations.
In particular, if $X$ is $U$-fibrant, then $\cobar(G,UG,UX)$ is injectively fibrant; so if \M has a $U$-fibrant replacement functor $R$, we obtain an injective fibrant replacement $X \to R X \to \cobar(G,UG,URX)$.

Importantly, we can also extend this to a factorization of morphisms.

\begin{defn}\label{defn:rel-cobar}
  Given a morphism $f:X\to Y$ in \M, its \textbf{relative cobar construction} is the pullback $E f$ shown below:
  \begin{equation}
    \begin{tikzcd}
      X \ar[dr,"{\lambda_f}" description] \ar[drr,"{\nu_X}"] \ar[ddr,"f"'] \\
      & E f \ar[r,"{\nu_f}"'] \ar[d,"{\rho_f}"] \ar[dr,phantom,"\lrcorner"] & \cobar(G,UG,UX) \ar[d,"{\cobar(G,UG,Uf)}"]\\
      & Y \ar[r,"{\nu_Y}"'] & \cobar(G,UG,UY)
    \end{tikzcd}\label{eq:Ef}
  \end{equation}
\end{defn}

We call the functorial factorization $f = \rho_f \circ \lambda_f$ the \textbf{cobar factorization}.

\begin{lem}\label{thm:cobarff-cart}
  If $U:\M\to\N$ preserves pullbacks and has a right adjoint $G$, then the cobar functorial factorization is cartesian in the sense of \cref{eg:cart-ff-local}.
\end{lem}
\begin{proof}
  The cobar construction $\cobar(G,UG,U-)$ preserves pullbacks, since $U$ and $G$ do, as does totalization of cosimplicial objects (being a limit construction).
  Thus any pullback square
  \[
    \begin{tikzcd}
      X' \ar[r] \ar[d,"f'"'] \ar[dr,phantom,near start,"\lrcorner"] & X \ar[d,"f"]\\
      Y' \ar[r] & Y
    \end{tikzcd}
  \]
  gives rise to a commutative cube:
  \[\begin{tikzcd}[row sep=small]
      E f' \arrow[dd, "{\rho_{f'}}"'] \arrow[rr] \arrow[rd] &  & {\cobar(G,UG,UX')} \arrow[dd] \arrow[rd] &  \\
      & E f \arrow[rr,crossing over] &  & {\cobar(G,UG,UX)} \arrow[dd] \\
      Y' \arrow[rr] \arrow[rd] &  & {\cobar(G,UG,UY')} \arrow[rd] &  \\
      & Y \arrow[rr] \arrow[from=uu, near start, "{\rho_f}"', crossing over] &  & {\cobar(G,UG,UY)}
    \end{tikzcd}\]
  in which the right-hand face as well as the front and back faces are pullbacks.
  Thus the left-hand face is also a pullback, i.e.\ $E$ is cartesian.
\end{proof}

\begin{lem}\label{thm:rho-injfib}
  Suppose $U:\M\to \N$ preserves pullbacks and has a simplicial right adjoint $G$ such that $U G$ is \qufibt.
  Then whenever $f$ is a $U$-fibration, $\rho_f$ is an injective fibration.
\end{lem}
\begin{proof}
  It is a pullback of $\cobar(G,UG,Uf)$, which is an injective fibration by the dual of \cref{thm:cobar-injfib}.
\end{proof}

\begin{lem}\label{thm:lambda-defret}
  If $U:\M\to\N$ has a simplicial right adjoint $G$ and preserves totalizations and pullbacks, then for any $f$ with cobar factorization $f = \rho_f \lambda_f$, the map $U\lambda_f$ in \N is the inclusion of a simplicial strong deformation retract over $U Y$.
\end{lem}
\begin{proof}
  The strong deformation retraction of \cref{thm:bar-repl} is natural (on functors with codomain \N, not \M).
\begin{concise}
  Thus, these retractions for $U\nu_Y$ and $U\nu_X$ induce one for $U\lambda_f$, using the enriched universal property of the pullback $U (E f)$ (which is also a pullback since $U$ preserves pullbacks).
\end{concise}
\begin{verbose}
  \fxwarning{Would need updating to general adjunctions.}
  That is, for fixed $d\in \D$, we have
  \begin{mathpar}
    r_X(d) : C(\D(d,-),\D,X) \to X(d) \and r_{Y}(d) : C(\D(d,-),\D,Y) \to Y(d)
  \end{mathpar}
  such that $f \circ r_X(d) = r_Y(d) \circ C(\D(d,-),\D,f)$, and simplicial homotopies $H_X(d)$ and $H_{Y}(d)$ from identities to $r_X(d) \circ \nu_X(d)$  and $r_Y(d) \circ \nu_Y(d)$ respectively, such that $C(\D(d,-),\D,f) \circ H_X(d) = H_Y(d) \circ C(\D(d,-),\D,f)$, and such that $H_X(d) \circ \nu_X(d)$ and $H_Y(d) \circ \nu_Y(d)$ are constant homotopies.

  Now the composite $r_X(d) \circ \nu_f : E f \to X$ is over $Y$ and satisfies $r_X(d) \circ \nu_f \circ \lambda_f = r_X(d) \circ \nu_X = 1_X$, so it suffices to give a homotopy from the identity of $E f$ to $\lambda_f \circ r_X(d) \circ \nu_f$.
  And by the universal property of the pullback $E f$, for this it suffices to give homotopies
  \begin{mathpar}
    \nu_f \sim \nu_f \circ \lambda_f \circ r_X(d) \circ \nu_f
    \and \rho_f \sim \rho_f \circ \lambda_f \circ r_X(d) \circ \nu_f
  \end{mathpar}
  that agree in $C(\D(d,-),\D,Y)$.
  But
  \[ \nu_f \circ \lambda_f \circ r_X(d) \circ \nu_f = \nu_X \circ r_X(d) \circ \nu_f \sim \nu_f\]
  by $H_X(d)\circ \nu_f$, while
  \begin{multline*}
  \rho_f \circ \lambda_f \circ r_X(d)\circ \nu_f
    = f \circ r_X(d)\circ \nu_f
  \\  = r_Y(d) \circ C(\D(d,-),\D,f) \circ \nu_f
    = r_Y(d) \circ \nu_Y \circ \rho_f
    = \rho_f
  \end{multline*}
  so we have the constant homotopy.
  And to see that these agree in $C(\D(d,-),\D,Y)$, we note that
  \[ C(\D(d,-),\D,f) \circ H_X(d)\circ \nu_f
    = H_Y(d) \circ C(\D(d,-),\D,f)\circ \nu_f
    = H_Y(d) \circ \nu_Y \circ \rho_f \]
  which is constant since $H_Y(d) \circ \nu_Y$ is constant.
  Finally, this homotopy from the identity of $E f$ to $\lambda_f \circ r_X(d) \circ \nu_f$ becomes constant when precomposed with $\lambda_f$ since $H_X(d)$ becomes so when composed with $\nu_X = \nu_f \circ \lambda_f$.
\end{verbose}
\end{proof}

We can now deduce our desired characterization of injective fibrations.

\begin{thm}\label{thm:injfib}
  Let $U:\M\to\N$ be a simplicial functor with both simplicial adjoints $F\adj U\adj G$, where \N is a simplicial model category, and suppose that:
  \begin{enumerate}
  \item $U G$ is \qufibt (equivalently, $U F$ is \qucoft).\label{item:im1}
  \item Inclusions of simplicial strong deformation retracts are cofibrations in \N (hence automatically acyclic cofibrations).\label{item:im3}
  \end{enumerate}
  Then a map $f:X\to Y$ in \M is an injective fibration if and only if it is a $U$-fibration and the map $\lambda_f : X\to E f$ has a retraction over $Y$.

  In particular, an object $X$ is injectively fibrant if and only if it is $U$-fibrant and the map $\nu_X : X \to \cobar(G,UG,UX)$ has a retraction.
\end{thm}
\begin{proof}
  Since $U$ is a simplicial right adjoint, it preserves totalizations and pullbacks, so \cref{thm:cobar-injfib,thm:lambda-defret} apply.
  In one direction, if $f$ is a $U$-fibration, then by \cref{thm:rho-injfib} $\rho_f$ is an injective fibration.
  Thus, if $f$ is a retract of it, then $f$ is also an injective fibration.

  Conversely, suppose $f:X\to Y$ is an injective fibration.
  Since~\ref{item:im1} implies that $U G$ is also \qfibt, the argument of \cref{thm:projmodel-gen} tells us that $f$ is a $U$-fibration.
  Moreover, by \cref{thm:lambda-defret} and~\ref{item:im3} $\lambda_f$ is a $U$-acyclic $U$-cofibration; thus $f$ has right lifting for it, yielding the desired retraction.
\end{proof}

Note that \cref{thm:injfib}\ref{item:im3} is automatically satisfied if the cofibrations in \N are the monomorphisms.
If in addition \N is right proper, an alternative proof of \cref{thm:lambda-defret} is to observe that the composite of $U\lambda_f$ with the weak equivalence $U\nu_f$ is the acyclic cofibration $U\nu_X$, hence it is an acyclic cofibration.

\begin{concise}
\begin{eg}\label{eg:cat}
  The dual of \cref{thm:injfib}\ref{item:im3} seems harder to satisfy in general, but one example where it holds is the ``trivial model structure''~\cite{lack:htpy-2monads} on a 2-category.
  In particular, if in \cref{eg:mat} we let $\V=\cCat$ with its canonical model structure and \E be a 2-category with its trivial model structure, then for any small 2-category \D, the dual of \cref{thm:injfib} applies to the projective model structure on the category $\pr\D\E$ of 2-functors and strict 2-natural transformations.
  In this case, as mentioned in \cref{rmk:2mnd,rmk:flexible}, the simplicial bar construction coincides with the \emph{codescent data} of~\cite{lack:codescent-coh}, and its realization with the corresponding \emph{pseudomorphism classifier} $\pscl A$.
  Thus we recover~\cite[Theorem 4.12]{lack:htpy-2monads} for the projective model structure: the cofibrant presheaves are the flexible ones.

  Dually, when \E has its ``cotrivial'' model structure, the injectively fibrant objects are the \emph{coflexible} presheaves: those that are retracts of their pseudomorphism coclassifier $\cocl A$, whose universal property is that strict 2-natural transformations $B \to \cocl A$ are bijective to pseudonatural ones $B\cohto A$.
  Explicitly, an element of $\cocl A_x$ is an element $a\in A_x$ together with, for each nonidentity $\xi:y\to x$ in \D, an object $a_\xi\in A_{y}$ and an isomorphism $a_\xi \cong \xi^*(a)$.
  (This construction is familiar to type theorists as the ``right adjoint splitting'' of a comprehension category~\cite{hofmann:ttinlccc,lw:localuniv}.)
  Thus, a presheaf $A$ is injectively fibrant if any such ``section with specified pseudo-restrictions'' can be ``rectified'' to a single section, in a way that is strictly natural and fixes the image of $A\to \cocl A$ (which is defined by taking $a_\xi=\xi^*(a)$ for all $\xi$).
  A similar explicit interpretation can be given for arbitrary \V. 
\end{eg}
\end{concise}

\begin{concise}
\begin{rmk}\label{rmk:reedy}
  When \D is a direct category, a generalized direct category~\cite{bm:extn-reedy}, or an elegant Reedy category~\cite{br:reedy}, the injective model structure coincides with the Reedy one.
  It is an interesting exercise to work out explicitly how in such cases \cref{thm:injfib} ends up describing the Reedy fibrations.
  Intuitively, in a Reedy fibration we can \emph{inductively} construct ``rectifications'' for families of pseudo-restrictions as described above, by successively applying the path lifting property for each matching object fibration $A_x \fib M_x A$.
  That is, in these special cases the ``natural rectifiability'' property of \cref{thm:injfib} can be reduced to a family of \emph{non-interacting} path lifting properties, one for each object of \D; but in the general case the naturality must be included in the operation.

  The only explicit characterization of a non-Reedy injective model structure that I am aware of is~\cite{bordg:thesis,bordg:injective} for \D the 2-element group.
  We leave it to the interested reader to relate this characterization directly to \cref{thm:injfib}.
\end{rmk}
\end{concise}

\begin{verbose}
\begin{eg}
  The injective fibrations in $\pr\D\E$ have previously known characterizations in special cases, such as when \D is a direct category, a generalized direct category~\cite{bm:extn-reedy}, an elegant Reedy category~\cite{br:reedy}, or the 2-element group~\cite{bordg:thesis,bordg:injective} (in all cases but the last, the characterization being that the injective model structure coincides with the Reedy one).
  Of course these characterizations must be equivalent to ours, but it is an interesting exercise to work out explicitly how that occurs.

  For example, consider the simplest case, when \D is the direct category $\dtwo = (0\to 1)$; the Reedy approach says that a presheaf $A\in\E^\dtwo$ is injectively fibrant if $A_0$ is fibrant and the map $A_1\to A_0$ is a fibration (hence $A_1$ is also fibrant).
  By contrast, \cref{thm:injfib} starts by requiring that $A_0$ and $A_1$ are separately fibrant and then imposes another condition.
  The simplicial cobar construction
  \[\scobar(G,U G,U A) = (\scobar(G,U G,U A)_1 \to \scobar(G,U G,U A)_0)\]
  has $\scobar(G,U G,U A)_0$ constant at $A_0$, while $\scobar(G,U G,U A)_1$ is the cosimplicial object
  \[
    \begin{tikzcd}[column sep=huge]
      A_0 \times A_1 \ar[r,shift left=2] \ar[r,shift right=2] \ar[r,<-] &
      A_0 \times A_0 \times A_1 \ar[r,shift left=4] \ar[r,shift right=4] \ar[r] \ar[r,shift left=2,<-] \ar[r,shift right=2,<-] &
      A_0 \times A_0 \times A_0 \times A_1 \mathrlap{\; \cdots}
    \end{tikzcd}
  \]
  in which the leftwards-pointing codegeneracy arrows discard factors of $A_0$ in the middle, and each rightwards-pointing arrow duplicates a different copy of $A_0$, except for the last one in each group which duplicates $A_1$ and then applies the given map $A_1\to A_0$.
  The totalization of this involves, at the beginning, $A_0\times A_1$ together with a simplicial homotopy between its two images in $A_0\times A_0\times A_1$; but this homotopy must become constant when the middle copy of $A_0$ is projected away, so it is really just a homotopy between the projection and the action $A_0\times A_1\toto A_0$.
  For similar reasons, the higher simplices are forced to be entirely degenerate, so that $\cobar(G,UG,UA)_1$ is the pullback
  \[
    \begin{tikzcd}
      \cobar(G,UG,UA)_1 \ar[d] \ar[r] \drpullback & \Delta[1] \pow A_0\ar[d] \\
      A_0\times A_1 \ar[r] & A_0\times A_0,
    \end{tikzcd}
  \]
  i.e.\ the mapping path space of the morphism $A_1\to A_0$.
  Since $A_0$ and $A_1$ are fibrant, the projection $\cobar(G,UG,UA)_1 \to \cobar(G,UG,UA)_0 = A_0$ is indeed a fibration; thus if $A$ is a retract of $\cobar(G,UG,UA)$ then $A_1\to A_0$ must also be a fibration.
  The converse is just the path lifting property of a fibration.
\end{eg}
\end{verbose}

\begin{verbose}
We end this section by considering some special cases of \cref{thm:injfib} which reproduce known results.

\begin{eg}\label{eg:cat}
  Let $\V=\cCat$ with its canonical model structure, in which the weak equivalences are the equivalences of categories, the cofibrations are the injective-on-objects functors, and the fibrations (``isofibrations'') are the functors with the isomorphism-lifting property.
  Let \E be a 2-category with its \emph{trivial model structure} in the sense of~\cite{lack:htpy-2monads}, whose weak equivalences are the internal equivalences and whose fibrations are the internal isofibrations.
  The factorizations in the latter are obtained from 2-categorical limits and colimits, hence are 2-functorial; thus \cref{thm:injmodel}\ref{item:im2} holds, while~\ref{item:im1} is clear.
  The dual of condition~\ref{item:im3} is straightforward to check, so \cref{thm:injmodel} yields the projective model structure on $\pr\D\E$.
  Note that $\pr\D\E$ is 2-monadic over $\E^{\ob\D}$ via the 2-monad
  \begin{equation}
    T A(x) = \coprod_{y\in\ob\D} \D(x,y)\times A(y),\label{eq:2mnd}
  \end{equation}
  and the projective model structure coincides with the one constructed in~\cite{lack:htpy-2monads} for algebras over a 2-monad.

  Now the geometric realization of a simplicial object in a 2-category coincides with the ``codescent object'' of its 2-truncation.
  In the case of $B(\D,\D,A)$, this codescent object is precisely the one used in~\cite{lack:codescent-coh} to construct the \emph{pseudo morphism classifier} $\pscl{A}$ for the 2-monad $T$, i.e.\ a left adjoint to the inclusion $T\algs \to T\alg$ of algebras and strict morphisms into algebras and pseudo morphisms.
  Thus \cref{thm:injfib} recovers the result of~\cite[Theorem 4.12]{lack:htpy-2monads} for the projective model structure: the cofibrant objects are those for which the map $q:\pscl{A}\to A$ has a section in $T\algs$, known in 2-category theory as \textbf{flexible algebras}.

  Dually, starting from the ``cotrivial model structure'' on \E, \cref{thm:injmodel} yields the injective model structure.
  The cobar construction $C(\D,\D,A)$ specializes to the descent object that constructs a pseudo morphism \emph{coclassifier} $\cocl{A}$ for the 2-\emph{comonad} right adjoint to $T$, whose coalgebras and their strict and pseudo morphisms coincide with the $T$-algebras and their morphisms.
  Now \cref{thm:injfib} yields a dual of~\cite[Theorem 4.12]{lack:htpy-2monads}: the fibrant objects in the injective model structure are the \textbf{coflexible coalgebras}, whose inclusion into their pseudo morphism coclassifier has a retraction.
  (Note that the construction of a model structure on $T\algs$ in~\cite{lack:htpy-2monads} does not dualize in general, as it appeals to local presentability.)

  Now suppose $\D$ is a 1-category and $\E=\cCat$; then we can describe $\pscl{A}$ and $\cocl{A}$ more explicitly.
  By~\cite[Theorem 4.10]{lack:codescent-coh}, since the 2-monad $T$ preserves bijective-on-objects morphisms, the pseudo morphism classifier $\pscl{A}$ can also be constructed by factoring the action map $T A \to A$ as a bijective-on-objects morphism followed by a fully-faithful one.
  Thus, the objects of $\pscl{A}(x)$ are those of $B_0(\D,\D,A) = T A$, i.e.\ pairs consisting of a morphism $\xi :x\to y$ in $\D$ and an object $a\in A(y)$; while its morphisms are the morphisms in $A$ between their images, i.e.\ a morphism $(\xi,a) \to (\ze,b)$ is a morphism $\xi^*(a) \to \ze^*(b)$ in $A(x)$.

  For the pseudo morphism coclassifier $\cocl{A}$ we can inspect its definition as a descent object to see that an object of $\cocl{A}(x)$ consists of the following:
  \begin{enumerate}[label=(\arabic*)]
  \item For every morphism $\xi:y\to x$ in \D, an object $a_\xi \in A(y)$.
  \item For every composable pair $z \xto{\ze} y \xto{\xi} x$ in \D, an isomorphism $\ze^*(a_\xi) \cong a_{\xi\ze}$.
  \item When $\ze=1_y$, this isomorphism is the identity of $a_\xi$.
  \item For every composable triple $w\xto{\chi} z \xto{\ze} y \xto{\xi} x$ in \D, the following square commutes:
    \[
      \begin{tikzcd}
        \chi^* \ze^*(a_\xi) \ar[r] \ar[d, equals] & \chi^*(a_{\xi\ze}) \ar[d]\\
        (\ze\chi)^*(a_\xi) \ar[r] & a_{\xi\ze\chi}.
      \end{tikzcd}
    \]
  \end{enumerate}
  We leave it to the reader to describe the morphisms in $\cocl{A}(x)$.
  The action of a morphism $\chi : w\to x$ in \D is the functor $\cocl{A}(x) \to \cocl{A}(w)$ sending the family of objects $(a_\xi)_{\xi : y\to x}$ to the family $(a_{\chi\xi})_{\xi:y\to w}$, with corresponding family of isomorphisms.

  However, this can be further simplified: taking $\xi=1_x$ in the above square shows that the isomorphisms $\ze^*(a_\xi) \cong a_{\xi\ze}$ are uniquely determined by the cases when $\xi=1_x$, in which case the squares impose no further condition.
  Thus, an object of $\cocl{A}(x)$ can be more simply described as:
  \begin{enumerate}[label=(\arabic*$'$)]
  \item For every morphism $\xi:y\to x$ in \D, an object $a_\xi \in A(y)$.
  \item For every morphism $\xi:y\to x$ in \D, an isomorphism $\xi^*(a_1) \cong a_\xi$.
  \item When $\xi=1_y$, this isomorphism is the identity of $a_1$.
  \end{enumerate}
  The presheaf $A$ is then injectively fibrant, i.e.\ coflexible, if the evident inclusion $A\to \cocl{A}$ has a retraction.
  Note that the obvious family of retractions $\cocl{A}(x) \to A(x)$ picking out $a_{1}$ is only \emph{pseudonatural} in $x$; for $A$ to be coflexible there must be a \emph{strictly} natural retraction.

  Note that these two constructions $\pscl{A}$ and $\cocl{A}$ are familiar to categorical type theorists under another name: they are the two ways to replace a comprehension category by a split one, called the \emph{left and right adjoint splittings} in~\cite{lw:localuniv}.
\end{eg}

\begin{eg}
  The primary special case in which an explicit description of injective fibrations is known is when \D is a \emph{direct category}, i.e.\ admits an identity-reflecting ``degree'' functor to the partially ordered class of ordinals.
  In this case the injective model structure coincides with the Reedy model structure, hence the injective fibrations coincide with the Reedy fibrations.
  The most direct proof of this is to simply observe that the two model structures have the same weak equivalences and cofibrations, but this conveys little insight into why it happens and how Reedy fibrations are related to injective ones in general.

  The underlying idea is clearest in the simple case $\E=\cCat$ described above.
  When \D is a direct category, to say that $A\in \pr\D\cCat$ is Reedy fibrant is to say that each matching map $p_x:A_x\to M_x A$ is an isofibration, i.e.\ that given $a\in A_x$ and $b\in M_x A$ and $\phi : p_x(a) \cong b$, there is an $\ahat \in A_x$ and $\phihat : a\cong \ahat$ with $p_x(\ahat) = b$ and $p_x(\phihat) = \phi$.
  Now $b\in M_x A$ consists of, for each nonidentity morphism $\xi : y\to x$ in \D, an object $b_\xi \in A(y)$, such that for any $\ze : z\to y$ we have $\ze^*(b_\xi) = b_{\xi\ze}$.
  Similarly, $\phi : p_x(a) \cong b$ consists of, for each nonidentity $\xi:y\to x$, an isomorphism $\phi_\xi : \xi^*(a) \cong b_\xi$, such that for any $\ze:z\to y$ we have $\ze^*(\phi_\xi) = \phi_{\xi\ze}$.
  But given these data, the objects $b_1 \coloneqq a$ and $b_\xi$ with the equalities $\ze^*(b_\xi) = b_{\xi\ze}$ and isomorphisms $\phi_\xi : \xi^*(a) \cong b_\xi$ can be regarded as an object of $\cocl{A}_x$.
  Thus, if $A$ is coflexible, the retraction $\cocl{A} \to A$ supplies an object $\ahat \in A_x$.
  Moreover, the isomorphisms $\phi$ also exhibit an isomorphism between this object and the image of $a$ in $\cocl{A}$, so the retraction also yields $\phihat : a\cong \ahat$.

  Conversely, if $A$ is Reedy fibrant and we are given an object of $\cocl{A}_x$, by inducting up the degree function we can successively transfer its objects along its isomorphisms to end up with an object of $A_x$ using the isofibration property of each matching map.
  For instance, when $\D = (z\xto{\ze} y\xto{\xi} x)$, an object of $\cocl{A}_x$ consists of objects $a_1\in A_x$, $a_\xi \in A_y$, and $a_{\xi\ze} \in A_z$ with isomorphisms $\xi^*(a_1) \cong a_\xi$ and $\ze^*\xi^*(a_1) \cong a_{\xi\ze}$.
  We can then lift $a_{\xi\ze}$ along the composite isomorphism $\ze^*(a_\xi) \cong \ze^*\xi^*(a_1) \cong a_{\xi\ze}$ to an object $\ahat_{\xi}$ with an isomorphism $a_\xi \cong \ahat_\xi$, and lift $\ahat_\xi$ along the composite isomorphism $\xi^*(a_1) \cong a_\xi \cong \ahat_\xi$ to an object $\ahat_1$ with an isomorphism $a_1\cong \ahat_1$, and define the retraction to be $\ahat_1\in A_x$.
  By choosing lifts of identity morphisms to be identities, this becomes a strict retraction, and by the careful choice of degreewise lifts it is strictly natural.
  (The most natural way to formulate this precisely for a general direct category \D is probably to construct the map $\cocl{A} \to A$ degree by degree in the usual Reedy way, using the lifting property of the matching isofibration against the levelwise acyclic cofibration $A_x \to \cocl{A}_x$ at each step.)

  A similar, but more complicated, analysis is possible for more general \E, as well as when \D is one of the more general \emph{elegant Reedy categories} of~\cite{br:reedy} which also have the property that their Reedy model structures coincide with injective ones (at least for some model categories \E).
\end{eg}

\begin{eg}
  The only other previously known explicit characterization of injective fibrations I am aware of was given in~\cite{bordg:thesis} in the case $\E=\cGpd$ and $\D$ the group $\dZ/2$ regarded as a 1-object groupoid.
  Thus, the objects of $\pr\D\E$ are groupoids equipped with a (strict) involution.
  Proposition 5.3.5 of~\cite{bordg:thesis} states that a map $f:(A,\alpha)\to (B,\beta)$ of groupoids with involution is an injective fibration if and only if it is an isofibration on underlying groupoids (i.e.\ a projective fibration), and given any $a\in A$ equipped with an isomorphism $\phi : \alpha(a) \cong a$ from its involute such that $\alpha(\phi) = \phi^{-1}$, and isomorphism $\psi : f(a) \cong b$ such that $\beta(b) = b$ and $\beta(\psi) = \psi\circ f(\phi)$, there exists an object $\ahat\in A$ with $\alpha(\ahat)=\ahat$ and $f(\ahat) = b$ and an isomorphism $\phihat:a\cong \ahat$ such that $\alpha(\phihat) = \phihat\circ \phi$ and $f(\phihat)=\psi$.
  This is certainly a lifting property against an injective acyclic cofibration, so any injective fibration must have this property; the hard part is showing that it is sufficient.
  In~\cite{bordg:thesis} this is done by decomposing any injective acyclic cofibration as a cell complex in an explicit way; we can instead derive it from \cref{thm:injfib} as follows.

  Given $(A,\alpha)$, an object of $C(\D,\D,A) = \cocl{A}$ is, by the explicit description in \cref{eg:cat}, two objects $a,a'\in A$ equipped with an isomorphism $\phi : \alpha(a) \cong a'$.
  The morphisms are pairs of morphisms in $A$ that commute with the isomorphisms, and the involution $\cocl{\alpha}$ maps $(a,a',\phi)$ to $(a',a,\alpha(\phi)^{-1})$.
  The injection $A\to \cocl{A}$ sends each object $a$ to the pair $(a,\alpha(a))$ with $\phi = 1_{\alpha(a)}$.

  Let $f:A\to B$ be a projective fibration satisfying the single explicit lifting property described above.
  We must exhibit $A$ as a retract of $E f$ over $B$, where the objects of $E f$ are triples $(a,a',\phi) \in \cocl{A}$ as above such that $f(a)=f(a')$ and $f(\phi)$ is an identity.
  We define $E f \to A$ on objects by three cases:
  \begin{enumerate}
  \item If $\phi$ is an identity, then $(a,\alpha(a),1)$ is the image of $a\in A$, so to have a retraction we must send it back to $a$.
    This respects the involution.
  \item If $\phi$ is not an identity and $\alpha(\phi)^{-1} \neq \phi$, then the same is true of $\cocl{\alpha}(a,a',\phi) = (a',a,\alpha(\phi)^{-1})$, which is unequal to $(a,a',\phi)$.
    Thus, objects of this sort come in pairs.
    In each pair we choose one object $(a,a',\phi)$ to send to $a\in A$, and send the other object $(a',a,\alpha(\phi)^{-1})$ to $\alpha(a)$ (as we must, to respect the involution).
  \item Finally, if $\phi$ is not an identity but $\alpha(\phi)^{-1} \neq \phi$, then the object $(a,a',\phi)$ is fixed by $\cocl{\alpha}$, so we must send it to a fixed object of $A$.
    However, since $f(a)=f(a')$ and $f(\phi)$ is an identity, we can let $\psi$ be an identity in the explicit lifting property; thus there is a fixed object $\ahat$ and an isomorphism $\phihat:a\cong \ahat$ such that $\alpha(\phihat) = \phihat\circ \phi$ and $f(\phihat)$ is an identity.
    We can then send $(a,a',\phi)$ to $\ahat\in A$.
  \end{enumerate}
  We leave it to the reader to extend this definition to morphisms, by composing with the given isomorphisms as needed.
\end{eg}
\end{verbose}

As an aside, if pointwise factorizations exist we can use \cref{thm:injfib} to \emph{construct} injective model structures, not requiring local presentability.

\begin{cor}\label{thm:injmodel}
  In the situation of \cref{thm:injfib}, suppose furthermore that
  \begin{enumerate}\setcounter{enumi}{2}
  \item Every morphism in \M factors as a $U$-cofibration followed by a $U$-acyclic $U$-fibration, and as a $U$-acyclic $U$-cofibration followed by a $U$-fibration.\label{item:im2}
  \end{enumerate}
  Then the injective model structure on \M exists.
\end{cor}
\begin{proof}
  Given a map $f$ to factor, first write $f = g i$ with $i$ a $U$-acyclic $U$-cofibration and $g$ a $U$-fibration, then factor $g = \rho_g \lambda_g$ as above.
  Then $\rho_g$ is an injective fibration, while $\lambda_g$ is (by \cref{thm:lambda-defret} and~\ref{item:im3}) a $U$-acyclic $U$-cofibration, hence so is $\lambda_g i$.
  The other factorization is analogous.
\end{proof}

\begin{cor}\label{thm:injfunctor}
  Suppose \V is a symmetric monoidal simplicial model category, \M is a \V-model category, and \D is a small \V-category such that each $\D(x,y)$ is cofibrant and each $\unit \to \D(x,y)$ is a cofibration in \V.
  Suppose moreover that inclusions of simplicial strong deformation retracts are cofibrations in \M, and that one of the following holds:
  \begin{enumerate}[label=(\alph*)]
  \item \M is cofibrantly generated.\label{item:if1}
  \item \M has \V-enriched functorial factorizations.\label{item:if2}
  \item \D is the free \V-category generated by an ordinary category.\label{item:if3}
  \end{enumerate}
  Then the injective model structure on $\pr\D\M$ exists.
\end{cor}
\begin{proof}
  \cref{thm:injfib}\ref{item:im1} holds by \cref{eg:mat} and~\ref{item:im3} by assumption, so it suffices to verify \cref{thm:injmodel}\ref{item:im2}.
  In cases~\ref{item:if2} and~\ref{item:if3} we can simply apply the factorizations of \M pointwise, while in case~\ref{item:if1} the projective model structure exists and its factorizations are in particular pointwise ones.
\end{proof}

We now return to our main goal of constructing \ttmts.

\begin{thm}\label{thm:mnd-ttmt}
  Let \E be a \ttmt and $T$ a simplicially enriched \qucoft monad on \E having a simplicial right adjoint $S$.
  Then the category $\E^T$ of $T$-algebras is again a \ttmt.
\end{thm}
\begin{proof}
  Since $T$ is a simplicially enriched monad, $\E^T$ is a simplicially enriched category.
  And since $T$ has a simplicial right adjoint $S$, there is an induced simplicial comonad structure on $S$ such that the $T$-algebras coincide with the $S$-coalgebras, and thus the forgetful functor $U:\E^T\to\E$ has both simplicial adjoints.
  Hence by \cref{thm:projmodel-gen}, the injective model structure on $\E^T$ exists and is proper, combinatorial, and simplicial.
  Moreover, $S$ is a simplicially enriched comonad that preserves finite limits, and the simplicial enrichment carries through the proof in~\cite[Theorem A4.2.1]{ptj:elephant} that the category of coalgebras for a finite-limit-preserving comonad inherits local cartesian closure; thus $\E^T$ is \slcc.
  And as $S$ is a right adjoint, it is accessible and left exact; hence its category of coalgebras is again a Grothendieck topos (it is an elementary topos by~\cite[Theorem A4.2.1]{ptj:elephant}, and locally presentable by the limit theorem of~\cite{mp:accessible}).

  Finally, if \F is the \local and \stratified \nfs for the fibrations of \E, let $\F^T= U^{-1}(\F) \times_\cE \dR_E$, where $\dR_E$ is as in \cref{eg:ff-fcos} for the cobar functorial factorization $E$ of \cref{defn:rel-cobar}.
  By \cref{thm:cobarff-cart,eg:cart-ff-local}, $\dR_E$ and hence $\F^T$ are \local, and by \cref{thm:injfib,eg:cof-ff-fcos} $\F^T$ is \stratified and the morphisms with $\F^T$-structure are the injective fibrations.
\end{proof}

\begin{rmk}
  The category of algebras for a cocontinuous monad on a presheaf category is again a presheaf category.
  Thus if \E in \cref{thm:mnd-ttmt} is a presheaf topos, so is $\E^T$ (cf.~\cref{rmk:design-space}).

  In fact, however, the proof of \cref{thm:mnd-ttmt} really relies more on the comonad $S$.
  The fact that $S$ has a left adjoint $T$ is only used to show that $S$ preserves finite limits and totalizations, that injective fibrations are $U$-fibrations, and that $S$-coalgebras satisfy the acyclicity condition to lift the model structure.
  The latter property could also be proven with the dual of ``Quillen's path object argument'', while the others could be assumed explicitly of a comonad, yielding a closer analogue of~\cite[Theorem A4.2.1]{ptj:elephant} for \ttmts. 
\end{rmk}

\begin{cor}\label{thm:simppre-ttmt}
  If \E is a \ttmt and \D is a small simplicially enriched category, then the injective model structure on $\pr\D\E$ is a \ttmt.
\end{cor}
\begin{proof}
  By \cref{eg:enrcat,eg:mat}, the monad on $\E^{\ob\D}$ whose category of algebras is $\pr\D\E$ is \qucoft, and has a simplicial right adjoint.
\end{proof}

\cref{thm:simppre-ttmt} is the main result we want: combined with left exact localizations (which we discuss in the next section) it will show that \ttmts suffice to present all \io-toposes.
(Note that it includes the models of~\cite[\sect 11]{shulman:invdia} and~\cite{shulman:elreedy}.)
However, \cref{thm:mnd-ttmt} yields many other examples as well.

\begin{cor}\label{thm:intpre-ttmt}
  If \E is a \ttmt and $\cD$ is an internal category in \E whose target map $t:D_1\to D_0$ is sharp, then the category $\pr \cD\E$ of internal presheaves on $\cD$ is a \ttmt.
\end{cor}
\begin{proof}
  By \cref{eg:intcat,eg:span}, the monad on $\E/D_0$ whose category of algebras is $\pr \cD\E$ is \qucoft, and since \E is \slcc it is a simplicial monad with a simplicial right adjoint.
\end{proof}

\begin{cor}\label{thm:enrpre-ttmt}
  If \E is a \ttmt and \D is a small \E-enriched category such that every object $\D(x,y)$ is sharp, then the injective model structure on $\pr\D\E$ is a \ttmt.
\end{cor}
\begin{proof}
  Although \E may not be a cartesian monoidal \emph{model} category, monomorphisms in a topos are always preserved by Leibniz products, and acyclicity of one factor is preserved if the other has sharp domain and codomain, similarly to \cref{eg:span}.
  Thus the monad on $\E^{\ob\D}$ whose category of algebras is $\pr\D\E$ is again \qucoft, and has a simplicial right adjoint.
\end{proof}

\begin{cor}\label{thm:glue-ttmt}
  If $\E_1$ and $\E_2$ are \ttmts and $F:\E_1 \toot \E_2 : G$ is a simplicial Quillen adjunction, then the comma category $(\E_1\dn G)$ (a.k.a.\ the ``Artin gluing'') is also a \ttmt.
\end{cor}
\begin{proof}
  By \cref{thm:prod-ttmt}, $\E_1\times \E_2$ is a \ttmt.
  Define a simplicial monad $T$ on $\E_1\times \E_2$ by $T(X_1,X_2) = (X_1, F X_1 \amalg X_2)$, with right adjoint $S(Y_1,Y_2) = (Y_1 \times G Y_2, Y_2)$.
  The Leibniz application of the unit $\Id \to T$ to a pointwise cofibration $i = (i_1,i_2) : (A_1,A_2) \cof (B_1,B_2)$ is
  \[ (\id, F i_1 \amalg \id) : (B_1, F A_1 \amalg B_2) \to (B_1, F B_1 \amalg B_2 ) \]
  which is a pointwise cofibration that is acyclic if $i$ is, since $F$ is left Quillen.
  Thus $T$ is \qucoft, so \cref{thm:mnd-ttmt} applies.
\end{proof}

We can obtain the models of~\cite{shulman:eiuniv} by applying \cref{thm:glue-ttmt} iteratively with $\E_1=\S$; hence they are \ttmts.
(In fact they are also instances of \cref{thm:intpre-ttmt}.)
More generally, \emph{colax limits} of diagrams of right Quillen functors can be constructed analogously, yielding \ttmts like the models of~\cite[\sect12]{shulman:invdia}.
As in~\cite[B3.4]{ptj:elephant} and~\cite[\sect6.3.2]{lurie:higher-topoi}, these should present \emph{colimits} of \io-toposes.



%% file: intloc.tex
\section{Internal localizations}
\label{sec:int-loc}

Let \E be a \ttmt.
In the next section we will show that any left exact localization of \E is again a \ttmt; in this section we begin by studying the wider class of \emph{internal} localizations.

For an object $A\in\E$, we write $A^*$ for the pullback functor $\E \to \E/A$, taking $X$ to the product $A\times X$.
Recall that $f_*$ denotes the right adjoint of pullback $f^*$.

\begin{defn}
  Let $S$ be a class of fibrations between fibrant objects in \E.
  \begin{itemize}
  \item $S\pb$ denotes the class of all pullbacks of morphisms in $S$:
    \[ S\pb = \setof{ g^*f | (f:A\to B)\in S,\; g:C\to B } \]
  \item For $X\in\E$, $X^*S$ denotes the class of pullbacks of morphisms in $S$ to $\E/X$:
    \[ X^*S = \setof{ X^*f | f\in S} \]
  \end{itemize}
\end{defn}

Of course even when $S$ is a set, $S\pb$ is a proper class, but it is nevertheless ``generated by a set'' in the relevant sense (cf.~~\cite[Proposition 6.2.1.2]{lurie:higher-topoi}):

\begin{lem}
  For any set $S$ of fibrations between fibrant objects, there is a set of morphisms $S'$ such that an object is $S'$-local if and only if it is $S\pb$-local.
\end{lem}
\begin{proof}
  By~\cite[Proposition 4.7]{dug:pres}, there is a regular cardinal $\la$ such that all canonical maps $\displaystyle\hocolim_{\E_\la\ni W\to X}(W) \to X$ are weak equivalences, where $\E_\la$ is the full subcategory of \la-compact objects.
  Let $S'$ be a set of representatives for the pullbacks of morphisms in $S$ to \la-compact objects, and consider some $f:A\fib B$ in $S$ and some $g:X\to B$.
  Then we have a commutative square as on the left:
  \[
    \begin{tikzcd}
      \displaystyle\hocolim_{\E_\la\ni W\to X}(W \times_B A) \ar[r,"\sim"] \ar[d] &
      X\times_B A \ar[d,two heads,"g^*f"'] \ar[r] \drpullback & A \ar[d,two heads,"f"]\\
      \displaystyle\hocolim_{\E_\la\ni W\to X}(W) \ar[r,"\sim"'] & X \ar[r] & B
    \end{tikzcd}
  \]
  inducing for any fibrant object $Z$ a commutative diagram:
  \[
    \begin{tikzcd}
      \ehom\E(X,Z) \ar[r,"\sim"'] &
      \displaystyle\ehom\E\left(\hocolim_{\E_\la\ni W\to X}(W),Z\right) \ar[r,"\sim"'] &
      \displaystyle\holim_{\E_\la\ni W\to X} \ehom\E(W,Z)\\
      \ehom\E(X\times_B A, Z) \ar[r,"\sim"] \ar[from=u,"{\ehom\E(g^*f,Z)}"'] &
      \displaystyle\ehom\E\left(\hocolim_{\E_\la\ni W\to X}(W \times_B A),Z\right) \ar[r,"\sim"] \ar[from=u] &
      \displaystyle\holim_{\E_\la\ni W\to X} \ehom\E(W \times_B A,Z) \ar[from=u]
    \end{tikzcd}
  \]
  Now if $Z$ is $S'$-local, each map $\ehom\E(W,Z) \to \ehom\E(W\times_B A,Z)$ is a weak equivalence, and homotopy limits preserve weak equivalences.
  Thus the right-hand map above is a weak equivalence, hence by 2-out-of-3 so is the left-hand map; so $Z$ is $S\pb$-local.
  The converse is obvious.
\end{proof}

Since \E is combinatorial, its left Bousfield localization at $S'$, hence at $S\pb$, exists.

\begin{lem}\label{thm:stloc}
  For a fibrant object $Z$ and a class $S$ of fibrations between fibrant objects, the following are equivalent.
  \begin{enumerate}
  \item $Z$ is $S\pb$-local.\label{item:sl1}
  \item For all $f:A\fib B$ in $S$, the induced morphism
    \begin{equation}
      \ftil_Z : B^* Z \to f_* A^* Z\label{eq:stloc}
    \end{equation}
    is a weak equivalence in $\E/B$.\label{item:sl2}
  \end{enumerate}
\end{lem}
\begin{proof}
  Since $f$ is a fibration and $Z$ is fibrant, both $B^*Z$ and $f_* A^* Z$ are fibrant in $\E/B$.
  Thus~\ref{item:sl2} is equivalent to saying that for any $(f:A\to B)\in S$ and $g:X\to B$ in $\E/B$ the induced map of simplicial hom-spaces
  \[ \ehom{\E/B}(X,B^*Z) \to \ehom{\E/B}(X,f_* A^* Z) \]
  is an equivalence.
  However, by the simplicial adjunction $\E/B \toot \E$, we have $\ehom{\E/B}(X,B^*Z) \cong \ehom\E(X,Z)$.
  Similarly, since \E is \slcc, we have a simplicial adjunction $f^* : \E/B \toot \E/A : f_*$, so that
  \[\ehom{\E/B}(X,f_* A^* Z) \cong \ehom{\E/A}(f^*X, A^* Z) \cong \ehom\E(X\times_B A, Z). \]
  Under these isomorphisms~\eqref{eq:stloc} is identified with
  \[ \ehom\E(g^*f,Z) : \ehom\E(X,Z) \to \ehom\E(X\times_B A,Z) \]
  and to say that this is an equivalence for all $f\in S$ is precisely~\ref{item:sl1}.
\end{proof}

\begin{lem}\label{thm:fl-we}
  For any set $S$ of fibrations between fibrant objects and any square
  \[
    \begin{tikzcd}
      Z \ar[d,two heads] \ar[r,"\sim"] & W \ar[d,two heads] \\
      X \ar[r,"g","\sim"'] & Y
    \end{tikzcd}
  \]
  such that the vertical maps are fibrations and the horizontal maps are weak equivalences, $W$ is $(Y^*S)\pb$-local in $\E/Y$ if and only if $Z$ is $(X^*S)\pb$-local in $\E/X$.
\end{lem}
\begin{proof}
  By right properness, the induced map $Z \to g^*W$ is a weak equivalence between fibrant objects of $\E/X$.
  Since $B^*$, $A^*$, and $f_*$ in~\eqref{eq:stloc} are right Quillen functors, they preserve weak equivalences between fibrant objects.
  Thus we have
  \[
    \begin{tikzcd}
      B^*Z \ar[r] \ar[d,"\sim"'] & f_*A^* Z \ar[d,"\sim"] \\
      B^*g^*W \ar[r] & f_* A^* g^* W 
    \end{tikzcd}
  \]
  in $\E/X$, so that by 2-out-of-3 $Z$ is $(X^*S)\pb$-local if and only if $g^* W$ is.
  On the other hand, we have $g^*(Y^*S) = X^*S$, and the construction of~\eqref{eq:stloc} commutes with pullback. 
  Thus we have a triangle of pullback squares
  \[
    \begin{tikzcd}[row sep=small,column sep=small]
      B^*g^*W \ar[ddr,two heads] \ar[drr] \ar[rrr,"\sim"] &&& B^*W \ar[drr]\ar[ddr,two heads] \\
      && f_* A^* g^* W \ar[dl,two heads] \ar[rrr,near start,"\sim",crossing over] &&& f_* A^* W \ar[dl,two heads]\\
      & X \ar[rrr,"g","\sim"'] &&& Y
    \end{tikzcd}
  \]
  in which all the rightwards-pointing arrows are weak equivalences; so by 2-out-of-3 again, $g^*W$ is $(X^*S)\pb$-local if and only if $W$ is $(Y^*S)\pb$-local.
\end{proof}

\begin{prop}\label{thm:flf-nfs}
  For any set $S$ of fibrations between fibrant objects in a \ttmt \E, there is a \local, \stratified, and homotopy invariant \nfs $\dL_S$ such that the maps admitting an $\dL_S$-structure are the fibrations $Z\fib Y$ that are $(Y^*S)\pb$-local in $\E/Y$.
\end{prop}
\begin{proof}
Since $B^*Z$ and $f_* A^* Z$ in~\eqref{eq:stloc} are fibrations over $\E/B$, by~\cite[Lemma 4.3]{shulman:elreedy} 
for each $f$ and $Z$ there is a fibrant object
\( \iLoc_f Z \coloneqq B_* \isequiv_{B}(\ftil_Z) \)
such that $\ftil_Z$ is an equivalence if and only if $\iLoc_f Z$ has a global element, in which case it is acyclic (i.e.\ $\iLoc_f Z \to 1$ is an equivalence).
This is because inside homotopy type theory the property of ``being an equivalence'' is a proposition.
Hence, if we define $\iLoc_S Z \coloneqq \prod_{f\in S} \iLoc_f Z$, then $Z$ is $S\pb$-local if and only if $\iLoc_S Z$ has a global element, in which case it is also acyclic. 

Furthermore, the construction of $\iLoc_S Z$ is stable under pullback: we have $X^*(\iLoc_S Z) \cong \iLoc_{X^* S}(X^* Z)$ coherently.
Thus $\iLoc_S$ is a fibred core-endofunctor of \E, so \cref{thm:sec-afib-strat} applies; hence $\iLoc_S^*(\cEp)$ is a \local and \stratified \nfs.
  So if we define $\dL_S = \F \times_\cE \iLoc_S^*(\cEp)$, where \F is the given \nfs for fibrations in \E, then $\dL_S$ is a \local and \stratified \nfs.
  The previous paragraph, applied in $\E/Y$, shows that the $\uly{\dL_S}$-algebras over $Y$ are the $(Y^*S)\pb$-local fibrations.
  Finally, homotopy invariance follows from \cref{thm:fl-we}.
\end{proof}

\begin{rmk}\label{rmk:modal-univ}
  In particular, by \cref{thm:uf-fibrant}, $\dL_S$ admits fibrant and univalent universes whether or not $S$-localization is left exact.
  These are an external version of the ``universes of modal types'' from~\cite{rss:modalities} (for modalities only, although it should be possible to generalize them to reflective subuniverses), and are also a strict model-categorical version of the object classifiers for stable factorization systems in \io-categories obtainable from~\cite[\sect3-4]{gk:univlcc}.
\end{rmk}

We now record a few more properties of these localizations, paralleling facts proven internally in type theory in~\cite{rss:modalities}.

\begin{lem}\label{thm:pb-fl-fle}
  For any set $S$ of fibrations between fibrant objects and any $g:X\to Y$, the functor $g^*: \E/Y\to \E/X$ takes $(Y^*S)\pb$-local objects to $(X^*S)\pb$-local objects and takes $(Y^*S)\pb$-local equivalences between fibrations to $(X^*S)\pb$-local equivalences.
\end{lem}
\begin{proof}
  The first statement follows since $g^*(Y^*S) = X^*S$ and~\eqref{eq:stloc} commutes with pullback, as in \cref{thm:fl-we}.
  For the second, by factoring $g$ we can assume separately that it is a weak equivalence and that it is a fibration.

  If $g$ is a weak equivalence, then the adjunction $g_! \adj g^*$ is a Quillen equivalence by~\cite[Proposition 2.5]{rezk:proper}.
  Thus, given a $(Y^*S)\pb$-local equivalence $f:A\to B$ in $\E/Y$, to show that $\ehom{\E/X}(g^*B,Z) \to \ehom{\E/X}(g^*A,Z)$ is an equivalence for all $(X^*S)\pb$-local objects $Z\in\E/X$, it suffices to show that $\ehom{\E/Y}(B,W) \to \ehom{\E/Y}(A,W)$ is an equivalence for all fibrant objects $W\in \E/X$ such that $g^*W$ is $(X^*S)\pb$-local.
  But by \cref{thm:fl-we}, the latter is equivalent to $W$ being $(Y^*S)\pb$-local, so this is true since $f$ is a $(Y^*S)\pb$-local equivalence.

  If $g$ is a fibration, then $g^*:\E/Y\to \E/X$ is a left Quillen functor since it preserves cofibrations and weak equivalences, and $g^*(Y^*S) = X^*S$ implies $g^*((Y^*S)\pb) \subseteq (X^*S)\pb$.
  Thus, by~\cite[Proposition 3.3.18]{hirschhorn:modelcats}, $g^*$ is also a left Quillen functor from $(\E/Y)_{(Y^*S)\pb}$ to $(\E/X)_{(X^*S)\pb}$, hence (since all objects are cofibrant) it takes all $(Y^*S)\pb$-local equivalences to $(X^*S)\pb$-local equivalences.
\end{proof}

\begin{lem}[{cf.~\cite[Theorem 2.17]{rss:modalities}}]\label{thm:sigma-closed}
  For any set $S$ of fibrations between fibrant objects, any $S\pb$-local object $Y$, and any fibration $g:X\fib Y$ that is $(Y^*S)\pb$-local in $\E/Y$, the object $X$ is $S\pb$-local.
\end{lem}
\begin{proof}
  We factor the naturality square for~\eqref{eq:stloc} at $g$ through its pullback:
  \[
    \begin{tikzcd}
      \mathllap{B^*X =\;} (Y^*B)^*X \ar[r,"\sim"] \ar[d,two heads] &
      (Y^*f)_*(Y^*A)^*X \ar[r,"\sim"] \ar[d,two heads] \drpullback & f_* A^* X\ar[d,two heads]\\
      B^*Y \ar[r,equals] & B^*Y \ar[r,"\sim"'] & f_* A^* Y.
    \end{tikzcd}
  \]
  Comparing universal properties, we see that this pullback is $(Y^*f)_*(Y^*A)^*X$ as shown, and the upper-left horizontal map is~\eqref{eq:stloc} for $X$ in $\E/Y$, hence a weak equivalence.
  Since the lower-right horizontal map is~\eqref{eq:stloc} for $Y$, it is also a weak equivalence, and thus so is its pullback.
  Therefore, the top composite $B^*X \to f_* A^* X$, which is~\eqref{eq:stloc} for $X$, is also a weak equivalence; so $X$ is $S\pb$-local.
\end{proof}

\begin{lem}[{cf.~\cite[Theorem 1.32]{rss:modalities}}]\label{thm:unit-connected}
  Let $S$ be a set of fibrations between fibrant objects and let $\eta: X\to \Xhat$ be an $S\pb$-localization, i.e.\ an $S\pb$-local equivalence to an $S\pb$-local object.
  Then $\eta$ is also an $(\Xhat^*S)\pb$-local equivalence in $\E/\Xhat$.
\end{lem}
\begin{proof}
  Let $p:Z\fib \Xhat$ be an $(\Xhat^*S)\pb$-local object of $\E/\Xhat$.
  Then by \cref{thm:sigma-closed}, $Z$ is an $S\pb$-local object of \E, hence $p$ is a fibration in the $S\pb$-local model structure $\E_{S\pb}$.
  Therefore, $Z$ is a fibrant object of the slice model structure $\E_{S\pb}/\Xhat$ on $\E/\Xhat$.
  But $\eta$ is a weak equivalence (between cofibrant objects) in that same model structure, so the induced map $\ehom{\E/\Xhat}(\Xhat,Z) \to\ehom{\E/\Xhat}(X,Z)$ is a weak equivalence.
  Since this is true for any $(\Xhat^*S)\pb$-local object $Z$, $\eta$ is a $(\Xhat^*S)\pb$-local equivalence.
\end{proof}


%% file: lexloc.tex
\section{Left exact localizations}
\label{sec:lex-loc}

For a fibration $f:A\fib B$ in \E between fibrant objects, by its \textbf{fibrant diagonal} we mean a replacement of its strict diagonal $A\to A\times_B A$ by a fibration $\diag f : A' \fib A\times_B A$.
Strictly speaking this depends on the choice of fibrant replacement, but since the weak equivalence $A\to A'$ is between fibrant (and cofibrant) objects, it is a simplicial homotopy equivalence.
Thus the ambiguity is irrelevant for localization, in that whether an object is $\diag f$-local is independent of the choice of $\diag f$.

We write $\diag^n f$ for the $n$-fold iterate of $\diag$, with $\diag^0 f = f$ and $\diag^{n+1} f = \diag(\diag^n f)$.
And we write $S\dia$ for the class of iterated fibrant diagonals of morphisms in $S$:
\[ S\dia = \setof{ \diag^n f | f\in S,\; n\in \dN}.\]
By the remarks above, $X^*(S\dia)$-locality coincides with $(X^*S)\dia$-locality.
And when $S$ is a set, so is $S\dia$, so we can also localize at $(S\dia)\pb$.

We say that \textbf{$S$-localization is left exact} if it preserves homotopy pullbacks.
This implies 
that $S$-local equivalences are stable under pullback along fibrations; hence in particular a left exact localization is again right proper.
We will rely on the following characterization of left exact localizations from~\cite{abfj:lexloc}:

\begin{thm}[\cite{abfj:lexloc}]\label{thm:abfj}
  For any set $S$ of fibrations between fibrant objects:
  \begin{enumerate}
  \item $(S\dia)\pb$-localization is left exact.
  \item If $S$-localization is left exact, then it coincides with $(S\dia)\pb$-localization.
  \end{enumerate}
  In particular, every accessible left exact localization of \E is a $(S\dia)\pb$-localization.
\end{thm}

Since the $(Y^*(S\dia))\pb$-local objects of $\E/Y$ underlie a good \nfs by \cref{thm:flf-nfs}, it remains only to show that these coincide with the $(S\dia)\pb$-local fibrations over $Y$.
The idea of this proof is that the families of $(Y^*(S\dia))\pb$-localizations form an \emph{accessible lex modality} in the sense of~\cite{rss:modalities}, and the \io-categorical stable factorization system induced by any lex modality is automatically a \emph{reflective factorization system} in the sense of~\cite{chk:reflocfact}, hence coincides with the (acyclic cofibration, fibration) factorization system of the localized model structure.
We use~\cref{thm:abfj} to bridge a gap that was left open in~\cite{rss:modalities}, by relating the internal and external notions of accessibility for lex modalities.

\begin{lem}[{cf.~\cite[Theorem 3.1(vii)]{rss:modalities}}]\label{thm:conn-modal-pb}
  Let $S$ be a set of fibrations between fibrant objects, and
  suppose given a commutative square
  \[
    \begin{tikzcd}
      A \ar[r] \ar[d,two heads] & X \ar[d,two heads] \\
      B \ar[r] & Y
    \end{tikzcd}
  \]
  such that
  \begin{itemize}
  \item $X\to Y$ is $(Y^*(S\dia))\pb$-local in $\E/Y$.
  \item $A\to B$ is $(B^*(S\dia))\pb$-local in $\E/B$.
  \item $B\to Y$ is a $(Y^*(S\dia))\pb$-local equivalence in $\E/Y$.
  \item $A\to X$ is an $(X^*(S\dia))\pb$-local equivalence in $\E/B$.
  \end{itemize}
  Then the square is a homotopy pullback.
\end{lem}
\begin{proof}
  First factor the square in the Reedy fashion:
  \[
    \begin{tikzcd}
      A \ar[r,"\sim"] \ar[d,two heads] & A' \ar[dr,two heads] \ar[r,two heads] &
      P \ar[d,two heads] \ar[r] \drpullback & X \ar[d,two heads] \\
      B \ar[rr,"\sim"'] & & B' \ar[r,two heads] & Y
    \end{tikzcd}
  \]
  Since local equivalences are invariant under weak equivalence, as is locality by \cref{thm:fl-we}, if we replace $A$ and $B$ by $A'$ and $B'$ the hypotheses still hold.
  Thus, we may assume that the maps $B\to Y$ and $A\to P$, hence also $A\to X$, are fibrations.

  Now consider the following $3\times 3$ square in $\E/P$:
  \[
    \begin{tikzcd}
      A \ar[r] \ar[d] & P\times_B A \ar[d] \ar[r] & P\times_{Y} X \ar[d] \\
      P \times_X A \ar[r] \ar[d] & P\times A \ar[d] \ar[r] & P\times X \ar[d]\\
      P\times_Y B \ar[r] & P\times B \ar[r] & P\times Y
    \end{tikzcd}
  \]
  Here the rows and columns are defined by the following pullbacks:
  \[
    \begin{tikzcd}
      P\times_Y B \ar[r] \ar[d,two heads] \drpullback &
      P\times B \ar[d,two heads] \ar[r] \drpullback & B \ar[d,two heads]\\
      P \ar[r] & P\times Y \ar[r] & Y
    \end{tikzcd}
    \qquad
    \begin{tikzcd}
      P\times_X A \ar[r] \ar[d,two heads] \drpullback &
      P\times A \ar[d,two heads] \ar[r] \drpullback & A \ar[d,two heads]\\
      P \ar[r] & P\times X \ar[r] & X
    \end{tikzcd}
  \]
  \[
    \begin{tikzcd}
      P\times_Y X \ar[r] \ar[d,two heads] \drpullback &
      P\times X \ar[d,two heads] \ar[r] \drpullback & X \ar[d,two heads]\\
      P \ar[r] & P\times Y \ar[r] & Y
    \end{tikzcd}
    \qquad
    \begin{tikzcd}
      P\times_B A \ar[r] \ar[d,two heads] \drpullback &
      P\times A \ar[d,two heads] \ar[r] \drpullback & A \ar[d,two heads]\\
      P \ar[r] & P\times B \ar[r] & B
    \end{tikzcd}
  \]
  Since the right-hand maps in these rectangles are fibrations, by \cref{thm:pb-fl-fle} and the assumptions we have that
  \begin{itemize}
  \item $P\times_Y X\to P$ is $(P^*(S\dia))\pb$-local in $\E/P$.
  \item $P\times_B A\to P$ is $(P^*(S\dia))\pb$-local in $\E/P$.
  \item $P\times_Y B\to P$ is a $(P^*(S\dia))\pb$-local equivalence.
  \item $P\times_X A \to P$ is a $(P^*(S\dia))\pb$-local equivalence.
  \end{itemize}
  However, we also have pullbacks
  \[
    \begin{tikzcd}
      A \ar[d,two heads] \ar[r] \drpullback &
      P\times_B A \ar[d,two heads] \ar[rr] \drpullback && A \ar[d,two heads] \\
      P \ar[r] & P\times_Y X \ar[r] & B\times_Y X \ar[r,equals] & P
    \end{tikzcd}
  \]
  \[
    \begin{tikzcd}
      A \ar[d,two heads] \ar[r] \drpullback &
      P\times_X A \ar[d,two heads] \ar[rr] \drpullback && A \ar[d,two heads] \\
      P \ar[r] & P\times_Y B \ar[r] & X\times_Y B \ar[r,equals] & P
    \end{tikzcd}
  \]
  exhibiting $A$ as the simultaneous fiber in $\E/P$ of the fibrations $P\times_B A \fib P\times_Y X$ and $P\times_X A \fib P\times_Y B$.
  The former is a fibration between $(P^*(S\dia))\pb$-local objects, hence a $(P^*(S\dia))\pb$-local fibration, so that $A$ is also $(P^*(S\dia))\pb$-local.
  And the latter is a fibration between $(P^*(S\dia))\pb$-acyclic objects;
  thus since $(P^*(S\dia))\pb$-localization is left exact by \cref{thm:abfj} (and the fact that it coincides with $((P^*S)\dia)\pb$-localization), its fiber is also $(P^*(S\dia))\pb$-acyclic.
  So $A$ is both $(P^*(S\dia))\pb$-local and $(P^*(S\dia))\pb$-acyclic, hence the map $A\to P$ is a weak equivalence in \E, i.e.\ the given square is a homotopy pullback.
\end{proof}

\begin{prop}\label{thm:locfib}
  Let $S$ be a set of fibrations between fibrant objects and $p:Z\fib Y$ a fibration; the following are equivalent.
  \begin{enumerate}
  \item $p:Z\fib Y$ is a fibration in the $(S\dia)\pb$-local model structure on \E.\label{item:lf1}
  \item $p:Z\fib Y$ is $(Y^*(S\dia))\pb$-local in $\E/Y$.\label{item:lf2}
  \end{enumerate}
\end{prop}
\begin{proof}
  Since the $(S\dia)\pb$-local model structure is right proper,~\cite[Proposition 3.4.8]{hirschhorn:modelcats} tells us that~\ref{item:lf1} is equivalent to the $(S\dia)\pb$-localization square:
  \begin{equation}
    \begin{tikzcd}
      Z \ar[d,two heads,"f"'] \ar[r] & \Zhat \ar[d,two heads,"\fhat"]\\
      Y \ar[r] & \Yhat
    \end{tikzcd}\label{eq:locfib}
  \end{equation}
  being a homotopy pullback.

  Thus, if we assume~\ref{item:lf1}, then since pullback takes $(\Yhat^*(S\dia))\pb$-local objects to $(Y^*(S\dia))\pb$-local ones by \cref{thm:pb-fl-fle}, and locality is preserved by weak equivalences of fibrant objects, it will suffice to show that $\Zhat$ is $(\Yhat^*(S\dia))\pb$-local in $\E/\Yhat$.
  Now by assumption $\Zhat$ and $\Yhat$ are $S\pb$-local objects, so $\fhat$ is a fibration in the $S\pb$-local model structure $\E_{S\pb}$, hence $\Zhat$ is a fibrant object of the slice model structure $\E_{S\pb}/\Yhat$.
  On the other hand, any morphism $f:A\to B$ in $(\Yhat^*(S\dia))\pb$ is a pullback (in \E) of a morphism in $S\dia$, hence lies in $(S\dia)\pb$.
  Thus in particular it is an $(S\dia)\pb$-local equivalence in $\E$, hence a weak equivalence in $\E_{S\pb}/\Yhat$, and so the induced map
  \( \ehom{\E/\Yhat}(B,\Zhat) \to \ehom{\E/\Yhat}(A,\Zhat) \)
  is a weak equivalence for any such $f$; hence~\ref{item:lf2} holds.

  Conversely, suppose~\ref{item:lf2}.
  Note that $\fhat$, being a fibration between $(S\dia)\pb$-local objects, is a fibration in the $(S\dia)\pb$-local model structure; so by what we just proved, it is $(\Yhat(S\dia))\pb$-local in $\E/\Yhat$.
  On the other hand, by \cref{thm:unit-connected}, $Y\to \Yhat$ is a $(\Yhat^*(S\dia))\pb$-local equivalence and $Z\to \Zhat$ is a $(\Zhat^*(S\dia))\pb$-local equivalence.
  So by \cref{thm:conn-modal-pb}, the square~\eqref{eq:locfib} is a homotopy pullback, i.e.~\ref{item:lf1} holds.
\end{proof}

\begin{thm}\label{thm:lexloc-ttmt}
  For any set $S$ of morphisms in a \ttmt \E such that $S$-localization is left exact, the localized model structure $\E_S$ is again a \ttmt.
\end{thm}
\begin{proof}
  Since $\E_S$ has the same underlying category and cofibrations as \E, it is a \slcc Grothendieck 1-topos with cofibrations being the monomorphisms.
  It is combinatorial and simplicial by the usual constructions (e.g.~\cite[Theorem 4.1.1]{hirschhorn:modelcats}), and right proper since the localization is left exact.

  Let $S'$ be the result of replacing each morphism of $S$ by a weakly equivalent fibration between fibrant objects.
  Then $\E_S = \E_{S'}$, and by \cref{thm:abfj} $\E_{S'} = \E_{(S'{}\dia)\pb}$.
  Let $\dL_{S'{}\dia}$ be the \local and \stratified \nfs from \cref{thm:flf-nfs} applied to $S'{}\dia$;
  then by \cref{thm:locfib}, the morphisms admitting $\dL_{S'{}\dia}$-structure are precisely the fibrations of $\E_S$.
\end{proof}


%% file: summary.tex
\section{Summary}
\label{sec:summary}

Putting together our conclusions from \cref{sec:ttmt,sec:injmodel,sec:lex-loc}, we have the following.

\begin{thm}\label{thm:iotop-ttmt}
  Every model topos~\cite{rezk:homotopy-toposes} is Quillen equivalent to a \ttmt.
  Therefore, every Grothendieck \io-topos~\cite{lurie:higher-topoi} can be presented by a \ttmt.
\end{thm}
\begin{proof}
  A model topos is, by definition, Quillen equivalent to a left exact localization of a projective model structure on simplicial presheaves.
  Thus it is also equivalent to the corresponding localization of the \emph{injective} model structure, which is a \ttmt by \cref{thm:enrpre-ttmt,thm:lexloc-ttmt}.
\end{proof}

\begin{thm}\label{thm:models}
  For any Grothendieck \io-topos \fE, there is a regular cardinal \la such that \fE can be presented by a Quillen model category that interprets Martin-L\"{o}f type theory with the following structure:
  \begin{enumerate}
  \item $\Sigma$-types, a unit type, $\Pi$-types with function extensionality, identity types, and binary sum types.\label{item:sigpiid}
  \item The empty type, the natural numbers type, the circle type $S^1$, the sphere types $S^n$, and other specific ``cell complex'' types such as the torus $T^2$.\label{item:unparam-hits}
  \item As many universe types as there are inaccessible cardinals larger than $\la$, all closed under the type formers~\ref{item:sigpiid} and containing the types~\ref{item:unparam-hits}, and satisfying the univalence axiom.\label{item:univ}
  \item $\mathsf{W}$-types, pushout types, truncations, localizations, James constructions, and many other recursive higher inductive types.\label{item:hits}
  \end{enumerate}
\end{thm}
\begin{proof}
  By \cref{thm:iotop-ttmt}, $\fE$ can be presented by a \ttmt \E, to which we can apply \cref{thm:ttmt-models}.
\end{proof}

We also have the following (cf.~\cref{rmk:pfthy}):

\begin{prop}\label{thm:propresizing}
  Any \ttmt satisfies the propositional resizing principle for sufficiently large universes.
\end{prop}
\begin{proof}
  As stated in~\cite[Axiom 3.5.5]{hottbook},
  the propositional resizing principle\footnote{Not to be confused with the propositional resizing \emph{rule} of Voevodsky, which is not known to have any univalent models at all.} says that for universes $U^{\ka_1}$ and $U^{\ka_2}$ with $\ka_1<\ka_2$, the corresponding inclusion between universes of $(-1)$-truncated maps $\Prop^{\ka_1}\into \Prop^{\ka_2}$ is an equivalence.
  Since by univalence this is always a $(-1)$-truncated map, it suffices to construct a homotopy right inverse.
  For this it will suffice to show that there is a regular cardinal $\ka$ such that any $(-1)$-truncated fibration $f:X\to Y$ is equivalent over $Y$ to a relatively \ka-presentable one, since the latter will be classified by a map to $\Prop^\ka$.

  We will prove this first when \E is a left exact localization of a simplicial presheaf category, $\pr\C\S_{(S\dia)\pb}$.
  In this case it will suffice to show that any $(-1)$-truncated fibration $f:X\to Y$ is equivalent over $Y$ to a monomorphism, since by \cref{eg:pshf-relpres} a monomorphism is $\ka$-small for any sufficiently large $\ka$.
  Given such an $f$, for each $c\in \C$ let $Z(c)$ be the union of all connected components of $Y(c)$ that contain any points in the image of $X(c)$.
  The functorial action of $\C$ preserves these components, since $f$ is a natural transformation, so $Z$ becomes a presheaf as well and we have a factorization $X \xto{j} Z \xto{p} Y$ in which $Z\to Y$ is a monomorphism.

  Now each map $p_c : Z(c) \to Y(c)$ is also a Kan fibration, and each of its fibers is either empty or contractible.
  The same is true of $j_c : X(c)\to Y(c)$ since $f$ is $(-1)$-truncated, and by construction if a given fiber of $Z(c)$ is inhabited then so is the corresponding fiber of $X(c)$.
  Thus $j_c$ is a weak equivalence on each fiber, hence a weak equivalence, and so $j$ is a pointwise weak equivalence.

  It remains to show that $p:Z\to Y$ is a fibration in $\pr\C\S_{(S\dia)\pb}$.
  But the inclusion a union of connected components actually has the \emph{unique} right lifting property against all weak equivalences of simplicial sets.
  Thus $p$ is an injective fibration, since we can lift against any injective acyclic cofibration at each $c$ separately and fit the lifts together by uniqueness.
  Finally, $p$ is an $(S\dia)\pb$-local fibration since it is fiberwise equivalent to the $(S\dia)\pb$-local fibration $f$.

  This completes the proof when $\E=\pr\C\S_{(S\dia)\pb}$.
  Hence any Grothendieck \io-topos has a subobject classifier (this is also remarked after~\cite[Proposition 6.1.6.3]{lurie:higher-topoi}): a monomorphism of which every monomorphism is uniquely a (homotopy) pullback.
  In particular, the homotopy \io-category of any \ttmt \E has such a subobject classifier, which is presented by some $(-1)$-truncated fibration $\pi_{-1}:\widetilde{\Omega} \fib\Omega$ in \E.
  Let \ka be such that $\pi_{-1}$ is relatively \ka-presentable; then any $(-1)$-truncated morphism in \E will be a homotopy pullback of $\pi_{-1}$, hence equivalent to the corresponding strict pullback, which is relatively \ka-presentable.
\end{proof}

The main remaining open question is therefore whether the parametrized higher inductive types in \cref{thm:models}\ref{item:hits} can be constructed in such a way that the universes are closed under them.
Additionally,~\cite{ls:hits} does not construct all possible higher inductive types (only those without ``fibrant structure in the constructors''), and has not yet been generalized to indexed higher inductive types, inductive-inductive types, etc.


%% file: coherence.tex
\section{Coherence with universes}
\label{sec:coherence}

As noted in \cref{thm:ttmt-models}, the existing coherence theorems in the literature~\cite{klv:ssetmodel,lw:localuniv,awodey:natmodels} do not include an arbitrary family of universes.
Thus, here we sketch an extension of the coherence theorem of~\cite{lw:localuniv} to universe types.
For simplicity, we consider only non-cumulative Tarski universes.
The material in this appendix owes a great deal to conversations with Peter Lumsdaine.

First we review the original theorem of~\cite{lw:localuniv}, reformulated to match our \crefrange{sec:2cat}{sec:nfs} using the ideas of~\cite{awodey:natmodels}.
Recall that we write $\Ehat = \cPsFun(\E\op,\cGPD)$ to denote the 2-category of pseudofunctors $\E\op\to\cGPD$, and the notions of strict discrete fibration and representability from \cref{defn:dfib,defn:rep}.

\begin{defn}
  A \textbf{natural pseudo-model} is a category \E with a terminal object and a representable strict discrete fibration $\varpi:\dTm\to\dTy$ in \Ehat.
  It is a \textbf{natural model}~\cite{awodey:natmodels} if $\dTy$ (hence also \dTm) is discrete, i.e.\ a presheaf $\E\op\to\nSET$.
\end{defn}

It is shown in~\cite{awodey:natmodels} that natural models are equivalent to \emph{categories with families}~\cite{dybjer:internal-tt}.
Analogously, we have:

\begin{lem}
  A natural pseudo-model is the same as a comprehension category~\cite{jacobs:compr-cat} whose fibers are all groupoids.
\end{lem}
\begin{proof}[Sketch of proof]
  This is most obvious when \E has all pullbacks.
  In this case, by \cref{thm:univrep} a representable strict discrete fibration $\varpi:\dTm\to\dTy$ is classified by an essentially unique morphism $\dTy\to\cE$.
  Reformulating pseudofunctors to $\cGPD$ as fibrations with groupoidal fibers, this morphism becomes exactly the usual notion of a comprehension category.
  Unwinding this argument explicitly, we see that it also works even if not all pullbacks exist in \E.
\end{proof}

The analogous reformulation of a natural model is known as a \emph{category with attributes}~\cite{cartmell:gatcc,hofmann:ssdts}, while if the morphism $\dTy\to\cE$ of a comprehension category is a full inclusion it is called a \emph{display map category}~\cite{taylor:pracfdn}.

\begin{rmk}
  By \cref{thm:fib-eqv}, any representable morphism in \Ehat is equivalent to a natural pseudo-model, and thus might be called a \textbf{pseudonatural pseudo-model}.
  As in \cref{thm:univrep}, pseudonatural pseudo-models are only ``bicategorically'' equivalent to comprehension categories, and I don't know of any naturally-occurring examples of such that are not natural pseudo-models.
\end{rmk}

Natural models are the ``fully algebraic''\footnote{To make this completely precise, a natural model must actually be equipped with \emph{specified} representing objects for each strict pullback of $\varpi$ to a representable.} models of type theory.
The objects $\Gamma\in\E$ represent contexts, the objects $A\in \dTy(\Gamma)$ represent types in context $\Gamma$, the fiber over such an $A$ in $\dTm(\Gamma)$ represent terms $a:A$ in context $\Gamma$, and the representability of $\varpi$ yields context extension:
\[
  \begin{tikzcd}
    \E(-,\Gamma\ce A) \ar[r] \ar[d] \drpullback & \dTm \ar[d,"\varpi"] \\
    \E(-,\Gamma) \ar[r,"A"] & \dTy.
  \end{tikzcd}
\]
We write $\ec$ for the terminal object of \E, regarding it as the empty context.

A natural pseudo-model is the category-theoretic input from which a coherence theorem constructs a natural model.

\begin{eg}
  The \textbf{trivial} natural pseudo-model on \E is $\cEp \to \cE$.
\end{eg}

\begin{eg}\label{eg:cc-can}
  Any model category \E gives rise to a \textbf{canonical} natural pseudo-model where $\dTy=\Fib$ with $\dTy\into\cE$ the inclusion, so that $\dTm = \dFib \times_\cE \cEp$.
\end{eg}

\begin{eg}\label{eg:rep-cwf}
  For any morphism $\pi:\Util\to U$ in \E, the induced map $\E(-,\pi) : \E(-,\Util) \to \E(-,U)$ defines a \textbf{represented} natural model.
  This is essentially the construction of a model of type theory from a universe as in~\cite[\sect1.3]{klv:ssetmodel}.
\end{eg}

\begin{eg}\label{eg:las}
  If \E is a natural pseudo-model, its \textbf{left adjoint splitting}~\cite{lw:localuniv} is the natural model defined by
  \[ \las\dTy = \coprod_{V_A\in \E \atop E_A \in \dTy(V_A)} \E(-,V_A). \]
  There is a (surjective) morphism $\las\dTy\to \dTy$ composed of the classifying maps $E_A : \E(-,V_A) \to \dTy$, and we set $\las\dTm = \las\dTy \times_{\dTy} \dTm$.
  The object $V_A$ is known as the \emph{local universe} of $(V_A,E_A,A) \in \las\dTy(\Gamma)$.
\end{eg}

\begin{rmk}
  Any \nfs (\cref{defn:fcos}) is of course a natural pseudo-model with $\dTy=\cE$.
  On the other hand, if the comprehension morphism $\dTy\to\cE$ of a natural pseudo-model is faithful (such as for a natural model or display map category), then it is equivalent to a strict discrete fibration, which is a \nfs if it has small fibers.
  A ``$\dTy$-structure'' on a morphism is ``a way to express it as the comprehension of a dependent type''.
\end{rmk}

Now we sketch the interpretation of the basic type forming operations.
For any natural pseudo-model, define $\dTy^\Pi,\dTy^\Sigma,\dTy^{\Idtype},\dTy^{+} \in\Ehat$ as follows:
\begin{itemize}
\item $\dTy^\Pi(\Gamma)$ is the groupoid of pairs $(A,B)$ with $A\in \dTy(\Gamma)$ and $B\in \dTy(\Gamma\ce A)$.
\item $\dTy^\Sigma = \dTy^\Pi$.
\item $\dTy^{\Idtype}(\Gamma)$ is the groupoid of triples $(A,x,y)$ with $A\in \dTy(\Gamma)$ and $x,y\in \dTm(\Gamma)$ with $\varpi(x)=\varpi(y)=A$. 
\item $\dTy^{+} = \dTy\times\dTy$, so $\dTy^+(\Gamma)$ is the groupoid of pairs $(A,B)$ with $A,B\in\dTy(\Gamma)$.
\end{itemize}
The following notions are mostly as in~\cite[Definition 3.4.2.8]{lw:localuniv}:
\begin{itemize}
\item A \textbf{pseudo-stable class of $\Pi$-types} in a natural pseudo-model is a morphism $\Pi:\dTy^\Pi \to \dTy$ (corresponding to the formation rule) together with appropriate extra structure on each type $\Pi(A,B)$ (corresponding to the introduction, elimination, and equality rules).
\item A \textbf{strictly stable class of $\Pi$-types} is a pseudo-stable class for which \dTy is a natural model (since then the morphism $\Pi:\dTy^\Pi \to \dTy$ must be strictly natural, and so on).
\item A \textbf{weakly stable class of $\Pi$-types} is a span $\dTy^\Pi \ot \dG \to \dTy$ in which the first leg $\dG \to \dTy^\Pi$ is a surjective strict fibration (\cref{defn:dfib}), together with appropriate extra structure on the images of the second leg $\dG\to\dTy$.
  This is closely related to~\cite[Definition 3.4.2.5]{lw:localuniv} but more functorial\footnote{I am indebted to Peter Lumsdaine for suggesting this rephrasing.}; the objects of $\dG$ are the ``good $\Pi$-types''.
\end{itemize}
Similarly, we define all three kinds of $\Sigma$-types, identity types, and binary sum types.
Formal definitions of the ``appropriate extra structure'' can be found in~\cite{lw:localuniv}, and extensions to higher inductive types can be found in~\cite{ls:hits}; here we are recalling only the parts of the construction that are relevant for our extension to universes.

\begin{eg}
  The canonical natural pseudo-model of a \ttmt has pseudo-stable $\Pi$-types and $\Sigma$-types.
  It has only weakly stable $\Idtype$-types, but even in this case the type-\emph{forming} operation can be taken to be pseudo-stable (powers with $\Delta[1]$ as in \cref{sec:ttmt}).
  The translation from category-theoretic structure to type-theoretic structure in these cases can be found in~\cite[\sect1.4 and Proposition 2.3.3]{klv:ssetmodel} (for simplicial sets), \cite[Theorem 3.1]{aw:htpy-idtype}, \cite[Theorem 2.17]{warren:thesis}, \cite[Theorem 26]{ak:htmtt}, \cite[Theorem 4.2.2]{lw:localuniv}, and~\cite[\sect3]{awodey:natmodels}.

  However, binary sums are an example where weak stability matters even for the formation rule: we take the objects of $\dG(Y)$ to consist of three fibrations $A\fib Y$, $B\fib Y$, and $C\fib Y$ together with an acyclic cofibration $A+B \acof C$ over $Y$.
  By~\cite[Theorem 3.3]{ls:hits}, such data is stable under pullback, hence defines an object $\dG\in\Ehat$.
  The projection $\dG \to \dTy^{+} = \dFib\times \dFib$ picks out $A$ and $B$; this is a surjective strict fibration since fibrant replacements exist and are stable under isomorphism.
  The other projection $\dG \to \dTy= \dFib$ picks out $C$, with the rest of the structure constructed as in~\cite[Theorem 3.3]{ls:hits}.
\end{eg}

The coherence theorem of~\cite{lw:localuniv} says that if \E is locally cartesian closed (technically, a slightly weaker condition suffices), then weakly stable structure on $\dTy$ induces strictly stable structure on its left adjoint splitting $\las\dTy$ (\cref{eg:las}).
The basic observation is that
\[ \las\dTy^\Pi(\Gamma) \cong \coprod_{V_A\in \E \atop E_A \in \dTy(V_A)} \coprod_{V_B\in \E \atop E_B \in \dTy(V_B)} \coprod_{A:\Gamma\to V_A} \E(\Gamma\ce A,V_B) \]
and that when \E is locally cartesian closed, there is a universal object $V_A \lu V_B$\footnote{As remarked in~\cite{lw:localuniv}, this notation is abusive in that $V_A\lu V_B$ depends on $E_A$ too, not just on $V_A$ and $V_B$.} such that $\coprod_{A:\Gamma\to V_A} \E(\Gamma\ce A,V_B)$ is naturally isomorphic to $\E(\Gamma,V_A\lu V_B)$.
Thus $\las\dTy^\Pi$ is again a coproduct of representables:
\begin{equation}
\las\dTy^\Pi \cong \coprod_{V_A\in \E \atop E_A \in \dTy(V_A)} \coprod_{V_B\in \E \atop E_B \in \dTy(V_B)} \E(-,V_A \lu V_B).\label{eq:lastypi}
\end{equation}
Now if $\dTy$ has weakly stable $\Pi$-types, then since representables are projective in $\Ehat$ we can lift the map $\las\dTy^\Pi \to \dTy^\Pi$ along the surjective fibration $\dG\to\dTy^\Pi$ to the dashed map below:
\begin{equation}
  \begin{tikzcd}
    \las\dTy^\Pi \ar[d] \ar[dr,dashed] \ar[rr, dotted] & \ar[dr,phantom,"\cong"] & \las\dTy \ar[d] \\
    \dTy^\Pi & \dG \ar[r] \ar[l,two heads] & \dTy.
  \end{tikzcd}\label{eq:las-lift}
\end{equation}
The composite map $\las\dTy^\Pi \to \dG \to \dTy$ assigns to each $(V_A,E_A,V_B,E_B)$ a type $P \in \dTy(V_A\lu V_B)$, so that $(V_A,E_A,V_B,E_B) \mapsto (V_A \lu V_B, P)$ defines the dotted lifting making the quadrilateral commute up to isomorphism.
The rest of the structure is treated similarly, as are $\Sigma$-types, $\Idtype$-types, and so on.

Now we define universes in natural models and pseudo-models.

\begin{defn}
  A \textbf{level structure} is a partially ordered set \cL with binary joins and a partial endofunction $\suc : \cL \rightharpoonup \cL$ (not required to respect the ordering).
\end{defn}

\begin{eg}
  Any ordinal $\mu$ is a level structure, with $\suc(\al) = \al+1$ whenever $\al+1<\mu$ (thus $\suc$ is total if $\mu$ is a limit ordinal).
  Important special cases are $\mu=0$, $1$, $\om$, and $\om+1$.
\end{eg}

\begin{defn}\label{defn:tu}
  Let $\varpi:\dTm\to\dTy$ be a natural pseudo-model on \E, and \cL a level structure.
  An \textbf{\cL-family of pseudo Tarski universes} consists of:
  \begin{enumerate}
  \item For each $\al\in\cL$, types $\U_\al \in \dTy(\ec)$ and $\El_\al \in \dTy(\ec\ce \U_\al)$.\label{item:tu1}
  \item An extension of the function $\al \mapsto (\pi_\al : \ec\ce \U_\al \ce \El_\al \to \ec\ce \U_\al)$ to a functor from \cL to the category whose objects are morphisms in \E and whose morphisms are pullback squares in \E, sending each $\al\le\be$ to a pullback square\label{item:tu2}
    \[
      \begin{tikzcd}
        \ec\ce \U_\al \ce \El_\al \ar[d,"{\pi_\al}"'] \ar[r]\drpullback &
        \ec\ce \U_\be \ce \El_\be \ar[d,"{\pi_\be}"] \\
        \ec\ce \U_\al \ar[r,"{\Lift}"'] & \ec\ce \U_\be.
      \end{tikzcd}
    \]
  \item Whenever $\suc(\al)$ is defined, a morphism $\iu_\al : \ec \to \ec\ce \U_{\suc(\al)}$ and an isomorphism $\iu_\al^*(\El_{\suc(\al)}) \cong \U_\al$ in $\dTy(\ec)$.\label{item:tu3}
  \end{enumerate}
  We call it \textbf{strict} if $\varpi$ is a natural model (hence $\iu_\al^*(\El_{\suc(\al)}) = \U_\al$).
\end{defn}

Inside of type theory,~\ref{item:tu1} says that each universe is a type $\U_\al$ equipped with a coercion $(X:\U_\al) \vdash (\El_\al(X) \type)$ allowing us to regard its elements as types.
(The presence of the coercion $\El_\al$ is what makes these ``Tarski'' universes, in constrast to ``Russell'' universes whose elements are \emph{themselves} types; see \cref{rmk:russell}.)

Condition~\ref{item:tu2} says that any type can be lifted to a larger universe, $(X:\U_\al) \vdash (\Lift(X) : U_\be)$, which is isomorphic to it via functions $\lift : \El_\al(X) \to \El_\be(\Lift(X))$ and $\low : \El_\be(\Lift(X)) \to \El_\al(X)$ such that $\lift\circ\low \jdeq \iid$ and $\low\circ\lift \jdeq \iid$ judgmentally.
Finally, condition~\ref{item:tu3} says (in the strict case) that the universe $\U_\al$ is an \emph{element} of its successor universe: we have $\iu_\al : \U_{\suc(\al)}$ such that $\El_{\suc(\al)}(\iu_\al) \jdeq \U_\al$.

\begin{eg}\label{eg:ttmt-tu}
  Let \E be a \ttmt, let \la satisfy \cref{thm:uf-fibrant} for \E, and let \cL be a set of regular cardinals \al such that $\al\shgt\la$, with $\pi_\al : \Util_\al\fib U_\al$ the fibration constructed for \al in \cref{thm:uf-fibrant}.
  The induced well-ordering on \cL makes it a level structure, with $\suc(\al)$ the least element of \cL such that $U_\al$ is $\suc(\al)$-presentable (if such exists).

  We will extend $\{\pi_\al\}$ to an \cL-family of pseudo Tarski universes for the canonical natural pseudo-model of \E.
  Let $\ka\in\cL$ and assume inductively that we have constructed pullback squares as in \cref{defn:tu}\ref{item:tu2} for all $\al<\be<\ka$, which are additionally \F-morphisms.
  Then since $\Fka$ is \local, the map
  \[\colim_{\be<\ka} \Util_\be \to \colim_{\be<\ka} U_\be
  \]
  is a relatively \ka-presentable \F-algebra which pulls back to $\Util_\be$ over each $U_\be$.
  Thus it is classified by some map to $U_\ka$, i.e.\ we have an \F-morphism
  \[
    \begin{tikzcd}
      \colim_{\be<\ka}\Util_\be \ar[d,"{\pi_\al}"'] \ar[r]\drpullback &
      \Util_\ka \ar[d,"{\pi_\be}"] \\
      \colim_{\be<\ka}U_\be \ar[r,"{\Lift}"'] & U_\ka
    \end{tikzcd}
  \]
  inducing \F-morphisms $\Lift : U_\be \to U_\ka$ for all $\be<\ka$.
  Finally, since $U_\al$ is $\suc(\al)$-presentable by assumption whenever $\suc(\al)$ is defined, it is classified by some map $1 \to U_{\suc(\al)}$ which we can take as $\iu_\al$.
\end{eg}

In general we expect universes to be closed under the type operations.
Recall from \cref{eg:rep-cwf} that each $\pi_\al$ induces a represented natural model; thus we can ask the latter to have strictly stable $\Pi$-types, $\Sigma$-types, and so on.
In fact, we have $\E(-,U)^\Pi \cong \E(-,U\lu U)$ and similarly for the domains of all the other type operations.
Thus $\Pi$-types on $\E(-,\pi)$ are given by a morphism $\Pi : U\lu U \to U$ in \E, and so on.
We also expect $\El$ to respect these operations.

This corresponds type-theoretically to an operation
\[(X:\U_\al), (Y:\El_\al(X) \to \U_\al) \vdash (\Pi_\al(X,Y): \U_\al)\]
such that (in the strict case) $\El_\al(\Pi_\al(X,Y)) \jdeq \prod_{x:\El_\al(X)} \El_\al(Y(x))$.
But especially since our universes are non-cumulative, it is more useful to have relative operations
\[(X:\U_\al), (Y:\El_\al(X) \to \U_\be) \vdash (\Pi_{\al,\be}(X,Y): \U_{\al\vee\be})\]
in which the base and fibers can lie in different universes.
This suggests the following.

\begin{defn}
  Let $\varpi:\dTm\to\dTy$ be a natural pseudo-model on \E that has weakly stable $\Pi$-types.
  An \cL-family of pseudo Tarski universes is \textbf{closed under $\Pi$-types} if for each $\al,\be\in\cL$ we have dashed and dotted morphisms making the following diagram commute (with an isomorphism in the quadrilateral):
  \begin{equation}
    \begin{tikzcd}
      \E(-,(\ec\ce\U_\al)\lu (\ec\ce\U_\be)) \ar[d] \ar[dr,dashed] \ar[rr, dotted] &
      \ar[dr,phantom,"\cong"] & \E(-,\ec\ce\U_{\al\vee\be}) \ar[d] \\
      \dTy^\Pi & \dG \ar[r] \ar[l,two heads] & \dTy.
    \end{tikzcd}\label{eq:tu-lift}
  \end{equation}
  We say it is \textbf{strictly closed} if $\varpi$ is a natural model (so that in particular the isomorphism is an identity).
\end{defn}

Of course, this is very similar to the construction of split structure on $\las\dTy$ in~\eqref{eq:las-lift}.
Note that, as there, the dashed lift can always be chosen since $\dG\to\dTy^\Pi$ is surjective and representables are projective.

We can similarly define closure under $\Sigma$-types, identity types, binary sum types, and so on.
Note that we do not have to say anything about the introduction or elimination rules, as these only happen after the coercions $\El$ are applied.

\begin{eg}
  The discussion in \cref{sec:ttmt} implies that the families of universes from \cref{eg:ttmt-tu} for the canonical natural pseudo-model of a \ttmt \E are always closed under $\Sigma$-types, can be chosen to be closed under identity types and binary sum types by taking $\la$ sufficiently large, and are closed under $\Pi$-types if each $\al\in\cL$ is inaccessible.
  
  Consider for instance the case of $\Pi$-types: we choose a dashed lift in~\eqref{eq:tu-lift}, and then the composite $\E(-,U_\al \lu U_\be) \to \dG\to\dTy$ classifies the universal dependent product of a relatively \be-presentable fibration along a relatively \al-presentable one.
  Since $\al\vee\be$ is inaccessible, this is relatively $(\al\vee\be)$-presentable, hence is classified by some map $U_\al \lu U_\be \to U_{\al\vee\be}$.
\end{eg}

Now suppose $\varpi:\dTm\to\dTy$ is a natural pseudo-model on \E with weakly stable $\Pi$-types and a family of universes closed under them; we would like $\las\dTy$ to inherit a family of universes \emph{strictly} closed under $\Pi$-types.
However, although the map $\E(-,(\ec\ce\U_\al)\lu (\ec\ce\U_\be)) \to \dTy^\Pi$ in~\eqref{eq:tu-lift} factors through $\las\dTy^\Pi$ by picking out the summand of~\eqref{eq:lastypi} corresponding to $(V_A,E_A,V_B,E_B)=(\ec\ce\U_\al, \ec\ce\U_\al\ce\El_\al, \ec\ce\U_\be, \ec\ce\U_\be\ce\El_\be)$, the square we might hope for does not commute:
\begin{equation*}
  \begin{tikzcd}
    \E(-,(\ec\ce\U_\al)\lu (\ec\ce\U_\be)) \ar[d] \ar[rr] \ar[drr,phantom,"\scriptstyle\text{(does not commute)}"] &&
    \E(-,\ec\ce\U_{\al\vee\be}) \ar[d]\\
    \las\dTy^\Pi \ar[d] \ar[dr] \ar[rr] & \ar[dr,phantom,"\cong"] & \las\dTy \ar[d] \\
    \dTy^\Pi & \dG \ar[r] \ar[l,two heads] & \dTy.
  \end{tikzcd}\label{eq:las-tu-lift}
\end{equation*}
The left-bottom composite is the inclusion into a summand of $\las\dTy$ with $V_A = (\ec\ce\U_\al)\lu (\ec\ce\U_\be)$, while the top-right composite maps into a summand with $V_A = \ec\ce\U_{\al\vee\be}$.
Thus we need to generalize the local universes construction.

\begin{defn}
  A \textbf{family of local universes} for a natural pseudo-model $\varpi : \dTm\to \dTy$ on \E is an indexed family $\cV$ of pairs $(V,E)$ with $V\in \E$ and $E\in \dTy(V)$.
  If $\varpi$ has weakly stable $\Pi$-types
  \( \dTy^\Pi \ot \dG \xto{\Pi} \dTy \),
  then \cV \textbf{supports $\Pi$-types} if it is equipped with, for every $(V_A,E_A)$ and $(V_B,E_B)$ in \cV, some $G$ in the fiber of $\dG(V_A \lu V_B)$ over the universal pair of types in $\dTy^\Pi(V_A\lu V_B)$ (this is $(E_A[\pi_A],E_B[\pi_B])$ in the notation of~\cite[Lemma 3.4.2.4]{lw:localuniv}) and a pullback square
  \begin{equation}
    \begin{tikzcd}
      (V_A\lu V_B)\ce \Pi(G) \ar[d] \ar[r] \drpullback & V_{\Pi(A,B)}\ce E_{\Pi(A,B)} \ar[d] \\
      V_A\lu V_B \ar[r] & V_{\Pi(A,B)}
    \end{tikzcd}\label{eq:lu-pi-pb}
  \end{equation}
  for some $(V_{\Pi(A,B)},E_{\Pi(A,B)})\in\cV$.
\end{defn}

\begin{eg}\label{eg:lu-triv}
  The \textbf{trivial} family of local universes is the family of \emph{all} such pairs $(V,E)$, with $V_{\Pi(A,B)} = V_A \lu V_B$ and~\eqref{eq:lu-pi-pb} the identity square (for some arbitrary choice of $G$).
\end{eg}

\begin{eg}\label{eg:lu-tu}
  An \cL-family of universes (\cref{defn:tu}) induces a family of local universes with \cV the \cL-indexed family of pairs $(\ec\ce\U_\al,\El_\al)$ for $\al\in\cL$.
  If the family of universes is closed under $\Pi$-types, then \cV supports $\Pi$-types with $V_{\Pi(\al,\be)} = \ec\ce\U_{\al\vee\be}$.
\end{eg}

\begin{eg}\label{eg:lu-mixed}
  We can also consider the disjoint union of the families of local universes from \cref{eg:lu-triv,eg:lu-tu}.
  (Thus the pairs $(\ec\ce\U_\al,\El_\al)$ appear twice in the indexed family \cV, once as an arbitrary pair $(V,E)$ and once as a ``universe'' pair.)
  We define the $\Pi(A,B)$ operation to restrict to those of \cref{eg:lu-triv,eg:lu-tu} when $A$ and $B$ are of the same kind, and when they are of different kinds we ``forget'' that the universe pair $(\ec\ce\U_\al,\El_\al)$ is special and proceed as in \cref{eg:lu-triv}.
\end{eg}

Now given a family of local universes \cV, we define 
\[ \lasv\dTy = \coprod_{(V,E)\in \cV} \E(-,V) \]
with $\lasv\dTm = \lasv\dTy \times_{\dTy} \dTm$ as in \cref{eg:las}.
We have
\[ \lasv\dTy^\Pi = \coprod_{(V_A, E_A) \in\cV} \coprod_{(V_B, E_B) \in\cV} \E(-,V_A \lu V_B),\]
so if \E has weakly stable $\Pi$-types and \cV supports $\Pi$-types, the squares~\eqref{eq:lu-pi-pb} are precisely what is needed to induce the top dotted map:
\begin{equation}
  \begin{tikzcd}
    \lasv\dTy^\Pi \ar[d] \ar[dr,dashed] \ar[rr, dotted] & \ar[dr,phantom,"\cong"] & \lasv\dTy \ar[d] \\
    \dTy^\Pi & \dG \ar[r] \ar[l,two heads] & \dTy.
  \end{tikzcd}\label{eq:lasv-lift}
\end{equation}
Thus $\lasv\dTy$ again has strictly stable $\Pi$-types.
(We do need to check that the strictification of the introduction and elimination rules also works.
But here we can use the ordinary local universes construction, with the reindexing of $E_{\Pi(A,B)}$ to $V_A \lu V_B$ along~\eqref{eq:lu-pi-pb} replacing its isomorph $\Pi(G)$.)

Applied to \cref{eg:lu-triv} this reproduces the original local universes model $\las\dTy$.
On the other hand, applied to \cref{eg:lu-tu} in the case $\cL=1$, it reproduces the ``global universe'' coherence theorem of~\cite{klv:ssetmodel}.
More generally, when we apply it to \cref{eg:lu-tu,eg:lu-mixed}, we obtain exactly the commutative squares we needed:
\begin{equation*}
  \begin{tikzcd}
    \E(-,(\ec\ce\U_\al)\lu (\ec\ce\U_\be)) \ar[d] \ar[rr] \ar[drr,phantom,"\scriptstyle\text{(commutes!)}"] &&
    \E(-,\ec\ce\U_{\al\vee\be}) \ar[d]\\
    \lasv\dTy^\Pi \ar[d] \ar[dr] \ar[rr] & \ar[dr,phantom,"\cong"] & \lasv\dTy \ar[d] \\
    \dTy^\Pi & \dG \ar[r] \ar[l,two heads] & \dTy.
  \end{tikzcd}\label{eq:las-tu-lift}
\end{equation*}
We are thus most of the way to the following.

\begin{thm}\label{thm:lasv}
  Suppose $\varpi:\dTm\to\dTy$ is a natural pseudo-model on \E with weakly stable $\Pi$-types, \cL a level structure, and we are given an \cL-family of pseudo Tarski universes for $\dTy$ that is closed under $\Pi$-types.
  Then for \cV as in \cref{eg:lu-mixed}, the natural model $\lasv\varpi:\lasv\dTm\to\lasv\dTy$ has strictly stable $\Pi$-types and an \cL-family of strict Tarski universes strictly closed under $\Pi$-types.
\end{thm}
\begin{proof}
  We have not yet actually constructed an \cL-family of universes in $\lasv\dTy$ itself.
  Given $\al$, if $\suc(\al)$ is defined, we take $\lasU{\al}\in \lasv\dTy(\ec)$ to be $(\ec\ce\U_{\suc(\al)}, \El_{\suc(\al)}, \iu_\al)$, i.e.\ as the ``delayed substitution'' corresponding to the given isomorphism $\iu_\al^*(\El_{\suc(\al)}) \cong \U_\al$.
  Then the projection $\lasv\dTy \to \dTy$ maps $\lasU{\al}$ to $\U_\al$, so we can take $\lasEl{\al} = \El_\al$.
  On the other hand, if $\suc(\al)$ is not defined, we take $\lasU{\al}$ to be $(\ec,\U_{\al},\id_\ec)$.

  Since the comprehension of $\lasv\dTy$ factors through $\dTy$, the lifting isomorphisms of the latter immediately induce ones for the former, and similarly for the morphisms $\lasu{\al}$.
  And the definition of substitution in $\lasv\dTy$ implies that $\lasu{\al}^*(\lasEl{\suc(\al)}) = \lasU{\al}$ when $\suc(\al)$ is defined.
  Thus we have a family of strict Tarski universes for $\lasv\dTy$, and the observations above show that it is strictly closed under $\Pi$-types.
\end{proof}

We can deal similarly with $\Sigma$-types, identity types, and binary sum types.
Thus we have the following more precise version of \cref{thm:ttmt-models}:

\begin{thm}\label{thm:ttmt-lasv}
  For any \ttmt \E, there is a regular cardinal \la such that if the inaccessible cardinals greater than $\la$ have order type $\mu$, for some ordinal $\mu$, then the canonical natural pseudo-model of \E can be replaced by a natural model $\lasv\Fib$ with:
  \begin{enumerate}
  \item Strictly stable $\Sigma$-types, a unit type, $\Pi$-types with function extensionality, identity types, and binary sum types.\label{item:last1}
  \item A strictly stable empty type, natural numbers type, circle type $S^1$, sphere types $S^n$, and other specific ``cell complex'' types such as the torus $T^2$.\label{item:last2}
  \item A $\mu$-family of strict Tarski universes, strictly closed under the type formers~\ref{item:last1} and containing the types~\ref{item:last2}, and satisfying the univalence axiom.\label{item:last3}
  \item Strictly stable $\mathsf{W}$-types, pushout types, truncations, localizations, James constructions, and many other recursive higher inductive types.\label{item:last4}
  \end{enumerate}
\end{thm}
\begin{proof}
  As in \cref{thm:ttmt-models}, we take \la large enough that the types in~\ref{item:last2} are $\la$-presentable.
  Then \cref{thm:lasv,eg:ttmt-tu} yield~\ref{item:last1} and~\ref{item:last3}.
  Finally, the local universes coherence arguments for higher inductive types from~\cite{ls:hits} apply just as well to the modified version $\lasv\dTy$, yielding~\ref{item:last2} and~\ref{item:last4}.
\end{proof}

\begin{rmk}\label{rmk:russell}
  If the successor in \cL is totally defined (e.g.\ if $\mu$ in \cref{thm:ttmt-lasv} is a limit ordinal, such as $\om$), then in \cref{thm:lasv} we can use \cref{eg:lu-tu} instead of \cref{eg:lu-mixed}.
  This gives a model like that of~\cite{coquand:pshf-model} in which every type belongs to a unique universe, which is the closest possible approximation to (non-cumulative) Russell-type universes obtainable in a natural model.
  Probably this construction could also be adapted to a semantic structure that models Russell-type universes more closely, such as the generalized algebraic theories of~\cite{coquand:cannorm-dtt,sterling:att-universes,chs:hocan-cubical}.
\end{rmk}

\begin{rmk}
  The universes we have constructed here are \emph{non-cumulative} in that $\El(\Lift(A))$ is judgmentally isomorphic to $\El(A)$ rather than judgmentally equal to it.
  Such universes are implemented in Agda and Lean, and are easy to use as long as we formulate the type formers appropriately using $\al\vee\be$.
  (Note that the universes in Agda have order type $\om+1$, with the $\om^{\mathrm{th}}$ universe not belonging to any larger universe; this is covered by \cref{thm:ttmt-lasv} as long as $\mu\ge \om+1$.)
  Other proof assistants such as Coq use cumulative universes instead, as does~\cite{hottbook}; it should be possible to model these in a \ttmt as well (some ideas are sketched in~\cite{shulman:invdia,shulman:elreedy}), but we leave this for future work.
\end{rmk}
